\documentclass[a4paper,reqno,11pt]{amsart}
\usepackage[T1]{fontenc}
\usepackage[utf8]{inputenc}
\usepackage{amssymb, amsmath, amsthm, graphicx, color}
\usepackage[inline]{enumitem}
\usepackage[usenames,dvipsnames,svgnames,table]{xcolor}
\usepackage[foot]{amsaddr}
\usepackage{comment}

\usepackage{booktabs}
\usepackage{multirow}

\usepackage[longnamesfirst,numbers,sort&compress]{natbib}

\usepackage{comment}

\makeatletter
\def\thm@space@setup{
\thm@preskip=4mm
\thm@postskip=0mm
}
\makeatother

\usepackage{hyperref}
\hypersetup{
    pdftitle={Weak coloring numbers of minor-closed graph classes},
    pdfauthor={Jędrzej Hodor, Hoang La, Piotr Micek, Clément Rambaud},
    colorlinks,
    linkcolor={RoyalBlue},
    citecolor={RubineRed},
    urlcolor={blue!80!black}
}

\usepackage{cleveref}

\usepackage{mathtools}

\usepackage{thmtools, thm-restate}

\setenumerate{label=\textup{(\roman*)}, noitemsep, topsep=0pt,
labelindent=.2em, leftmargin=*, widest=iii,}
\setitemize{noitemsep, topsep=-\parskip, labelindent=.2em, leftmargin=*, widest=iii,}

\newenvironment{enumerateOurAlph}{\begin{enumerate}
[label={\normalfont(\makebox[\mywidth]{\alph*})}]}{\end{enumerate}}

\newenvironment{enumerateOurAlphPrim}{\begin{enumerate}
[label={\normalfont(\makebox[\mywidthprim]{\alph*'})}]}{\end{enumerate}}

\newenvironment{enumerateOurAlphPrimPrim}{\begin{enumerate}
[label={\normalfont(\makebox[\mywidthprimprim]{\alph*''})}]}{\end{enumerate}}

\newenvironment{enumerateOurAlphCapital}{\begin{enumerate}
[label={\normalfont(\makebox[\mywidthA]{\Alph*})}]}{\end{enumerate}}

\newenvironment{enumerateOurAlphCapitalPrim}{\begin{enumerate}
[label={\normalfont(\makebox[\mywidthAprim]{\Alph*'})}]}{\end{enumerate}}

\newdimen\mywidth
\sbox0{m}%
\newdimen\mywidthprim
\sbox0{m}%
\newdimen\mywidthprimprim
\sbox0{m}%
\newdimen\mywidthA
\sbox0{m}%
\newdimen\mywidthAprim
\sbox0{m}%
\DeclarePairedDelimiter\set{\{}{\}}

\theoremstyle{plain}
\newtheorem{thm}{Theorem}
\newtheorem*{thm*}{Theorem}
\newtheorem{theorem}[thm]{Theorem}

\newtheorem{lemma}[thm]{Lemma}
\newtheorem*{lemma*}{Lemma}
\newtheorem{cor}[thm]{Corollary}
\newtheorem*{cor*}{Corollary}

\newtheorem*{corollary*}{Corollary}
\newtheorem{obs}[thm]{Observation}

\theoremstyle{remark}

\newtheorem{ques}{Question}

\newtheorem{claim1}{Claim}
\crefname{obs}{Observation}{Observations}
\theoremstyle{definition}

\newtheorem*{conj*}{Conjecture}
\crefname{lem}{Lemma}{Lemmas}
\crefname{thm}{Theorem}{Theorems}
\crefname{cor}{Corollary}{Corollaries}

\newenvironment{proofclaim}[1][]
	{\begin{proof}[Proof of the claim] }{\end{proof}}

\newcommand{\td}{\operatorname{td}}
\newcommand{\rtd}{\operatorname{rtd}}
\newcommand{\pw}{\operatorname{pw}}
\newcommand{\tw}{\operatorname{tw}}

\newcommand{\cgF}{\mathcal{F}}

\newcommand{\Oh}{\mathcal{O}}

\newcommand{\dist}{\mathrm{dist}}

\newcommand{\calC}{\mathcal{C}}

\let\leq\leqslant
\let\geq\geqslant

\let\subset\subseteq

\let\epsilon\varepsilon

\DeclareMathOperator\WReach{WReach}
\DeclareMathOperator\wcol{wcol}

\DeclareMathOperator\vc{vc}

\renewcommand{\setminus}{-}

\newcommand{\cc}{\bar{c}}

\title{Weak coloring numbers of minor-closed graph classes}

\begin{document}

\author[Hodor]{Jędrzej Hodor}
\address[J.~Hodor]{Theoretical Computer Science Department, 
Faculty of Mathematics and Computer Science and Doctoral School of Exact and Natural Sciences, Jagiellonian University, Krak\'ow, Poland}
\email{jedrzej.hodor@gmail.com}

\author[La]{Hoang La}
\address[H.~La]{LISN, Universit\'e Paris-Saclay, CNRS, Gif-sur-Yvette, France}
\email{hoang.la.research@gmail.com}

\author[Micek]{Piotr Micek}
\address[P.~Micek]{Theoretical Computer Science Department, 
Faculty of Mathematics and Computer Science, Jagiellonian University, Krak\'ow, Poland}
\email{piotr.micek@uj.edu.pl}

\author[Rambaud]{Clément Rambaud}
\address[C.~Rambaud]{Universit\'e C\^ote d'Azur, CNRS, Inria, I3S, Sophia-Antipolis, France}
\email{clement.rambaud@inria.fr}

\thanks{This research was funded by the National Science Center of Poland under grant UMO-2023/05/Y/ST6/00079 within the WEAVE-UNISONO program. J.\ Hodor was partially supported by a Polish Ministry of Education and Science grant (Perły Nauki; PN/01/0265/2022). C.\ Rambaud was partially supported by the French Agence Nationale de la Recherche under contract Digraphs ANR-19-CE48-0013-01}

\begin{abstract}
We study the growth rate of weak coloring numbers of graphs excluding a fixed graph as a minor. 
Van den Heuvel et al. (European J.\ of Combinatorics, 2017) showed that for a fixed graph $X$, the maximum $r$-th weak coloring number of $X$-minor-free graphs is polynomial in $r$.
We determine this polynomial up to a factor of $\Oh(r \log r)$.
Moreover, we tie the exponent of the polynomial to a structural property of $X$, namely, $2$-treedepth.
As a result, 
for a fixed graph $X$ and an $X$-minor-free graph $G$, we show that $\wcol_r(G)= \mathcal{O}(r^{\td(X)-1}\log r)$, which improves on the bound $\wcol_r(G) = \Oh(r^{g(\td(X))})$ given by Dujmović et al.\ (SODA, 2024), where $g$ is an exponential function. 
In the case of planar graphs of bounded treewidth, we show that the maximum $r$-th weak coloring number is in $\Oh(r^2\log r$), which is best possible.

\end{abstract}

\maketitle

\newpage
 
\section{Introduction}\label{sec:intro}

Let $G$ be a graph, let $\Pi(G)$ be the set of all vertex orderings of $G$, let $\sigma \in \Pi(G)$, and let $r$ be a nonnegative integer.
For all $u$ and $v$ vertices of $G$, we say that
$v$ is \emph{weakly $r$-reachable from $u$ in $(G,\sigma)$}, if there exists
a path between $u$ and $v$ in $G$
containing at most $r$ edges 
such that for every vertex $w$ on the path, $v\leq_{\sigma} w$. 
Let $\WReach_r[G, \sigma, u]$ be the set of vertices that are weakly $r$-reachable from $u$ in $(G,\sigma)$. 
The $r$-\emph{th weak coloring number} of $G$ is defined as
\[\wcol_r(G) = \min_{\sigma \in \Pi(G)}\max_{u \in V(G)}\ |\WReach_r[G, \sigma, u]|.\]
Let $X$ be a graph.
The \emph{treedepth} of $X$, denoted by $\td(X)$, 
is defined recursively as follows
\[
\td(X)=\begin{cases}
0 & \textrm{if $X$ is the null graph,}\\
\min_{v \in V(X)} \td(X-v) + 1 &\textrm{if $X$ is connected\footnotemark, and}\\
\max_{i \in [k]}\td(C_i)&\textrm{if $X$ consists of components $C_1,\dots,C_k$ and $k > 1$.}
\end{cases}
\]
\footnotetext{In this paper, connected graphs are nonnull, that is, they have at least one vertex. Note that a tree is defined as a connected forest, thus, trees and subtrees are also assumed to be nonnull.}

The following two theorems are among the main contributions of this paper.

\begin{thm}\label{thm:td}
For every positive integer $t$, for every graph $X$ with $\td(X) \leq t$, there exists an integer $c$ such that for every graph $G$, if $G$ is $X$-minor-free, then for every integer $r$ with $r \geq 2$,
\[
\wcol_r(G) \leq c \cdot r^{t-1} \log r.
\]    
\end{thm}

\begin{thm}\label{thm:td_tw}
For every integer $t$ with $t \geq 2$, for every graph $X$ with $\td(X) \leq t$, there exists an integer $c$ such that for every graph $G$, if $G$ is $X$-minor-free, then for every integer $r$ with $r \geq 2$,
\[
\wcol_r(G) \leq c \cdot (\tw(G)+1) \cdot r^{t-2} \log r.
\]    
\end{thm}

Weak coloring numbers were introduced by Kierstead and Yang~\cite{KY03} in 2003, though a parameter similar to $\wcol_2(G)$ is already present 
in the work of Chen and Schelp~\cite{CS93} from 1993.
This family of parameters gained considerable attention when Zhu~\cite{Zhu09} proved that it captures important and robust notions of sparsity, namely, bounded expansion and nowhere denseness. 
Specifically, a class of graphs $\calC$ has bounded expansion if and only if 
there exists a function $g$ such that for every graph $G$ in $\calC$ and every positive integer $r$, we have $\wcol_r(G) \leq g(r)$. 
Classes of bounded expansion include in particular, planar graphs, graphs of bounded treewidth, and proper minor-closed graph classes; 
see the book by Nešetřil and Ossona de Mendez~\cite{sparsity} or the recent lecture notes of Pilipczuk, Pilipczuk, and Siebertz~\cite{notes} for more information on this topic. 
Many algorithmic problems were solved using the weak coloring numbers characterization of sparse graphs.
Dvo\v{r}\'ak showed a constant-factor approximation for distance versions of domination number and independence number~\cite{Dvorak13},
with further applications in fixed-parameter algorithms and kernelization by Eickmeyer, Giannopoulou, Kreutzer, Kwon, Pilipczuk, Rabinovich, and Siebertz~\cite{EGKKPRS17}.
Grohe, Kreutzer, and Siebertz proved that deciding first-order properties is fixed-parameter tractable in nowhere dense graph classes~\cite{GroheKreutzerSiebertz17}. 
Reidl and Sullivan presented an algorithm counting the number of occurrences of a fixed induced subgraph in sparse graphs~\cite{ReidlSullivan23}. 
The time complexities of all these algorithms depend heavily on the asymptotics of $\wcol_r$ in respective classes of graphs.

The growth rate of $\wcol_r(G)$ when $G$ is in a fixed proper minor-closed class of graphs has been extensively studied.
In particular, 
Grohe, Kreutzer, Rabinovich, Siebertz, and Stavropoulos~\cite{Grohe15} proved  that
if $\tw(G)\leq t$\footnote{For a graph $G$, let $\tw(G)$, $\pw(G)$, and $\vc(G)$ stand for the treewidth, pathwidth, and vertex cover number of $G$ respectively.}, 
then $\wcol_r(G)\leq \binom{r+t}{t}$. 
This is tight as for all nonnegative 
integers $r,t$ they constructed a graph $G_{r,t}$ 
with $\tw(G_{r,t}) = t$ and $\wcol_r(G_{r,t}) = \binom{r+t}{t}$\footnote{We recall the construction in~\Cref{subsec:outine_the_key_parameter}.}.  
In the class of planar graphs, $\wcol_r(G) = \Oh(r^3)$ as proved by van den Heuvel, Ossona de Mendez, Quiroz, Rabinovich, and Siebertz~\cite{vdHetal17}.
On the other hand, for the family of stacked triangulations (i.e.\ planar graphs of treewidth at most $3$), we have $\wcol_r(G) = \Omega(r^2\log r)$, 
as shown by Joret and Micek~\cite{JM22}. 
The exact growth rate of maximum $r$-th weak coloring numbers of planar graphs is unknown.
\Cref{thm:td_tw} immediately implies that in the class of
planar graphs (or graphs of bounded Euler genus) of bounded treewidth, we have $\wcol_r(G) = \Oh(r^2\log r)$, which is tight. 
Indeed, it follows from Euler's formula that graphs of Euler genus at most $g$ exclude $K_{3,2g+3}$ as a minor 
and $\td(K_{3,2g+3}) = 4$ for all nonnegative integers $g$.

\begin{cor}\label{cor:genus}
For all nonnegative integers $g,w$, there exists an integer $c$ such that for every graph $G$ of Euler genus at most $g$ and with $\tw(G) \leq w$, and for every integer $r$ with $r \geq 2$,
\[
\wcol_r(G) \leq c \cdot r^{2} \log r.
\]    
\end{cor}

Since outerplanar graphs are $K_{2,3}$-minor-free and have bounded treewidth, \Cref{thm:td_tw} yields that $\wcol_r(G)=\Oh(r\log r)$ in the class of outerplanar graphs $G$. 
This was already proved by Joret and Micek~\cite{JM22}, who additionally showed that this bound is tight.

More generally, fix a graph $X$.
What is the growth rate with respect to $r$ of the maximum of $\wcol_r(G)$ for all $X$-minor-free graphs $G$?
Van den Heuvel et al.~\cite{vdHetal17} showed that $\wcol_r(G)=\Oh\left(r^{|V(X)|-1}\right)$.
Subsequently, van den Heuvel and Wood~\cite{vandenHeuvel2018} proved that $\wcol_r(G)=\Oh\left(r^{\vc(X)+1}\right)$.
Dujmović, Hickingbotham, Hodor, Joret, La, Micek, Morin, Rambaud, and Wood~\cite{DHHJLMMRW24} proved that there exists an exponential function $g$ such that $\wcol_r(G)=\Oh\left(r^{g(\td(X))}\right)$. 
We directly improve this result, namely, \Cref{thm:td} states that $\wcol_r(G) = \Oh\left(r^{\td(X)-1}\log r\right)$ and \Cref{thm:td_tw} states that $\wcol_r(G)=\Oh\left(\tw(G) \cdot r^{\td(X)-2} \log r\right)$.
Moreover, since $\td(X)-1 \leq \vc(G)$, we obtain $\wcol_r(G)=\Oh\left(r^{\vc(X)} \log r\right)$ and $\wcol_r(G)=\Oh\left(\tw(G) \cdot r^{\vc(X)-1} \log r\right)$.
In these cases, the construction of Grohe et al.~\cite{Grohe15} witnesses that our bounds are tight up to a factor of $\Oh(r\log r)$ in the general case and up to $\Oh(\log r)$ in the case of bounded treewidth.
Most of the known bounds on weak coloring numbers of minor-closed graph classes are summarized in \Cref{table:intro}.

All this previous work can be seen as an effort to understand the following graph parameter. 
For a given graph $X$, let
\begin{align*}
    f(X) = \inf\big\{\alpha \in \mathbb{R} \mid &\text{ there exists } c > 0 \text{ such that for every $X$-minor-free graph $G$}\\
    &\text{ and for every nonnegative integer $r$, } \wcol_r(G) \leq c \cdot r^\alpha\big\}.
\end{align*}
The question is whether $f$ is tied to%
\footnote{Two graph parameters $p,q$ are said to be \emph{tied} if there are two functions $\alpha,\beta$ such that $p(G) \leq \alpha(q(G))$ and $q(G) \leq \beta(p(G))$ for every graph $G$.}
some other well-established graph parameters.
Recall that for every graph $X$,
\[\tw(X) \leq \pw(X) \leq \td(X) - 1 \leq \vc(X) \leq |V(X)| - 1.\]
The aforementioned results imply that $\tw(X)-1 \leq f(X)\leq \td(X)-1$.
However, 
$f$ is not tied to any of these parameters. 
Indeed, neither pathwidth nor treedepth can lower bound $f$.
For every positive integer $k$, let $T_k$ be a complete ternary tree of vertex-height $k$.
Recall that there is a constant depending on $k$ bounding pathwidth of $T_k$-minor-free graphs by Robertson-Seymour Excluded Tree Minor Theorem~\cite{GM1}.
Also, it is easy to show\footnote{Proceed by induction on $\pw(G)$. We may assume that $G$ is connected. If $\pw(G)=0$, then $G$ has no edge and so $\wcol_r(G) \leq 1$. If $\pw(G)>0$, let $Q$ be a shortest path from the first bag to the last bag of an optimal path decomposition of $G$. Then $\pw(G-V(Q))<\pw(G)$ and so by induction $\wcol_r(G-V(Q)) \leq 1+(\pw(G)-1)(2r+1)$. Let $\sigma_0$ be an ordering of $V(G)$ witnessing this fact. Now, let $\sigma$ be an ordering of $V(G)$ extending $\sigma_0$ such that the vertices in $V(Q)$ appear first. By \Cref{lemma:intersection_ball_with_geodesic}, it follows that $\sigma$ witnesses $\wcol_r(G) \leq 1+(\pw(G)-1)(2r+1) + (2r+1) = 1+\pw(G)(2r+1)$.} that $\wcol_r(G)\leq 1 + \pw(G)(2r+1)$ for every 
graph $G$.
Thus, $f(T_k) \leq 1$ while $\pw(T_k) = k$ and $\td(T_k) = k+1$.
Next, we argue that neither treewidth nor pathwidth can upper-bound $f$.
For every positive integer $k$, let $L_k$ be a ladder with $k$ rungs. 
There is a graph $G_{r,t}$ (constructed in~\cite{Grohe15}) such that $\wcol_r(G_{r,t})=\Omega(r^t)$, and if $k = \Omega(\log t)$, then $G_{r,t}$ excludes $L_k$ as a minor.
Therefore, $f(L_k) = 2^{\Omega(k)}$, and $\tw(L_k) \leq \pw(L_k) \leq 2$.

Surprisingly, the key parameter to our problem is $2$-treedepth as defined by Huynh, Joret, Micek, Seweryn, and Wollan in~\cite{HJMSW22},
where they use it to characterize the structure of graphs excluding a fixed ladder as a minor. 
Let $X$ be a graph.
A {\em cut vertex} of \(X\) is a vertex \(v\in V(X)\) such that \(X - v\) has more components than \(X\).
A \emph{block} of \(X\) is a maximal connected subgraph of \(X\) without a cut vertex.\footnote{The blocks can be of three types: maximal $2$-connected subgraphs, cut edges together with
their endpoints, and isolated vertices.
Two blocks have at most one vertex in common, and such a vertex is always a cut vertex.}
The $2$-\emph{treedepth} of $X$, denoted by $\td_2(X)$, is defined recursively as follows
\[
\td_2(X)=\begin{cases}
0 & \textrm{if $X$ is the null graph,}\\
\min_{v \in V(X)} \td_2(X-v) + 1 &\textrm{if $X$ consists of one block, and}\\
\max_{i \in [k]}\td_2(B_i)&\textrm{if $X$ consists of blocks $B_1,\dots,B_k$ and $k > 1$.}
\end{cases}
\]
We show that $f$ is tied by a linear function to $\td_2$.
The first inequality in the theorem below is witnessed by the construction given in~\cite{Grohe15}.
\begin{thm}\label{thm:main_intro}
For every graph $X$ with at least one edge, we have
    \[\td_2(X)-2 \leq f(X) \leq 2\td_2(X) -3.\]
\end{thm}
    
To prove Theorems~\ref{thm:td}~and~\ref{thm:main_intro}, we prove that the value of $f$ is tied with the maximum $t$ such that $X$ is a subgraph of $G_{r,t}$ (as in \cite{Grohe15}).
In other words, we prove that $G_{r,t}$ is the obstruction for the growth of weak coloring numbers.
More precisely, we introduce a slightly modified version of $2$-treedepth, which we call \emph{rooted $2$-treedepth} and denote by $\rtd_2(\cdot)$.
Later, we show that for all graphs $X$ with at least one edge, $\rtd_2(X)$ is the minimum $t$ such that there exists $r$ such that $X$ is a subgraph of $G_{r,t-1}$.
See \Cref{sec:rtd2} for the definition of rooted $2$-treedepth and \Cref{lemma:rtd2_and_Grt} for the equivalence.

\begin{table}[!ht]
    \centering
    \setlength{\tabcolsep}{1.25ex}
    \def\arraystretch{1.5}
    \begin{tabular}{l@{\hspace{1.5ex}}c@{\hspace{1.5ex}}lc@{\hspace{1.5ex}}l}
        \toprule
        \bf Class $\mathcal{C}$ & \multicolumn{2}{c}{\bf lower bound} & \multicolumn{2}{c}{\bf upper bound} \\
        \midrule
        planar  & $\Omega(r^2\log r)$ & \cite{JM22} & $\Oh(r^3)$ & \cite{vdHetal17} \\
        \midrule
        planar and $\tw \leq k$ & $\Omega(r^2 \log r)$ & \cite{JM22} & $\Oh(r^2 \log r)$  &  \Cref{thm:td_tw} \\
        \midrule
        Euler genus $\leq g$ & $\Omega(r^2\log r)$ & \cite{JM22} & $\Oh(r^3)$ &   \cite{vdHetal17} \\
        \midrule
        Euler genus $\leq g$ and $\tw \leq k$ & $\Omega(r^2 \log r)$ & \cite{JM22} & $\Oh(r^2 \log r)$ & \Cref{thm:td_tw} \\
        \midrule
        outerplanar  & $\Omega(r \log r)$ & \cite{JM22} & $\Oh(r \log r)$ & \cite{JM22}  \\
        \midrule
        $K_{2,t}$-minor-free & $\Omega(r \log r)$ & \cite{JM22} & $\Oh(r \log r)$ & \Cref{thm:td_tw} \\
        \midrule
        $\tw \leq k$ & $\binom{r+k}{k}$ & \cite{Grohe15} & $\binom{r+k}{k}$ & \cite{Grohe15} \\
        \midrule
        $K_t$-minor-free & $\Omega(r^{t-2})$ & \cite{Grohe15} & $\Oh(r^{t-1})$ & \cite{vdHetal17} \\
        \midrule
        $K_{s,t}$-minor-free & $\Omega(r^{s-1} \log r)$ & \cite{JM22} & $\Oh(r^{s} \log r)$ & \Cref{thm:td} \\
        \midrule
        $K_{s,t}$-minor-free and $\tw \leq k$ & $\Omega(r^{s-1} \log r)$ & \cite{JM22} & $\Oh(r^{s-1} \log r)$ & \Cref{thm:td_tw} \\
        \midrule
        $X$-minor-free & $\Omega(r^{\rtd_2(X)-2})$ & \cite{Grohe15} & $\Oh(r^{\rtd_2(X)-1} \log r)$ &  \Cref{thm:rooted_main} \\
        \midrule
        $X$-minor-free and $\tw \leq k$ & $\Omega(r^{\rtd_2(X)-2})$ & \cite{Grohe15} & $\Oh(r^{\rtd_2(X)-2} \log r)$ & \Cref{thm:rooted_main_bounded_tw} \\
        \bottomrule
    \end{tabular}
    \medskip
    \caption{
    Lower and upper bounds on $\max_{G \in \mathcal{C}} \wcol_r(G)$ for some minor-closed graphs classes $\mathcal{C}$. 
    The variables $g,k,s,t$ are fixed positive integers with $s+3 \leq t \leq k$, and $X$ is a fixed nonnull graph.
    The weak coloring numbers of $K_{s,t}$-minor-free graphs were first studied by van den Heuvel and Wood in~\cite{vandenHeuvel2018}.
    In particular they proved the upper bound $\Oh(r^{s+1})$ and they conjectured $\Oh(r^s)$.
    \Cref{thm:td} implies that $K_{s,t}$-minor-free graphs have weak coloring numbers in $\Oh(r^{s}\log r)$.
    The lower bound $\Omega(r^{s-1} \log r)$ follows from the fact that 
    graphs of simple treewidth $s$ are $K_{s,t}$-minor-free and among them there 
    are graphs with weak coloring numbers in $\Omega(r^{s-1}\log r)$, see~\cite{JM22} for further details.    
    }
    \label{table:intro}
\end{table}

\begin{figure}[tp]
    \centering 
    \includegraphics[scale=1]{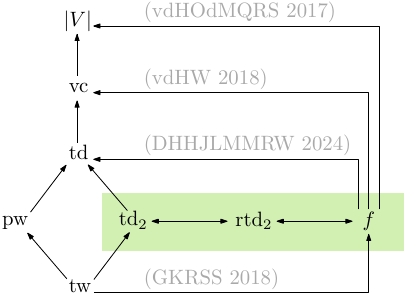} 
    \caption{Connections of $f$ to other graph parameters.
    An arrow from a parameter $p$ to a parameter $q$ indicates that there is a function $\alpha$ such that $p(X) \leq \alpha(q(X))$ for every graph $X$. We show that $f$ is tied to $\td_2$ and $\rtd_2$ but not to $\tw$, $\pw$, $\td$, $\vc$, or $|V|$.
    The results marked in the figure (top-to-bottom) are in \cite{vdHetal17}, \cite{vandenHeuvel2018}, \cite{DHHJLMMRW24}, and \cite{Grohe15} respectively.} \label{fig:results}
\end{figure}

Given a graph $X$ with at least one edge, we will show that $\td_2(X) \leq \rtd_2(X) \leq 2\td_2(X)-2$ and it will be clear from the definition that $\rtd_2(X) \leq \td(X)$.
As a consequence, Theorems~\ref{thm:td}~and~\ref{thm:main_intro} are implied by the following more accurate technical statement.
See also \Cref{fig:results}.

\begin{theorem}\label{thm:rooted_main}
For every positive integer $t$, for every graph $X$ with $\rtd_2(X) \leq t$, there exists an integer $c$ such that for every graph $G$, if $G$ is $X$-minor-free, then for every integer $r$ with $r \geq 2$,
\[
\wcol_r(G) \leq c \cdot r^{t-1} \log r.
\]    
\end{theorem}

Moreover, for all nonnegative integers $r,t$, the graph $G_{r,t}$ of~\cite{Grohe15} satisfies $\rtd_2(G_{r,t}) = t + 1$ and $\wcol_r(G_{r,t}) = \Omega(r^t)$.
Since rooted $2$-treedepth is minor-monotone\footnote{A graph parameter $p$ is said to be \emph{minor-monotone} if $p(H) \leq p(G)$ for all graphs $H$ and $G$ such that $H$ is a minor of $G$.} (see \Cref{lemma:rtd2_minor}), for every graph $X$ with at least one edge, $G_{r,\rtd_2(X)-2}$ is $X$-minor-free and $\wcol_r(G_{r,\rtd_2(X)-2}) = \Omega(r^{\rtd_2(X)-2})$ for every positive integer $r$.
This and \Cref{thm:rooted_main} imply that for every graph $X$ with at least one edge,
\[
\rtd_2(X) - 2 \leq f(X) \leq \rtd_2(X) - 1.
\]

Similarly, \Cref{thm:td_tw} is a direct consequence of the following technical statement.

\begin{theorem}\label{thm:rooted_main_bounded_tw}
For every integer $t$ with $t \geq 2$, for every graph $X$ with $\rtd_2(X) \leq t$, there exists an integer $c$ such that for every graph $G$, if $G$ is $X$-minor-free, then for every integer $r$ with $r \geq 2$,
\[
\wcol_r(G) \leq c \cdot (\tw(G)+1) \cdot r^{t-2} \log r.
\]    
\end{theorem}

When $X$ is a planar graph, $X$-minor-free graphs have bounded treewidth by the Grid-Minor Theorem~\cite{GM5}.
Hence, \Cref{thm:rooted_main_bounded_tw} implies that for every planar graph $X$ with at least one edge,
\[
f(X) = \rtd_2(X)-2.
\]

We conclude the introduction with our two favorite problems in the area. 

\begin{ques}
    What is the asymptotic of the maximum of $\wcol_r(G)$ when G is planar? As discussed, it is known to be $\Omega(r^2\log r)$ and $\Oh(r^3)$. 
    In this paper, we show that the lower bound is tight for planar graphs of bounded treewidth.
\end{ques}

For a positive integer $p$, a vertex coloring $\phi$ of a graph $G$ is \emph{$p$-centered} if for every connected subgraph $H$ of $G$, either $\phi$ uses more than $p$ colors in $V(H)$ or there is a color that appears exactly once in $H$. The \emph{$p$-centered chromatic number} of $G$, denoted by $\chi_p(G)$, is the least number of colors in a $p$-centered coloring of $G$. 
Centered colorings are tied with weak coloring numbers and therefore they also characterize classes of bounded expansion. 
However, we seemingly miss the right proof technique to get upper bounds on $\chi_p(G)$ when $G$ excludes a fixed graph as a minor.

\begin{ques}
    Is there a function $g$ such that for every fixed graph $X$, 
    for every $X$-minor-free graph $G$ and for every positive integer $p$, 
    \[
        \chi_p(G) =\Oh\left(p^{g(\td(X))}\right)?
    \]
\end{ques}
All we know is that $\chi_p(G) =\Oh\left(p^{g(|V(X)|)}\right)$ for some function $g$ as proved by Pilipczuk and Siebertz in~\cite{PS19}.

\section{Outline of the proofs}
In this section, we introduce several notions that we use in the proofs of \Cref{thm:rooted_main} and \Cref{thm:rooted_main_bounded_tw} and then we sketch their proofs.

First, we establish basic notation.
For a positive integer $k$, we write $[k]=\{1,\ldots,k\}$ and $[0] = \emptyset$.
The \emph{null graph} is the graph with no vertices. 
All graphs considered in this paper are finite, simple, and undirected. 
Let $G_1, G_2$ be two graphs.
We denote by $G_1 \sqcup G_2$ the disjoint union of $G_1$ and $G_2$,
and by $G_1 \oplus G_2$ the graph obtained from $G_1 \sqcup G_2$ by adding every edge with one endpoint in $V(G_1)$ and the other in $V(G_2)$.
For every positive integer $k$, for every graph $G$, we write $k \cdot G$ for the union of $k$ disjoint copies of $G$.

\subsection{The key parameter}\label{subsec:outine_the_key_parameter} 
Let $G$ be a graph and let $k$ be a nonnegative integer.
A \emph{separation of order $k$} of $G$ is a pair $(A, B)$ of subgraphs of $G$ such that\footnote{For all graphs $A,B$, let $A \cup B = (V(A) \cup V(B), E(A) \cup E(B))$ and $A \cap B = (V(A) \cap V(B), E(A) \cap E(B))$.} $A\cup B = G$,
$E(A \cap B) = \emptyset$, and $|V (A \cap B)| = k$. 
We define recursively a new graph parameter called \emph{rooted $2$-treedepth},
denoted by $\rtd_2$, as follows.
For every graph $G$,
\begin{enumerate}[label={\normalfont (r\arabic*)}]
    \item $\rtd_2(G) = 0$ if $G$ is the null graph, 
    \item $\rtd_2(G) = 1$ if $G$ is a one vertex graph, and otherwise 
    \item $\rtd_2(G)$ is the minimum of $\max\{\rtd_2(A),\rtd_2(B \setminus V(A)) + |V(A) \cap V(B)|\}$ over all separations $(A,B)$ of $G$ of order at most one with $V(A) \neq \emptyset$ and $V(B) \setminus V(A) \neq \emptyset$. 
\end{enumerate}

Another way to understand $\rtd_2$ is through ``natural'' separations of the graph defined by its block decomposition. When $G$ is not connected, $\rtd_2(G)$ is realized by $\rtd_2(C)$, where $C$ is a component of $G$ for which the value of $\rtd_2$ is the greatest. When $G$ consists of a single block, a separation $(A,B)$ of $G$ with $V(A)\neq \emptyset$ and $V(B)-V(A)\neq \emptyset$ of order at most one is such that $V(A)=\{v\}$ and $B = G$. Therefore, $\rtd_2(G) = \rtd_2(G-v)+1$. When $G$ consists of multiple blocks, the minimum of $\max\{\rtd_2(A),\rtd_2(B \setminus V(A)) + |V(A) \cap V(B)|\}$ is reached for separations where $V(A)\cap V(B)$ consists of exactly one cut-vertex of $G$. The above can be summarized as the following properties.
For every graph $G$,
\begin{enumerate}[resume*]
    \item $\rtd_2(G)$ is the maximum of $\rtd_2(C)$ over all components $C$ of $G$ when $G$ is not connected, 
    \item $\rtd_2(G)$ is the minimum of $\rtd_2(G-v)+1$ over all vertices $v$ of $G$ when $G$ consists of one block, 
    \item $\rtd_2(G)$ is the minimum of $\max\{\rtd_2(A),\rtd_2(B \setminus V(A)) + 1\}$ over all separations $(A,B)$ of $G$ of order one with $V(A) \cap V(B)$ consisting of a cut-vertex, when $G$ is connected and consists of more than one block. 
\end{enumerate}
Moreover, observe that $\rtd_2(G) \leq \max\{\rtd_2(A),\rtd_2(B \setminus V(A)) + |V(A) \cap V(B)|\}$ for every separation $(A,B)$ of $G$ of order at most one.
In particular, for every $u \in V(G)$, $(G[\{u\}],G)$ is a separation of $G$ of order one, and so,
\begin{enumerate}[resume*]
    \item $\rtd_2(G) \leq 1+\rtd_2(G-u)$. 
\end{enumerate}
Finally, vertices of degree $1$ can not increase rooted $2$-treedepth of a graph.
\begin{enumerate}[resume*]
    \item $\rtd_2(G) \leq \max\{2,\rtd_2(G-u)\}$ for every $u \in V(G)$ of degree at most $1$. 
\end{enumerate}

Rooted $2$-treedepth has several interesting properties: 
it is minor-monotone, see \Cref{lemma:rtd2_minor};
and it is also tied to $2$-treedepth.
More precisely, for every graph $G$, we have $\td_2(G) \leq \rtd_2(G) \leq \max\{1, 2\td_2(G)-2\}$ -- see~\Cref{lemma:rtd_2-and-td_2}, and these inequalities are tight -- see~\Cref{lemma:construction_rtd2_td2}.

\begin{figure}[tp]
    \centering 
    \includegraphics[scale=1]{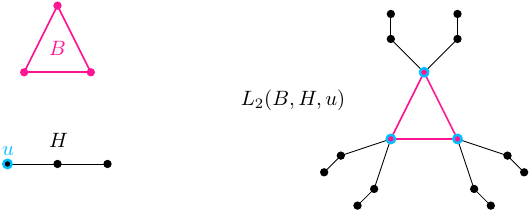} 
    \caption{An construction of $L_2(B,H,u)$, where $B$ is a triangle and $H$ is a path on three vertices with $u$ being one of its endpoints.
    } \label{fig:L-construction}
\end{figure}

The parameter $\rtd_2$ originates from a construction by Grohe et al.~\cite{Grohe15} of graphs $G_{r,t}$ for all nonnegative integers $r$ and $t$ such that
\[\tw(G_{r,t}) = t\quad\textrm{and}\quad\wcol_r(G_{r,t}) = \binom{r+t}{t}.\]
We now recall this construction.
Let $d$ be a positive integer, let $B$ and $H$ be two graphs, and let $u$ be a vertex of $H$.
We define $L_d(B,H,u)$ as the graph obtained in the following process.
Take a copy of $B$ and $d|V(B)|$ copies of $H$.
Label the latter $H_{i,x}$ for each $i \in [d]$ and $x \in V(B)$.
Next, for each $x \in V(B)$ identify $x$ and $u$ in each $H_{i,x}$ for $i \in [d]$.
See \Cref{fig:L-construction}.

For all nonnegative 
integers $r,t$ the graph $G_{r,t}$ is defined recursively for all nonnegative integers $r$ and $t$ by
\[
\left\{
\begin{aligned}
    G_{0,t} &= G_{r,0} = K_1, \\
    G_{r,t} &= L_{\binom{r+t}{t}}(G_{r-1,t}, K_1 \oplus G_{r,t-1}, u) \text{ if $r,t > 0$, }
\end{aligned}
\right.
\]
where $u$ is the vertex of $K_1$ in $K_1 \oplus G_{r,t-1}$.
Observe that $G_{r,1}$ is a tree for all nonnegative integers~$r$.

One can show that for all nonnegative integers $r,t$ and positive integer $d$,
\[\rtd_2\big(L_d(G_{r-1,t}, K_1 \oplus G_{r,t-1}, u)\big) \leq \max\{\rtd_2(G_{r-1,t}), 1 + \rtd_2(G_{r,t-1})\}.\] 
Therefore, by induction, $\rtd_2(G_{r,t})=t+1$. 
In fact, this construction is universal for graphs of rooted $2$-treedepth at most $t+1$ (see~\Cref{lemma:rtd2_and_Grt}).
Namely, for every graph $G$, $\rtd_2(G)\leq t+1$ if and only if $G$ is isomorphic to a subgraph of $G_{r,t}$ for some nonnegative integer $r$.
Since rooted $2$-treedepth is minor-monotone, we deduce that for every graph $X$, $G_{r,\rtd_2(X)-2}$ is $X$-minor-free.
It follows that there are $X$-minor-free graphs with $r$-th weak coloring numbers in $\Omega(r^{\rtd_2(X)-2})$.
Hence, \Cref{thm:rooted_main} yields that for every nonnegative integer $t$, if a minor-closed class of graphs contains graphs with $r$-th weak coloring in $\omega(r^{t} \log r)$,
then it contains $G_{r,t}$ for every nonnegative integer $r$. 
As a consequence, in the setting of minor-closed graphs classes, our results imply that the family constructed by Grohe et al.\
is, up to an $\Oh(r\log r)$ factor, the unique construction of graphs with large weak coloring numbers. 
Similarly, in the setting of minor-closed graphs classes of bounded treewidth, the family constructed by Grohe et al.\ is, up to an $\Oh(\log r)$ factor, the unique construction of graphs with large weak coloring numbers.

\subsection{Weak coloring numbers of \texorpdfstring{$(G,S)$}{(G,S)}}

Another key ingredient in our method is a notion of weak coloring numbers focused on a given subset $S$ of vertices of a graph $G$.
Intuitively, we want to order the vertices of $S$ and place them first in the ordering of $V(G)$ so that, whatever the ordering of the other vertices is, 
every vertex weakly reaches a small number of vertices in~$S$.

\begin{figure}[tp]
    \centering 
    \includegraphics[scale=1]{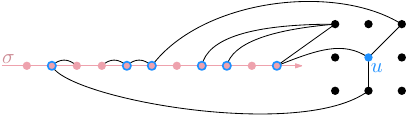} 
    \caption{The pink vertices correspond to the set $S$. The vertices in $S$ highlighted blue are in $\WReach_3[G,S,\sigma,u]$.} \label{fig:wcols}
\end{figure}

Let $G$ be a graph, let $r$ be a nonnegative integer, let $S \subseteq V(G)$, let $\sigma$ be an ordering of $S$, let $u \in V(G)$, and let $v \in S$.
We say that \emph{$v$ is weakly $r$-reachable from $u$ in $(G,S,\sigma)$}
if there is an $u$-$v$ path $P$ in $G$ of length at most $r$ such that $\min_\sigma (V(P) \cap S) = v$.
We denote by $\WReach_r[G,S,\sigma,u]$ the set of all the weakly $r$-reachable vertices from $u$ in $(G,S,\sigma)$ and we write $\wcol_r(G,S,\sigma) = \max_{u \in V(G)} |\WReach_r[G,S,\sigma,u]|$.
Finally, let $\wcol_r(G,S)$ be the minimum value of $\wcol_r(G,S,\sigma)$ among all $\sigma$ orderings of $S$.
For each of the defined objects, we drop $S$ when $S=V(G)$.
Namely, $v$ is weakly $r$-reachable from $u$ in $(G,\sigma)$ whenever $v$ is weakly $r$-reachable from $u$ in $(G,V(G),\sigma)$, $\WReach_r[G,\sigma,u] = \WReach_r[G,V(G),\sigma,u]$, $\wcol_r(G,\sigma) = \wcol_r(G,V(G),\sigma)$, and $\wcol_r(G) = \wcol_r(G,V(G))$.
This matches the definition given in~\Cref{sec:intro}.
See an illustration in \Cref{fig:wcols}.
In \Cref{sec:prelim}, we give many properties of this notion.

\subsection{\texorpdfstring{$\mathcal{F}$}{F}-rich models}
For a graph $H$, a \emph{model} of $H$ in a graph $G$ is a family $\big(B_x \mid x \in V(H)\big)$ of disjoint subsets of $V(G)$ such that
\begin{enumerate}
    \item $G[B_x]$ is connected, for every $x \in V(H)$; and
    \item there is an edge between $B_x$ and $B_y$ in $G$, for every $xy \in E(H)$.
\end{enumerate}
If $G$ has a model of $H$, then we say that $H$ is a \emph{minor} of $G$.
Let $\mathcal{F}$ be a family of connected subgraphs of $G$.
Such a model of $H$ is said to be \emph{$\mathcal{F}$-rich} if we have the additional following property.
\begin{enumerate}[resume*]
    \item For every $x \in V(H)$, there exists $F \in \mathcal{F}$ such that $F \subseteq G[B_x]$.
\end{enumerate}

For example, if $H$ has $k$ vertices and no edges, then $G$ has an $\mathcal{F}$-rich model of $H$ if and only if 
$G$ contains $k$ pairwise disjoint members of $\mathcal{F}$.
Another extreme case is when $\mathcal{F}$ contains every one-vertex subgraph of $G$.
Then every model of $H$ in $G$ is $\mathcal{F}$-rich.

\subsection{Plan of the proof}\label{sec:outline}
We now present the main ideas behind the proofs of \Cref{thm:rooted_main} and \Cref{thm:rooted_main_bounded_tw}.
In order to prove the theorems, we strengthen the statement and instead of considering graphs with no models of $X$, 
we consider graphs with no $\mathcal{F}$-rich models of $X$, 
given a family $\mathcal{F}$ of connected subgraphs of $G$.
This turns out to be very helpful in keeping induction invariant if we carefully choose $\mathcal{F}$.
For some choices of $\mathcal{F}$, excluding an $\mathcal{F}$-rich model of a graph may be a local property.
Therefore, we keep a global property that $G$ excludes $K_k$ as a minor, where $k = |V(X)|$.
The locality of the property suggests that the bound on weak coloring numbers should be local in some sense as well.
To this end, instead of bounding $\wcol_r(G)$, we bound $\wcol_r(G,S)$, where $S$ is a hitting set for $\mathcal{F}$, that is $V(F) \cap S \neq \emptyset$ for every $F \in \mathcal{F}$.
A result in this spirit, that corresponds to excluding every $\mathcal{F}$-rich model of an edgeless graph $X$, is already present in~\cite{DHHJLMMRW24}.
We restate this result below with adjusted notations (see~\Cref{lemma:E-P_Kk-free_original_statement} for the original statement).
\begin{lemma}\label{lemma:E-P_in_Kk_minor_free_graphs_outline}
    There exists a function $\delta$ such that for all positive integers $k,d$, for every connected $K_k$-minor-free graph $G$, 
    for every family $\cgF$ of connected subgraphs of $G$, 
    if there are no $d$ pairwise vertex-disjoint subgraphs in $\mathcal{F}$, 
    then there exists $S \subseteq V(G)$ such that 
    \begin{enumerateOurAlph}
        \item $V(F) \cap S \neq \emptyset$ for every $F \in \mathcal{F}$; 
        \item $G[S]$ is connected;
        \item $\wcol_r(G,S) \leq \delta(k,d) \cdot r$ for every positive integer $r$. 
    \end{enumerateOurAlph}
\end{lemma}

Pushing this idea further, we show the following technical version of \Cref{thm:rooted_main}, which can be seen as the induction setup.

\begin{theorem}\label{thm:main} 
    Let $k$ and $t$ be positive integers with $t \geq 2$.
    Let $X$ be a graph with $\rtd_2(X) \leq t$.
    There exists an integer $c(t,X,k)$ such that
    for every connected $K_k$-minor-free graph $G$, 
    for every family $\mathcal{F}$ of connected subgraphs of $G$,
    if $G$ has no $\mathcal{F}$-rich model of $X$, then
    there exists $S \subseteq V(G)$ such that
    \begin{enumerateOurAlphCapital}
        \item $V(F) \cap S \neq \emptyset$ for every $F \in \mathcal{F}$; \label{thm:main_outline:item:hit}
        \item $G[S]$ is connected; \label{thm:main_outline:item:connected}
        \item $\wcol_r(G,S) \leq c(t,X,k) \cdot r^{t-1} \log r$ for every integer $r$ with $r \geq 2$. \label{thm:main_outline:item:wcol}
    \end{enumerateOurAlphCapital}
\end{theorem}

To recover \Cref{thm:rooted_main} from this last statement, apply it for a family $\mathcal{F}$ containing each one-vertex subgraph of $G$.
Then \ref{thm:main_outline:item:hit} implies that $S = V(G)$, and so, $\wcol_r(G) = \wcol_r(G,S) \leq c(t,X,k) \cdot r^{t-1} \log r$ by \ref{thm:main_outline:item:wcol}.
The assumption that $G$ is $K_k$-minor-free is necessary: for every positive integer $n$, if $\mathcal{F}$ is the family of all the subgraphs of $K_n$ with more than $\frac{n}{2}$
vertices, then $K_n$ has no $\mathcal{F}$-rich model of $K_2$, 
but every hitting set $S$ of $\mathcal{F}$ satisfies $\wcol_r(G,S) \geq \frac{n}{2}$ for every nonnegative integer $r$, 
which is not bounded by a function of $r$. 
\Cref{thm:main_outline:item:connected} is a technical condition that supports the induction. 

\begin{figure}[tp]
    \centering 
    \includegraphics[scale=1]{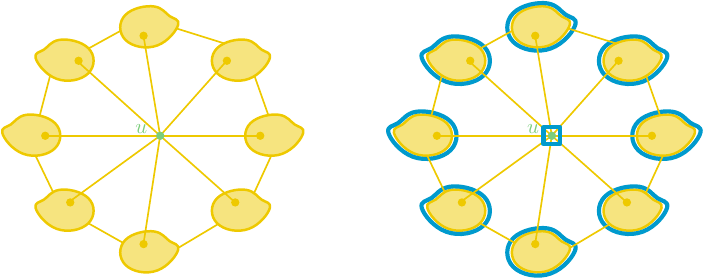} 
    \caption{
    On the left-hand side, we depict an $\mathcal{F}$-rich model of $X$, where $X$ is a cycle on $8$ vertices and $\mathcal{F}$ is the family of all connected subgraphs of $G - \{u\}$ containing a neighbor of $u$ in $G$.
    On the right-hand side, we show how to construct, given an $\mathcal{F}$-rich model of $X$, a model of $K_1 \oplus X$.
    } \label{fig:F-rich-models}
\end{figure}

The proof of \Cref{thm:main} is by induction on $t$. 
Within the inductive step, 
given the result for all $X$ with $\rtd_2(X)= t-1$ we argue that the result holds for graphs of rooted $2$-treedepth equal $t$.
First, we prove it for graphs of the form $K_1\oplus X$ where $\rtd_2(X)=t-1$.
Let $G$ be a $K_1 \oplus X$-minor-free graph.
Let $u$ be a vertex of $G$ and let $\mathcal{F}$ be the family of all the connected subgraphs $H$ of $G-\{u\}$ that contain a neighbor of $u$ in $G$.
Then observe that any $\mathcal{F}$-rich model $\big(B_x \mid x \in V(X)\big)$ of $X$ in $G-\{u\}$ yields a model $\big(C_x \mid x \in V(X)\cup\{s\}\big)$ of $K_1 \oplus X$ 
defined by $C_x = B_x$ for every $x \in V(X)$ and $C_s = \{u\}$ -- see \Cref{fig:F-rich-models}.
Therefore, $G-\{u\}$ has no $\mathcal{F}$-rich model of $X$.
Hence, choosing $\mathcal{F}$ carefully, we can deduce that $G$ has no $\mathcal{F}$-rich model of $X$ knowing that $G$ has no model of $K_1 \oplus X$.
This technique will allow us to prove \Cref{thm:main} for $K_1 \oplus X$, assuming the result for $X$ (see \Cref{claim:adding_apices_J} and \Cref{claim:adding_apices_less_technical_J} in the proof of \Cref{thm:main}).
When we have the result for graphs of the form $K_1\oplus X$, we follow the inductive definition of rooted $2$-treedepth and conclude the full statement of~\Cref{thm:main}, see~\Cref{claim:last_claim}.

Here is the summary of the plan of the proof of \Cref{thm:main}, with the bounds on $\wcol_r(G,S)$ obtained at each step:

\begin{enumerate}[label=\arabic*.]
    \item {\bf Pre-base case: $\rtd_2(X) \leq 1$}, \textit{i.e.} $X$ has no edges; \; \Cref{lemma:E-P_in_Kk_minor_free_graphs}. \hfill $\mathcal{O}(r)$
    \item {\bf Base case: $\rtd_2(X) \leq 2$}, \textit{i.e.} $X$ is a forest; \; \Cref{sec:base}. 
        \begin{enumerate}[label=\alph*.]
            \item {\bf $X$ is a star}; \; Lemmas~\ref{lemma:excluding_a_star_patch} and~\ref{lemma:excluding_a_star_wcol}. \hfill $\mathcal{O}(r \log r)$
            \item {\bf $X$ is a forest}; \; \Cref{lem:new_base_case}. \hfill $\mathcal{O}(r \log r)$
        \end{enumerate}
    \item {\bf Induction: $\rtd_2(X) \leq t$}, assuming the result for $t-1$; \; \Cref{sec:induction}.
        \begin{enumerate}[label=\alph*.]
            \item {\bf $X=K_1 \oplus Y$ for some $Y$ with $\rtd_2(Y) \leq t-1$}; \; Claims~\ref{claim:adding_apices_J} and~\ref{claim:adding_apices_less_technical_J}. \hfill $\mathcal{O}(r^{t-1} \log r)$
            \item {\bf $X$ is any graph with $\rtd_2(X) \leq t$}; \; \Cref{claim:last_claim}. \hfill $\mathcal{O}(r^{t-1} \log r)$
        \end{enumerate}
\end{enumerate}

The logarithmic factor appears already in the base case.
This phenomenon can be explained as follows.
Graphs excluding a fixed tree as a minor have bounded pathwidth~\cite{GM1},
and it is known that paths have logarithmic weak coloring numbers (see~\cite{JM22} or~\Cref{fig:path}).
Combining these ideas, one can show that graphs excluding a fixed tree as a minor have logarithmic weak coloring numbers.
Here, we are working in the more general setting of graphs with no $\mathcal{F}$-rich model of a fixed tree, and so we do not have bounded pathwidth. 
However, using a similar strategy together with \Cref{lemma:E-P_in_Kk_minor_free_graphs_outline}, we prove in our base case an $\Oh(r \log r)$ bound 
for weak coloring numbers of a hitting set of $\mathcal{F}$, in a $K_k$-minor-free graphs with no $\mathcal{F}$-rich model of a fixed tree.

The proof of \Cref{thm:rooted_main_bounded_tw} is very similar, except that a factor $\mathcal{O}(r)$ is saved in the first step.
This comes from the fact that for every graph $G$, for every family $\mathcal{F}$ of connected subgraphs of $G$,
if there are no $d+1$ disjoint members of $\mathcal{F}$, then there exists a hitting set $S \subseteq V(G)$ of $\mathcal{F}$ of size at most $d \cdot (\tw(G)+1)$.
In particular, $\wcol_r(G,S) \leq |S| \leq d \cdot (\tw(G)+1)$ for every nonnegative integer $r$.
This fact (see~\Cref{lemma:helly_property_tree_decomposition}) will replace \Cref{lemma:E-P_in_Kk_minor_free_graphs_outline} in the first step of the proof,
which improves by a factor $\Oh(r)$ the bound obtained in the more general setting of $K_k$-minor-free graphs.
In the proof of~\Cref{thm:rooted_main}, for technical reasons we keep as an invariant that $G[S]$ is connected.
The techniques used in the proof~\Cref{thm:rooted_main_bounded_tw} force us to relax this condition slightly.
Now, we assert that for every component $C$ of $G - S$, $N_G(V(C))$ intersects only few components of $G - V(C)$.
This makes the proof slightly more technical.
However, the main ideas are the same, and substantial parts of the proofs overlap.

\subsection{Organisation of the paper}
In \Cref{sec:prelim}, we introduce notation and some simple properties that will be used throughout the paper.
In \Cref{sec:rtd2}, we prove several properties of the rooted $2$-treedepth, and show its connection with the construction $\{G_{r,t}\}_{r,t \geq 0}$~\cite{Grohe15}.
In \Cref{sec:base}, we prove the base case of \Cref{thm:main} when $X$ is a forest.
In \Cref{sec:induction}, we prove \Cref{thm:main}.
Finally, in \Cref{sec:bounded_tw}, we prove \Cref{thm:rooted_main_bounded_tw} with the same method.

\section{Preliminaries}\label{sec:prelim}

Let $G$ be a graph and let $A,B,Z \subseteq V(G)$.
We say that $Z$ \emph{separates $A$ and $B$ in $G$} if no component of $G-Z$ intersects both $A$ and $B$.

A collection $\mathcal{P}$ of subsets of a non-empty set $S$ is a \emph{partition} of $S$ 
if elements of $\mathcal{P}$ are non-empty, pairwise disjoint, and $\bigcup \mathcal{P} = S$.
A sequence $(P_0,\dots,P_m)$ of subsets of a set $S$ is an \emph{ordered partition} of $S$ if $\{P_i\}_{i \in \{0, \dots, m\}}$ is a partition of $S$.
Given a graph $G$ and a partition $\mathcal{P}$ of $V(G)$, the \emph{quotient graph} $G/\mathcal{P}$ is the graph with the vertex set $\mathcal{P}$ 
and two distinct $P,P' \in \mathcal{P}$ 
are adjacent in $G/\mathcal{P}$ if there are $u \in P$ and $u' \in P'$ such that $uu'$ is an edge in $G$.

A \emph{layering} of a graph $G$ is an ordered partition $(P_0,\dots,P_\ell)$ of $V(G)$ such that for every edge $uv$ in $G$ either there is $i \in \{0,\dots,\ell\}$ with $u,v \in P_i$ or there is $i \in \{0,\dots,\ell-1\}$ with $u \in P_i$ and $v \in P_{i+1}$.
A \emph{tree partition} of a graph $G$ is a pair $(T,\mathcal{P})$, where $T$ is a tree and $\mathcal{P} = (P_{x} \mid x \in V(T))$ is a partition of $V(G)$ such that for every edge $uv$ in $G$ either there is $x \in V(T)$ with $u,v \in P_x$ or there is an edge $xy$ in $T$ with $u \in P_x$ and $v \in P_y$.

Let $G$ be a graph.
For $X, Y \subset V(G)$, an \emph{$X$-$Y$ path} is a path in $G$ that is either a one-vertex path with the vertex in $X \cap Y$ or a path with one endpoint in $X$ and the other endpoint in $Y$ such that no internal vertices are in $X \cup Y$.
When $u,v \in V(G)$, instead of $\{u\}$-$\{v\}$ path we write $u$-$v$ path for short.
The \emph{length} of a path $P$ 
is the number of edges of $P$.
A path $P$ is a \emph{geodesic} in $G$ if it is a shortest path between its endpoints in $G$.
The \emph{distance} between two vertices $u$ and $v$ in $G$, denoted by $\dist_G(u,v)$, is the length of a $u$-$v$ geodesic in $G$ when it exists, and $+\infty$ otherwise.
Let $u$ be a vertex of $G$.
The \emph{neighborhood} of $u$ in $G$, denoted by $N_G(u)$, is the set $\{v \in V(G) \mid uv \in E(G)\}$.
For every set of vertices $X$ of $G$, let $N_G(X)=\bigcup_{u \in X} N_G(u) \setminus X$.
For every positive integer $r$, we denote by $N_G^r[u] = \{v \in V(G) \mid \dist_G(u,v) \leq r\}$.
The following lemma is folklore, see e.g.\ \cite[Lemma 23]{DHHJLMMRW24} for a proof.

\begin{lemma}
\label{lemma:intersection_ball_with_geodesic}
    Let $G$ be a graph and $r$ be a nonnegative integer. 
    For every geodesic $Q$ in $G$ and for every vertex $v \in V(G)$,
    \[
    |N^r[v] \cap V(Q)| \leq 2r+1.
    \]
\end{lemma}

An \emph{ordering} $\sigma$ of a finite set $E$ is a sequence $(x_1, \dots, x_{|E|})$ of all the elements of $E$. 
For all $x,y \in E$, we write $x \leq_\sigma y$ when there are $i,j \in [|E|]$ such that $x_i = x$, $x_j = y$, and $i\leq j$. 
We also write $\min_{\sigma} E = x_1$ and $\max_{\sigma} E = x_{|E|}$.
When $F \subseteq E$, we write $\sigma\vert_F$ for the \emph{restriction} of $\sigma$ to $F$ that is defined as the ordering of $F$ such that $x \leq_{\sigma\vert_F} y$ if and only if $x \leq_\sigma y$ for all $x,y \in F$. 
For every nonempty $F \subseteq E$, we define 
$\min_{\sigma}F = \min_{\sigma|F} F$ and $\max_{\sigma}F = \max_{\sigma|F} F$.
For all $F,F' \subset E$, we write $F <_\sigma F'$ whenever for all $x \in F$ and $y \in F'$, we have $x <_\sigma y$.
If $\sigma'$ is an ordering of $F \subseteq E$, we say that $\sigma$ \emph{extends} $\sigma'$ if $\sigma\vert_F = \sigma'$.

When $H$ is a subgraph of a graph $G$ and $\mathcal{F}$ is a family of subgraphs of $G$, we denote by $\mathcal{F}\vert_H$ the family $\{F \in \mathcal{F} \mid F \subseteq H\}$. 

Next, we state a bunch of simple observations concerning the notion of weak coloring numbers in the version presented above.

\begin{obs}\label{obs:wcol_components}
    Let $G$ be a graph and let $S \subset V(G)$.
    Let $G'$ consist of all the components of $G$ that contain a vertex from $S$.
    For every nonnegative integer $r$, we have
        \[\wcol_r(G,S) = \wcol_r(G',S).\]
\end{obs}

\begin{obs}\label{obs:wcol_components2}
    Let $G$ be a graph, $S \subset V(G)$, and $\mathcal{C}$ be the family of components of $G$.
    For every nonnegative integer $r$, we have
        \[\wcol_r(G,S) = \max_{C \in \mathcal{C}}\wcol_r(C,S \cap V(C)).\]
\end{obs}

\begin{obs}\label{obs:wcol_union}
    Let $G$ be a graph and let $S,S' \subset V(G)$.
    For every nonnegative integer $r$, we have
        \[\wcol_r(G,S\cup S') \leq \wcol_r(G,S) + \wcol_r(G-S,S'-S).\]
\end{obs}

Observations \ref{obs:wcol_components} and \ref{obs:wcol_components2} are clear from the definition and to see~\Cref{obs:wcol_union}, it suffices to order all the vertices of $S$ before all the vertices of $S'$.

Geodesics are a useful tool when bounding weak coloring numbers.
For instance, \Cref{lemma:intersection_ball_with_geodesic} implies the following.

\begin{obs}\label{obs:geodesics}
    Let $G$ be a graph, $S \subset V(G)$, $\ell$ be a positive integer, and $Q_1,\dots,Q_\ell$ be geodesics in $G$.
    For every nonnegative integer $r$, we have
        \[\wcol_r(G,S \cup V(Q_1) \cup \dots \cup V(Q_\ell)) \leq \wcol_r(G,S) + \ell(2r+1).\]
\end{obs}
This inequality is witnessed by an ordering of $S\cup V(Q_1)\cup\cdots V(Q_{\ell})$ obtained from the ordering $\sigma$ of $S$ witnessing $\wcol_r(G,S)$ by putting vertices from $(V(Q_1)\cup\cdots \cup V(Q_{\ell})) \setminus S$ arbitrarily.

Note that $\wcol_r(G,S)$ is not monotone with respect to $S$.
For example, let $G$ be a star with the root $v$.
We have $V(G) \setminus \{v\} \subset V(G)$, however, $\wcol_1(G,V(G) \setminus \{v\}) = |V(G)|$ and $\wcol_1(G,V(G)) = 2$.
On the other hand, our version of weak coloring numbers is monotone in the following sense.

\begin{obs}\label{obs:monotone}
    Let $G$ be a graph, $S \subset V(G)$, and $U \subset V(G)$.
    For every nonnegative integer~$r$, we have
        \[\wcol_r(G - U,S - U) \leq \wcol_r(G,S).\]
\end{obs}

The ideas of Observations~\ref{obs:geodesics}~and~\ref{obs:monotone} can be combined to obtain another property.

\begin{obs}\label{obs:geodesics_plus_monotone}
    Let $G$ be a graph, $S \subset V(G)$, $U \subset V(G)$,  $\ell$ be a positive integer, and $Q_1,\dots,Q_\ell$ be geodesics in $G$.
    For every nonnegative integer $r$, we have
        \[\wcol_r(G-U,(S \cup V(Q_1) \cup \dots \cup V(Q_\ell))-U) \leq \wcol_r(G-U,S-U) + \ell(2r+1).\]
\end{obs}

Finally, we can decide to place a fixed subset $A$ of vertices first in the ordering and then consider the geodesic paths in the remaining graph $G-A$. 
This turns out to be an important trick.
\begin{obs}\label{obs:geodesics_in_G-A}
    Let $G$ be a graph, $A \subset V(G)$, $\ell$ be a positive integer, and $Q_1,\dots,Q_\ell$ be geodesics in $G - A$.
    For every nonnegative integer $r$, we have
        \[\wcol_r(G,A \cup V(Q_1) \cup \dots \cup V(Q_\ell)) \leq |A| + \ell(2r+1).\]
\end{obs}
To see that the observation holds, 
we just take an arbitrary ordering of $A\cup V(Q_1)\cup\cdots\cup V(Q_{\ell})$ with vertices of $A$ preceding vertices of $V(Q_1)\cup\cdots\cup V(Q_{\ell})$.

Now, we present one of the key basic tools in the proof of \Cref{thm:rooted_main}.


\begin{lemma}[\Cref{lemma:E-P_in_Kk_minor_free_graphs_outline} restated]\label{lemma:E-P_in_Kk_minor_free_graphs}
    There exists a function $\delta$ such that for all positive integers $k,d$, for every connected $K_k$-minor-free graph $G$, for every family $\cgF$ of connected subgraphs of $G$, 
    either there are $d$ pairwise vertex-disjoint subgraphs in $\mathcal{F}$, 
    or there exists $S \subseteq V(G)$ such that 
    \begin{enumerateOurAlph}
        \item $V(F) \cap S \neq \emptyset$ for every $F \in \mathcal{F}$;\label{lem:EP:item:hit}
        \item $G[S]$ is connected; \label{lem:EP:item:connected}
        \item $\wcol_r(G,S) \leq \delta(k,d) \cdot r$ for every positive integer $r$. \label{lem:EP:item:wcol}
    \end{enumerateOurAlph}
\end{lemma}

\Cref{lemma:E-P_in_Kk_minor_free_graphs} is a consequence of the following statement from~\cite{DHHJLMMRW24}, which relies on the Graph Minor Structure Theorem by Robertson and Seymour.

\begin{lemma}[{\cite[Lemma~21]{DHHJLMMRW24}}]\label{lemma:E-P_Kk-free_original_statement}
    There exists a function $\gamma$ such that for all positive integers $k,d$, for every $K_k$-minor-free graph $G$, for every family $\cgF$ of connected subgraphs of $G$ either
    \begin{enumerate}[label=\normalfont(\arabic*)]
        \item there are $d$ pairwise vertex-disjoint subgraphs in $\mathcal{F}$, or \label{lemma:E-P_Kk-free_original_statement:item:i}
        \item there exists $A \subseteq V(G)$ with $|A| \leq (d-1) \gamma(k)$ and there exists a subgraph $X$ of $G$ which is the union of at most $(d-1)^2\gamma(k)$ geodesics in $G-A$,
        such that for every $F \in \cgF$ we have $V(F) \cap (V(X) \cup A) \neq \emptyset$.  \label{lemma:E-P_Kk-free_original_statement:item:ii}
    \end{enumerate}
\end{lemma}

\begin{proof}[Proof of \Cref{lemma:E-P_in_Kk_minor_free_graphs}]
    Let $\delta(k,d) = 12(d-1)^2\gamma(k)$
    where $\gamma$ is the function from \Cref{lemma:E-P_Kk-free_original_statement}.
    Let $G$ be a $K_k$-minor-free graph and let $\mathcal{F}$ be a family of connected subgraphs of $G$.
    Suppose that there are no $d$ pairwise disjoint members of $\mathcal{F}$, and hence,~\Cref{lemma:E-P_Kk-free_original_statement}.\ref{lemma:E-P_Kk-free_original_statement:item:ii} holds, yielding $A \subset V(G)$ and a subgraph $X$ of $G$ such that $|A| \leq (d-1) \gamma(k)$ and     
    $X$ is the union of at most $(d-1)^2 \gamma(k)$ geodesics in $G-A$.
    Note that $G[A \cup V(X)]$ has at most $|A| + (d-1)^2 \gamma(k)$ components.
    Let $Q_1, \dots, Q_\ell$ be a family of at most $(d-1) \gamma(k) + (d-1)^2 \gamma(k) -1$ geodesics in $G$ such that the set $S = A \cup V(X) \cup \bigcup_{i \in [\ell]} V(Q_i)$ induces a connected subgraph in $G$.
    In particular, $\ell \leq 2(d-1)^2\gamma(k)$.
    By Observations~\ref{obs:geodesics}~and~\ref{obs:geodesics_in_G-A}, for every positive integer $r$,
    \begin{align*}
        \wcol_r(G,S) &\leq \wcol_r(G,A \cup V(X)) + \ell(2r+1)\\
        &\leq |A| + (d-1)^2\gamma(k)(2r+1) + \ell(2r+1) \\
        &\leq (d-1)\gamma(k) + (d-1)^2\gamma(k)(2r+1) + 2(d-1)^2\gamma(k)(2r+1)\\
        &\leq 4(d-1)^2\gamma(k)(2r+1) \leq \delta(k,d)r. \qedhere
    \end{align*}    
\end{proof}

\begin{figure}[tp]
    \centering 
    \includegraphics[scale=1]{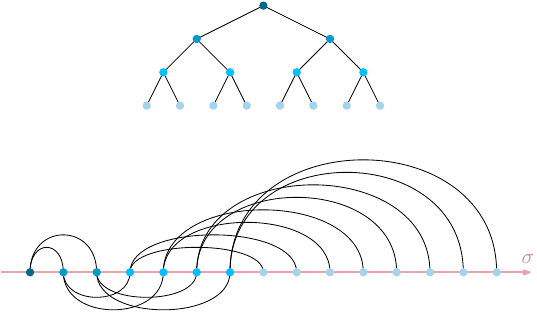} 
    \caption{An example of an eliminating ordering of a complete binary tree of height $3$.} 
    \label{fig:tree}
\end{figure}

It is easy to derive an upper bound on weak coloring numbers of trees.
It suffices to root a given tree and order the vertices in an elimination order.
Namely, an \emph{elimination ordering} of a tree $T$ rooted in $s \in V(T)$ is an ordering $(x_1, \dots, x_{|V(T)|})$ of $V(T)$
such that $x_1 = s$ and for every $i \in \{2,\dots,|V(T)|\}$, $N(x_i) \cap \{x_j \mid j \in [i-1]\} = \{y\}$ where $y$ is the parent of $x_i$.
Note that in such an ordering a vertex weakly reaches only its ancestors. 
See~\Cref{fig:tree}.

\begin{obs}\label{obs:wcols:trees}
    Let $T$ be a tree. For every positive integer $r$, we have
        \[\wcol_r(T) \leq r+1.\]
    More precisely, for every elimination ordering $\sigma$ of $T$, for every $u \in V(T)$, and for every positive integer $r$ we have
        \[\WReach_r[T,\sigma,u] \subseteq A_{r,x}\]
    where $A_{r,x}$ is the set of ancestors of $x$, including $x$ itself, in distance at most $r$.
    Note that $|A_{r,x}| \leq r+1$.
\end{obs}

This idea can be generalized to elimination orderings of tree decompositions -- see~\cite{Grohe15}.

\section{Rooted \texorpdfstring{$2$}{2}-treedepth}\label{sec:rtd2}
First, we recall the definition of rooted $2$-treedepth.
For every graph $G$,
\begin{enumerate}[label={\normalfont (r\arabic*)}]
    \item $\rtd_2(G) = 0$ if $G$ is the null graph,\label{def:rtd2:item:null-graph}
    \item $\rtd_2(G) = 1$ if $G$ is a one vertex graph, and otherwise \label{def:rtd2:item:one-vertex}
    \item $\rtd_2(G)$ is the minimum of $\max\big\{\rtd_2(A),\rtd_2\big(B \setminus V(A)\big) + |V(A) \cap V(B)|\big\}$ over all separations $(A,B)$ of $G$ of order at most one with $V(A) \neq \emptyset$ and $V(B) \setminus V(A) \neq \emptyset$.\label{def:rtd2:item:separations}
\end{enumerate}

As mentioned in \Cref{subsec:outine_the_key_parameter}, the following properties are direct consequences of the definition.
For every graph $G$,
\begin{enumerate}[resume*]
    \item $\rtd_2(G)$ is the maximum of $\rtd_2(C)$ over all components $C$ of $G$ when $G$ is not connected;\label{rtd2:item:components}
    \item $\rtd_2(G)$ is the minimum of $\rtd_2(G-v)+1$ over all vertices $v$ of $G$ when $G$ consists of one block;\label{rtd2:item:one-block}
    \item $\rtd_2(G)$ is the minimum of $\max\{\rtd_2(A),\rtd_2(B \setminus V(A)) + 1\}$ over all separations $(A,B)$ of $G$ of order one with $V(A) \cap V(B)$ consisting of a cut-vertex, when $G$ is connected and consists of more than one block;\label{rtd2:item:more-blocks}
    \item $\rtd_2(G) \leq 1+\rtd_2(G-u)$;\label{rtd2:item:removing-vertex}
    \item $\rtd_2(G) \leq \max\{2,\rtd_2(G-u)\}$ for every $u \in V(G)$ of degree at most $1$.\label{rtd2:item:degree1}
\end{enumerate}

Item~\ref{def:rtd2:item:separations} in the definition can be in fact strengthened in the following way.
For every graph $G$,
\begin{enumerate}[resume*]
    \item $\rtd_2(G)$ is the minimum of $\max\big\{\rtd_2(A),\rtd_2\big(B \setminus V(A)\big) + |V(A) \cap V(B)|\big\}$ over all separations $(A,B)$ of $G$ of order at most one with $V(A) \neq \emptyset$ and $V(B) \setminus V(A) \neq \emptyset$ such that $B$ is a block.\label{def:rtd2:item:separations-stronger}
\end{enumerate}
To see that, consider a separation $(A,B)$ of $G$ of order one with $V(A) \neq \emptyset$ and $V(B) \setminus V(A) \neq \emptyset$ such that $B$ is a block.
By \ref{def:rtd2:item:separations}, $\rtd_2(G) \leq \max\big\{\rtd_2(A),\rtd_2\big(B \setminus V(A)\big) + 1\big\}$.
For the other inequality,
observe that we can assume $G$ connected by \ref{rtd2:item:components}.
Suppose that $\rtd_2(G) = \max\big\{\rtd_2(A),\rtd_2\big(B \setminus V(A)\big) + 1\big\}$ where $(A,B)$ is a separation of $G$ of order one with $V(A) \neq \emptyset$ and $V(B) \setminus V(A) \neq \emptyset$.
Let $(A',B')$ be a separation of $G$ of order one with $V(A') \neq \emptyset$ and $V(B') \setminus V(A') \neq \emptyset$ such that $B'$ is a block, $A \subset A'$ and $B' \subseteq B$.
Then $\rtd_2(A') \leq \rtd_2(G)$ and $\rtd_2(B' \setminus V(A'))+1 \leq \rtd_2(B \setminus V(A))+1 \leq \rtd_2(G)$.

As an illustration of the definition of rooted $2$-treedepth, we characterize graphs having the values of $\rtd_2$ in $\{1,2\}$.
First, note the following straightforward observation.
\begin{obs}\label{obs:char:rtd_2=1}
    For every graph $G$, we have $\rtd_2(G) \leq 1$ if and only if $G$ has no edges.
\end{obs}
Next, we show that for every tree $T$, we have $\rtd_2(T) \leq 2$.
We proceed by induction on the number of vertices of $T$.
For the base case, $\rtd_2(K_1)=1 \leq 2$.
In general, if $x$ is a leaf of $T$ whose parent is $y$, then the separation $(T-x,T[\{x,y\}])$ witnesses $\rtd_2(T) \leq \max\{\rtd_2(T-x), \rtd_2(K_1)+1\} \leq 2$.

We are about to show that $\rtd_2$ is minor-monotone. 
Note that this yields a characterization of graphs with $\rtd_2$ at most two.
Namely, $\rtd_2(G) \leq 2$ if and only if $G$ is a forest.
Indeed, observe that $\rtd_2(K_3) = 3$.

\begin{obs}\label{obs:char:rtd_2=2}
    For every graph $G$, we have $\rtd_2(G) \leq 2$ if and only if $G$ is a forest.
\end{obs}

\begin{lemma}\label{lemma:rtd2_minor}
    For all graphs $G,H$, if $H$ is a minor of $G$, then
    \[
    \rtd_2(H) \leq \rtd_2(G).
    \]
\end{lemma}
\begin{proof}
    We proceed by induction on $|V(G)|$.
    When $|V(G)| \leq 1$, then the assertion holds. 
    Hence, let $G$ be a graph on at least two vertices.
    There is a separation $(A,B)$ of $G$ of order at most one such that $\rtd_2(G) = \max\{\rtd_2(A),\rtd_2(B \setminus V(A)) + |V(A) \cap V(B)|\}$, $V(B) \setminus V(A) \neq \emptyset$, and $V(A) \neq \emptyset$.
    In particular, $\rtd_2(A) \leq \rtd_2(G)$ and $\rtd_2(B \setminus V(A)) + |V(A) \cap V(B)|\leq \rtd_2(G)$.

    We claim that $H$ has a separation $(A',B')$ such that
    $A'$ is a minor of $A$, $B' \setminus V(A')$ is a minor of $B \setminus V(A)$, and $|V(A') \cap V(B')| \leq |V(A) \cap V(B)|$.
    Indeed, let $\big(C_x \mid x \in V(H)\big)$ be a model of $H$ in $G$, and let
        \begin{align*}
            A' &= H[\{x \in V(H) \mid C_x \cap V(A) \neq \emptyset\}]\\
            B' &= H[\{x \in V(H) \mid C_x \cap V(B) \neq \emptyset\}].
        \end{align*}
    Note that $(A',B')$ is a separation of $H$. 
    Moreover, by construction $|V(A') \cap V(B')| \leq |V(A) \cap V(B)|$.
    In addition, since $|V(A) \cap V(B)| \leq 1$, $A'$ is a minor of $A$ and $B'-V(A')$ is a minor of $B - V(A)$.
    Observe that since $|V(A') \cap V(B')| \leq 1$, $\big(C_x \cap V(A) \mid x \in V(A')\big)$ is a model of $A'$ in $A$, and $\big(C_x \cap (V(B) \setminus V(A)) \mid x \in V(B') \setminus V(A')\big)$ is a model of $B' - V(A')$ in $B$.

    Since $V(B) \setminus V(A) \neq \emptyset$, and $V(A) \neq \emptyset$, we have $|V(A)| < |V(G)|$ and $|V(B-V(A))| < |V(G)|$.
    Therefore, by induction hypothesis, since $A'$ is a minor of $A$ and $B'-V(A')$ is a minor of $B - V(A)$,
    $\rtd_2(A') \leq \rtd_2(A)$ and $\rtd_2(B' - V(A')) \leq \rtd_2(B - V(A))$.
    We deduce that
    \begin{align*}
        \rtd_2(H) 
        &\leq \max\{\rtd_2(A'),\rtd_2(B' \setminus V(A')) + |V(A')\cap V(B')|\} \\
        &\leq \max\{\rtd_2(A), \rtd_2(B \setminus V(A)) + |V(A) \cap V(B)|\} \\
        &= \rtd_2(G). \qedhere
    \end{align*}
\end{proof}

It is evident from the definitions that the parameters $\td_2$ and $\rtd_2$ are closely related.
Indeed, for instance, $\td_2$ is minor-monotone too -- one can see this by following the above proof with a slight modification (we do not give an explicit proof since we never use this fact).
In the next part of this section, we discuss relations between $\rtd_2$ and $\td_2$.
Namely, we show that the two parameters are linearly tied.

\begin{lemma}\label{lemma:rtd_2-and-td_2}
    For every graph $G$ with at least one edge,
    \[
    \td_2(G) \leq \rtd_2(G) \leq 2 \td_2(G) - 2.
    \]
\end{lemma}

\begin{proof}
    First, we prove that $\td_2(G) \leq \rtd_2(G)$ for every graph $G$.
    We proceed by induction on $|V(G)|$.
    When $G$ is a null graph, we have $\td_2(G) = \rtd_2(G) = 0$ and when $G$ is a one-vertex graph, we have $\td_2(G) = \rtd_2(G) = 1$.
    Thus, we assume that $|V(G)| \geq 2$.
    If $G$ consists of one block, then by~\ref{rtd2:item:one-block} and induction hypothesis,
        \[\td_2(G) = \min_{v\in V(G)} \td_2(G - v) + 1 \leq \min_{v\in V(G)} \rtd_2(G - v) + 1 = \rtd_2(G).\]
    If $G$ consists of blocks $B_1,\dots,B_k$ for $k > 1$, then by induction hypothesis,
        \[\td_2(G) = \max_{i \in [k]} \td_2(B_i) \leq \max_{i \in [k]} \rtd_2(B_i) \leq \rtd_2(G).\]

    Now, we prove the other inequality for every graph $G$ with at least one edge. 
    We again proceed by induction on $|V(G)|$.
    If $\td_2(G)=2$, then $G$ is a forest with at least one edge, and so as mentioned earlier $\rtd_2(G) = \td_2(G) = 2$.
    Now assume that $\td_2(G) \geq 3$, and so in particular $|V(G)| \geq 3$, and that the result holds for smaller graphs.
    In particular, for every nonnull graph $H$ with $|V(H)|<|V(G)|$, either $H$ has no edge and so $\rtd_2(H)=\td_2(H)=1$,
    or $\rtd_2(H) \leq 2\td_2(H)-2$. In both cases, $\rtd_2(H) \leq \max\{1,2\td_2(H)-2\}$.
    By \ref{def:rtd2:item:separations-stronger},
    there is a separation $(A,B)$ of $G$ of order at most one such that $\rtd_2(G) = \max\{\rtd_2(A),\rtd_2(B \setminus V(A)) + |V(A) \cap V(B)|\}$, $V(B) \setminus V(A) \neq \emptyset$, $V(A) \neq \emptyset$, and $B$ is a block of $G$.
    If $|V(A) \cap V(B)| = 0$, then $B-V(A) = B$ and so
    \begin{align*}
        \rtd_2(G) 
        &= \max\{\rtd_2(A),\rtd_2(B)\} \\
        &\leq \max\{\max\{1,2\td_2(A)-2\},\max\{1,2\td_2(B)-2\}\} \\
        &= \max\{1,2\max\{\td_2(A),\td_2(B)\} - 2\} \\
        &= 2\td_2(G)-2.
    \end{align*}
    Therefore, we assume that $|V(A) \cap V(B)|=1$ and $V(A) \cap V(B) = \{u\}$.
    There exists $v \in V(B)$ such that $\td_2(B-v) = \td_2(B)-1$.
    Then, by~\ref{rtd2:item:removing-vertex},
    \begin{align*}
        \rtd_2(B-u) &\leq \rtd_2(B-u-v) + 1 \leq \rtd_2(B-v)+1 \\
        &\leq \max\{1,2\td_2(B-v)-2\} + 1 \\
        &= \max\{2,2\td_2(B-v)-1\} \\
        &\leq \max\{2, 2\td_2(B)-3\}.
    \end{align*}
    Finally, since $\td_2(G) \geq 3$,
    \begin{align*}
        \rtd_2(G) &= \max\{\rtd_2(A),\rtd_2(B \setminus u) + 1\}\\
        &\leq \max\{\max\{1,2\td_2(A)-2\}, \max\{2, 2\td_2(B)-3\}+1\}\\
        &= \max\{3,2\td_2(A)-2, 2\td_2(B)-3 + 1\}\\
        &\leq \max\{3,2\td_2(G)-2\}\\
        &= 2\td_2(G) - 2.\qedhere
    \end{align*}    
\end{proof}

The bounds in \Cref{lemma:rtd_2-and-td_2} are tight. 
Indeed, for every positive integer $n$, we have $\td_2(K_n) = \rtd_2(K_n) = n$, which witnesses that the first inequality is tight.
For the second one, see \Cref{lemma:construction_rtd2_td2}, which we precede with a simple observation.
Note that this observation is also true for $\td_2$, namely, $\td_2(K_1 \oplus G) = 1 + \td_2(G)$ -- again, the proof is very similar and we omit it.

\begin{obs}\label{lemma:add_apex_increases_rtd2}
    For every graph $G$, 
    \[\rtd_2(K_1 \oplus G) = 1 + \rtd_2(G).\]
\end{obs}

\begin{proof}
    Let $G$ be a graph and let $s$ the vertex of $K_1$ in $K_1 \oplus G$.
    By definition, $\rtd_2(K_1 \oplus G) \leq 1 + \rtd_2(G)$.
    For the other inequality, we proceed by induction on $\rtd_2(G)$.
    The assertion is clear when $G$ is the null graph, thus, assume that $G$ is not the null graph.
    If $G$ is not connected, then $\rtd_2(G) = \rtd_2(C)$ for some component $C$ of $G$, and since $\rtd_2(K_1 \oplus G) \geq \rtd_2(K_1 \oplus C)$, it suffices to show $\rtd_2(K_1 \oplus C) \geq 1 + \rtd_2(C)$.
    Therefore, we assume that $G$ is connected.
    Since $K_1 \oplus G$ is also connected, there is a separation $(A,B)$ of $K_1 \oplus G$ of order one such that $\rtd_2(G) = \max\{\rtd_2(A),\rtd_2(B \setminus V(A)) + 1\}$, $V(B) \setminus V(A) \neq \emptyset$, and $V(A) \neq \emptyset$.
    Since $s$ is adjacent to all other vertices in $K_1 \oplus G$, the only possibility is that $V(A) = \{s\}$ and $V(B) = V(K_1 \oplus G)$.
    It follows that $B-V(A)$ contains a subgraph isomorphic to $G$, and thus, $\rtd_2(K_1 \oplus G) \geq 1 + \rtd_2(B \setminus V(A)) \geq 1 + \rtd_2(G)$.
\end{proof}

\begin{lemma}\label{lemma:construction_rtd2_td2}
    For every integer $k$ with $k \geq 2$, there is a graph $G$ with
    $\td_2(G)\leq k$ and $\rtd_2(G)\geq 2k-2$.
\end{lemma}

\begin{figure}[tp]
        \centering 
        \includegraphics[scale=1]{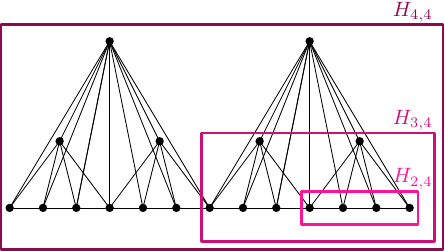} 
        \caption{The proof of \Cref{lemma:construction_rtd2_td2} implies that $\td_2(H_{4,4}) \leq 4$ and $\rtd_2(H_{4,4}) \geq 6$.
        } \label{fig:lower_bound_rtd2_td2}
    \end{figure}

\begin{proof}
    We define inductively graphs $H_{k,\ell}$ with two distinguished vertices $u_{k,\ell}$ and $v_{k,\ell}$ for every integers $k,\ell$ with $k,\ell \geq 2$.
    For $k=2$, $H_{k,\ell}$ is a path on $\ell$ vertices and $u_{k,\ell},v_{k,\ell}$ are its endpoints.
    For $k \geq 3$, $H_{k,\ell}$ is obtained from two disjoint copies $H_1,H_2$ of $K_1\oplus H_{k-1,\ell}$ by identifying the copy of $v_{k-1,\ell}$ in $H_1$ with the copy of $u_{k-1,\ell}$ in $H_2$. The vertices $u_{k,\ell}, v_{k,\ell}$ are then respectively the copy of $u_{k-1,\ell}$ in $H_1$ and the copy of $v_{k-1,\ell}$ in $H_2$.
    See Figure~\ref{fig:lower_bound_rtd2_td2}.
    
    By induction on $k$, we show that $\td_2(H_{k,\ell}) \leq k$ and $\rtd_2(H_{k,\ell}) \geq 2k-2$ for all integers $k,\ell$ with $\ell \geq k \geq 2$.
    When $k=2$, $H_{2,\ell}$ is a path on at least two vertices and so $\td_2(H_{2,\ell}) = \rtd_2(H_{2,\ell}) = 2$.
    Now suppose that $k\geq 3$.
    First, observe that $H_{k,\ell}$ has exactly two blocks $H_1,H_2$, both isomorphic to $K_1 \oplus H_{k-1,\ell}$.
    Hence, $\td_2(H_{k,\ell}) \leq \td_2(K_1 \oplus H_{k-1,\ell}) \leq 1 + \td_2(H_{k-1,\ell}) \leq k$ by induction hypothesis.
    Let $v$ be the unique cut-vertex of $H_{k,\ell}$.
    Since $H_{k,\ell}$ is connected and consists of more than one block, by~\ref{rtd2:item:more-blocks}, there is a separation $(A,B)$ of $H_{k,\ell}$ such that $V(A) \cap V(B) = \{v\}$ and $\rtd_2(H_{k,\ell}) = \max\{\rtd_2(A), \rtd_2(B-v)+1\}$.
    It follows that the graph $B-v$ contains $K_1 \oplus H_{k-1,\ell-1}$ as a subgraph, and so, applying \Cref{lemma:add_apex_increases_rtd2},
    \[\rtd_2(H_{k,\ell-1}) \geq \rtd_2(B-v) + 1 \geq \rtd_2(K_1 \oplus H_{k-1,\ell-1}) + 1 \geq \rtd_2(H_{k-1,\ell-1}) + 2 \geq 2k-2.\]
    This concludes the proof of the lemma.
\end{proof}

In \Cref{sec:intro}, we mentioned several times the construction given in \cite{Grohe15}.
Let us now finally introduce it properly, and later show the equivalent description of rooted $2$-treedepth using this construction.

Let $d$ be a positive integer, let $B, H$ be two graphs, and let $u$ be a vertex of $H$. 
Recall that $L_d(B,H,u)$ is the graph obtained in the following process.
Take a copy of $B$ and $d|V(B)|$ copies of $H$.
Label the latter $H_{i,x}$ for each $i \in [d]$ and $x \in V(B)$.
Next, for each $x \in V(B)$ identify $x$ and $u$ in each $H_{i,x}$ for $i \in [d]$.
See \Cref{fig:L-construction}.

Moreover, for all nonnegative integers $r,t$ the graph $G_{r,t}$ is defined recursively by
\[
\left\{
\begin{aligned}
    G_{0,t} &= G_{r,0} = K_1, \\
    G_{r,t} &= L_{\binom{r+t}{t}}(G_{r-1,t}, K_1 \oplus G_{r,t-1}, u) \text{ if $r,t > 0$,}
\end{aligned}
\right.
\]
where $u$ is the vertex of $K_1$ in $K_1 \oplus G_{r,t-1}$.

\begin{lemma}\label{lemma:rtd2_and_Grt}
For every graph $G$ with at least one edge, $\rtd_2(G)$ is equal to the least $t$ such that there exists a nonnegative integer $r$ with $G \subset G_{r,t-1}$.
\end{lemma}
\begin{proof}
    First, we show that for every 
    $G$ with at least one edge and every integer $t\geq0$ such that there exists $r\geq0$ with $G\subseteq G_{r,t-1}$, we have $\rtd_2(G)\leq t$.
    Since $G$ has an edge, we have $t\geq2$. 
    By \Cref{lemma:rtd2_minor}, it is enough to show that $\rtd_2(G_{r,t-1}) \leq t$ for all integers $r,t$ with $r \geq 0$ and $t \geq 2$.
    We proceed by induction on $t+r$. 
    If $t=2$ or $r=0$, then $G_{r,t-1}$ is a tree and the result follows. 
    Now suppose $t \geq 3, r \geq 1$. 
    By the definition of $G_{r,t-1}$, there is a sequence $H_0, \dots, H_m$ 
    of graphs such that $H_0 = G_{r-1,t-1}$, $H_m = G_{r,t-1}$, and for every $i \in [m]$, $H_i$ contains a block $B_i$ containing a unique 
    cut vertex $r_i$ of $H_i$ in $B_i$, which is such that $B_i - r_i$ is isomorphic to $G_{r,t-2}$ and $H_{i-1} = H_i - V(B_i-r_i)$.

    We claim that $\rtd_2(H_i) \leq t$ for every $i \in \{0,\dots,m\}$.
    We proceed by induction on $i$.
    The main induction hypothesis gives $\rtd_2(H_0) = \rtd_2(G_{r-1,t-1}) \leq t$.
    Now let $i \in [m]$.
    The graph $H_{i}$ has a block $B_{i}$ containing a unique cut vertex $r_{i}$ of $H_{i}$ such that $B_{i} - r_{i}$ is isomorphic to $G_{r,t-2}$ and $H_{i-1} = H_{i} - V(B_{i}-r_{i})$.
    Then $(H_{i-1}, B_{i})$ is a separation of $H_{i}$ of order at most one, and hence by induction hypothesis $\rtd_2(H_{i}) \leq \max\{\rtd_2(H_{i-1}), \rtd_2(G_{r,t-2})+1\} \leq t$.
    For $i=m$, this proves that $\rtd_2(G_{r,t-1}) \leq t$ as asserted.
    
    Next, we prove the opposite implication, i.e., if $\rtd_2(G)\leq t$ then there exists $r\geq0$ such that $G\subseteq G_{r,t-1}$.
    We proceed by induction on $(t, |V(G)|)$ in the lexicographic order. 
    Since $G$ has an edge, we must have $|V(G)|\geq2$ and $t\geq2$.
    Suppose that $t=2$, and so, $\rtd_2(G) \leq 2$, which implies by \Cref{obs:char:rtd_2=2} that $G$ is a forest. 
    There is a vertex $x$ of $G$ of degree at most one.
    Since $\rtd_2(G-x) \leq \rtd_2(G) \leq t$, by induction hypothesis, there exists a nonnegative integer $r$ such that
    $G-x \subseteq G_{r,1}$, and so $G \subseteq G_{r+1,1}$.

    Now, suppose that $t \geq 3$. 
    There is a separation $(A,B)$ of $G$ of order at most one such that $\rtd_2(A) \leq \rtd_2(G)$, $\rtd_2(B \setminus V(A)) \leq \rtd_2(G) - |V(A) \cap V(B)|$, $V(B) \setminus V(A) \neq \emptyset$ and $V(A) \neq \emptyset$.
    We take such a separation $(A, B)$ with $|V(A) \cap V(B)| = 1$ if possible.
    Assume that it is impossible, and so, $|V(A) \cap V(B)| = 0$.
    Then we claim that every component of $G$ has only one vertex.
    Indeed, by~\ref{rtd2:item:components} there is a component $C$ of $G$ such that $\rtd_2(G)=\rtd_2(C)$.
    If $C$ has more than one vertex, there is a separation
    $(A',B')$ of $C$ of order at most one such that
    $\rtd_2(A') \leq \rtd_2(C)$, $\rtd_2(B' \setminus V(A')) \leq \rtd_2(C) - |V(A') \cap V(B')|$, $V(B') \setminus V(A') \neq \emptyset$ and $V(A') \neq \emptyset$.
    Note that $|V(A') \cap V(B')| = 1$ as $C$ is connected.
    However, $\big((G\setminus V(C)) \cup A',B'\big)$ contradicts the impossibility assumption.
    This proves that every component of $G$ has one vertex, and so, $G$ has no edges, which is a contradiction.
    
    We have $|V(A) \cap V(B)| = 1$, $|V(A)| < |V(G)|$, and $|V(B-V(A))| < |V(G)|$.
    Moreover, $\rtd_2(A) \leq t$ and $\rtd_2(B-V(A)) \leq t-1$, therefore, by induction hypothesis,
    there exist nonnegative integers $r$ and $r'$ such that $A \subseteq G_{r,t-1}$ and $B \setminus V(A) \subseteq G_{r',t-2}$.
    It follows that $B \subseteq K_1 \oplus G_{r',t-2}$, and so $G \subseteq G_{\max\{r+1,r'\},t-1}$. 
\end{proof}

We finish this section with another universal construction for graphs of rooted $2$-treedepth at most $t$,
which will be useful in the proofs of \Cref{thm:rooted_main} and \Cref{thm:rooted_main_bounded_tw}.
Given a graph $G$, we will define a graph $T_{h,d}(G)$ for all positive integers $h,d$, whose blocks are all isomorphic to $K_1 \oplus G$.

\begin{figure}[tp]
    \centering 
    \includegraphics[scale=1]{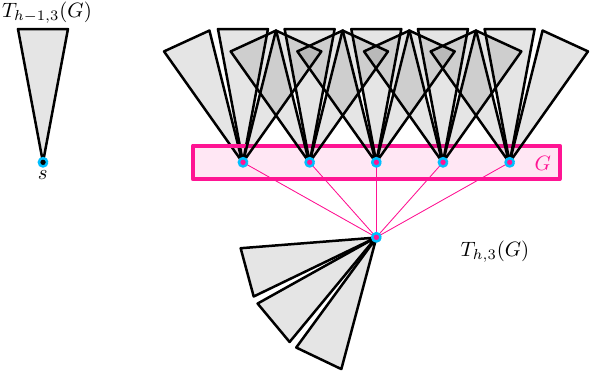} 
    \caption{An illustration of the construction of $T_{h,d}(G)$, which is obtained from $K_1 \oplus G$ by gluing on every vertex $d$ copies of $T_{h-1,d}(G)$.}
    \label{fig:thd}
\end{figure}

Let $G$ be a graph and let $d$ be a positive integer.
For every positive integer $h$, we define the graph $T_{h,d}(G)$ with one distinguished vertex, which we call the \emph{root} of $T_{h,d}(G)$.
When $h = 1$, let $T_{1,d}(G) = K_1 \oplus G$, and let the vertex of $K_1$ be the root.
When $h > 1$, fix a copy of $T_{h-1,d}(G)$ with the root $s$, and let
    \[T_{h,d}(G) = L_d(K_1 \oplus G,T_{h-1,d}(G),s).\]
As the root of $T_{h,d}(G)$, we distinguish the vertex of $K_1$ in the copy of $K_1 \oplus G$ given as the first argument to $L_d$.
See \Cref{fig:thd}.

\begin{lemma}\label{obs:rtd_and_Thd}
    For every graph $G$, for all positive integers $h,d$, we have 
    \[\rtd_2(T_{h,d}(G)) = \rtd_2(G) + 1. \]
\end{lemma}
\begin{proof}
    Let $G$ be a graph and let $h,d$ be positive integers.
    First, $T_{h,d}(G)$ contains $K_1 \oplus G$ as a subgraph, and so $\rtd_2(T_{h,d}(G)) \geq \rtd_2(G)+1$ by \Cref{lemma:add_apex_increases_rtd2} and \Cref{lemma:rtd2_minor}.
    We now prove $\rtd_2(T_{h,d}(G)) \leq \rtd_2(G)+1$.
    Let $\mathcal{B}$ be the family of all the blocks of $T_{h,d}(G)$.
    For every $B \in \mathcal{B}$, $B$ is isomorphic to $K_1 \oplus G = T_{1,d}(G)$, and we denote by $s_B$ the root of $B$.
    By the definition of $T_{h,d}(G)$, one can inductively construct a tree $T$ with $V(T) = \mathcal{B}$
    such that, for every $B \in \mathcal{B}$ which is not the root,
    if $B'$ is the parent of $B$ in $T$,
    then $s_B$ belongs to $V(B')$ .
    By induction on $|V(T)|$, applying \ref{rtd2:item:more-blocks} to a leaf of $T$, we conclude that $\rtd_2(T_{h,d}(G)) \leq \rtd_2(K_1 \oplus G) \leq \rtd_2(G)+1$.
\end{proof}

\begin{lemma} \label{lemma:Thd_universal_for_rtd}
    For every nonnull graph $G$, 
    there exists a graph $H$ such that
    \begin{enumerateOurAlph}
        \item $\rtd_2(H) \leq \rtd_2(G)-1$ and 
        \item $G \subseteq T_{h,d}(H)$ for some positive integers $h$ and $d$.
    \end{enumerateOurAlph}
\end{lemma}

\begin{proof}
    If $G$ has no edges, then we set $H = K_1$ and the assertion follows, thus assume that $G$ has at least one edge.
    Let $t = \rtd_2(G)-1$.
    By \Cref{lemma:rtd2_and_Grt}, there exists a nonnegative integer $r$ such that $G \subseteq G_{r,t}$.
    Hence, it is enough to show that there exists a graph $H$ with $\rtd_2(H) \leq t$ and $G_{r,t} \subseteq T_{h,d}(H)$ for some positive integers $h,d$.
    Let $\mathcal{B}$ be the family of all the blocks of $G_{r,t}$.
    Recall that $G_{r,t} = L_{\binom{r+t}{t}}(G_{r-1,t}, K_1 \oplus G_{r,t-1}, u)$ if $r,t > 0$. Thus, 
    each $B\in\mathcal{B}$ is either isomorphic to $K_1\oplus G_{r,t-1}$, or is a block of $G_{r-1,t}$.
    Thus, for each $B\in\mathcal{B}$ there is a nonnegative integer $r'$ with $r' \leq r$ such that every $B$ is isomorphic to $K_1 \oplus G_{r',t-1}$. 
    This implies that for each $B\in\mathcal{B}$ we have $B \subseteq K_1 \oplus G_{r,t-1}$.

    We denote by $s_B$ the vertex of $B$ corresponding to $K_1$.
    By the definition of $G_{r,t}$, one can inductively construct a rooted tree $T$
    with $V(T) = \mathcal{B}$ such that for every $B \in \mathcal{B}$ which is not the root,
    the vertex $s_B$ belongs to $V(B')$ if $B'$ is the parent of $B$ in $T$ (in the same way as in~\Cref{obs:rtd_and_Thd}). 
    Moreover, $B-s_B$ is isomorphic to $G_{r',t-1}$ for some integer $r'$ with $r' \leq r$.
    Let $S$ be the root of $T$.
    Let $h$ be the vertex-height of $T$, and let $d$ be the maximum number of children of a vertex in $T$.
    We claim that $H = G_{r,t-1}$ satisfies the required conditions.
    The first assertion is clear by \Cref{lemma:rtd2_and_Grt}, thus, it suffices to prove that $G_{r,t} \subseteq T_{h,d}(G_{r,t-1})$.
    To this end, we show by induction the following property.
    For every subtree $T'$ of $T$ rooted in $S$ of vertex-height $h'$, we have $\bigcup_{B \in V(T')} V(B) \subseteq T_{h',d}(G_{r,t-1})$.

    When $h'=1$, $|V(T')| = 1$ and so $G\left[\bigcup_{B \in V(T')} V(B)\right] \subseteq K_1 \oplus G_{r,t-1} = T_{1,d}(G_{r,t-1})$.
    Now suppose that $h'>1$.
    Let $L$ be the set of all vertices of $T'$ at distance $h'-1$ from $S$.
    By induction hypothesis and because $T'-L$ has vertex-height $h'-1$, we have $\bigcup_{B \in V(T'-L)} V(B) \subseteq T_{h'-1,d}(G_{r,t-1})$.
    Then, $G\left[\bigcup_{B \in V(T')} V(B)\right]$ is obtained from $G\left[\bigcup_{B \in V(T'-L)} V(B)\right]$ by gluing on every vertex at most $d$ blocks isomorphic to a subgraph of $K_1 \oplus G_{r,t-1}$.
    Hence $G\left[\bigcup_{B \in V(T')} V(B)\right] \subseteq T_{h',d}(G_{r,t-1})$.
    Applying the above to $T'=T$ ends the proof.
\end{proof}

Having a model of $T_{h,d}(X)$ in some graph $G$, sometimes it will be handy to insist that a given vertex $u \in V(G)$ is in the branch set of the root of $T_{h,d}(X)$.
To this end, we introduce another auxiliary construction and we prove \Cref{lemma:root_a_model_of_T'hd} below.

\phantomsection\label{def:T'hd}
For every graph $X$ and for all positive integers $d,h$, let $T'_{h,d}(X)$ be the result of taking two disjoint copies of $T_{h,d}(X)$ and identifying their roots -- we call this new vertex the root of $T'_{h,d}(X)$.

\begin{lemma}\label{lemma:root_a_model_of_T'hd}
    Let $X$ be a graph and let $h,d$ be positive integers.
    Let $G$ be a connected graph and let $u \in V(G)$.
    If $\big(A_x \mid x \in V(T'_{h,d}(X))\big)$ is a model of $T'_{h,d}(X)$ in $G$,
    then there exists a model $\big(B_x \mid x \in V(T_{h,d}(X))\big)$ of $T_{h,d}(X)$ in $G$ such that
    \begin{enumerateOurAlph}
        \item $u \in B_s$ where $s$ is the root of $T_{h,d}(X)$ and
        \item for every $x \in V(T_{h,d}(X))$, there exists $y \in V(T'_{h,d}(X))$ such that $A_y \subseteq B_x$.
    \end{enumerateOurAlph}
\end{lemma}

\begin{proof}
    Observe that $T'_{h,d}(X)$ has a separation $(H_1,H_2)$ such that
    $V(H_1) \cap V(H_2) = \{s'\}$ where $s'$ is the root of $T'_{h,d}(X)$,
    and $H_i$ is isomorphic to $T_{h,d}(X)$ for each $i \in \{1,2\}$.

    Let $A$ be the union of all branch sets $A_x$ for $x \in V(T'_{h,d}(X))$.
    Fix a path $P$ in $G$ from $u$ to any vertex of $A$ with no internal vertices in $A$.
    Suppose that the endpoint of $P$ in $A$ is contained in $A_x$ for some $x_0 \in V(T'_{h,d}(X))$
    Let $\{i,j\} = \{1,2\}$ be such that $x_0 \in V(H_i)$.
    Let $B_{s} = V(P) \cup \bigcup_{x \in V(H_i)} A_x$ and
    $B_y = A_y$ for every $y \in V(H_{j}) \setminus \{s\}$.
    Then $\big(B_y \mid y \in V(H_{j})\big)$ is a model of $T_{h,d}(X)$ in $G$
    satisfying the conclusion of the lemma.
\end{proof}

\section{The base case: Graphs with no \texorpdfstring{$\mathcal{F}$}{F}-rich model of a given tree} \label{sec:base}

The proof of \Cref{thm:rooted_main} is by induction.
The technical statement of the induction is stated in \Cref{thm:main}.
In this section, we provide the base case for the induction, that is, the case where $X$ is a forest.
It turns out that methods used in the so-called theory of product structure of graphs are useful in the study of weak coloring numbers.
For instance, the main result of~\cite{DHHJLMMRW24} is actually a product structure result and the weak coloring numbers bound follows from a slight adjustment of the argument.
In the material of this section, we are strongly inspired by another product structure paper by Dujmović, Hickingbotham, Joret, Micek, Morin, and Wood~\cite{Dujmovi2023}.

For all positive integers $h$ and $d$, we denote by $F_{h,d}$ the (rooted) complete $d$-ary tree of vertex-height $h$.
In particular, $F_{2,d}$ is the star with $d$ leaves.
Note that for every tree $X$, there exist positive integers  $h,d$ such that $X \subseteq F_{h,d}$.

We start with the case where $X$ is a star.
Let $\delta$ be the function given by \Cref{lemma:E-P_in_Kk_minor_free_graphs}.

In the following proof, we need the notion of path decompositions.
General tree decompositions are discussed in detail in~\Cref{sec:bounded_tw}.
Let $G$ be a graph.
A sequence of subsets $(W_0,\dots,W_\ell)$ of $V(G)$ is a \emph{path decomposition} of $G$ if 
\begin{enumerate}
    \item for every $u \in V(G)$, the set $\{i \in \{0,\dots,\ell\} : u \in W_i\}$ is a nonempty interval, and
    \item for every edge $uv \in E(G)$, there exists $i \in \{0,\dots,\ell\}$ with $u,v \in W_i$.
\end{enumerate}

\begin{lemma}\label{lemma:excluding_a_star_patch}
    Let $k,d$ be positive integers.
    Let $G$ be a connected $K_k$-minor-free graph, and
    let $\mathcal{F}$ be a family of connected subgraphs of $G$
    such that
    $G$ has no $\mathcal{F}$-rich model of $F_{2,d}$.
    For every nonempty $U \subset V(G)$ such that $G[U]$ is connected, 
    there is a path decomposition $(W_0, \dots, W_\ell)$ of $G$ with $\ell \geq 1$
    and sets $R_2, \dots, R_\ell \subseteq V(G)$
    such that
    for $S = U \cup \bigcup_{i\in\{2, \dots,\ell\}} (W_{i-1} \cap W_i)$,
    \begin{enumerateOurAlph}
        \item $W_0 = U$; \label{item:excluding_a_star_patch_i}
        \item $V(F) \cap S \neq \emptyset$ for every $F \in \mathcal{F}$; \label{item:excluding_a_star_patch_ii}
        \item $G[S]$ is connected; \label{item:excluding_a_star_patch_iii}
        \item $G[R_i]$ is connected for every $i \in \{2, \dots, \ell\}$; \label{item:excluding_a_star_patch_iv}
        \item $W_{i-1} \cap W_i \subseteq R_{i} \subseteq \bigcup_{j \in \{0, \dots, i-1\}} W_j$ for every $i \in \{2, \dots, \ell\}$; \label{item:excluding_a_star_patch_v}
        \item $W_{i}$ and $W_{i+2}$ are disjoint for every $i \in \{0, \dots, \ell-2\}$; and \label{item:excluding_a_star_patch_vi}
        \item $\wcol_r(G,R_i) \leq (\delta(k,d+1)+3) \cdot r$ for every $i \in \{2, \dots, \ell\}$ and for every positive integer~$r$. \label{item:excluding_a_star_patch_vii}
    \end{enumerateOurAlph}
\end{lemma}

    \begin{figure}[tp]
        \centering 
        \includegraphics[scale=1]{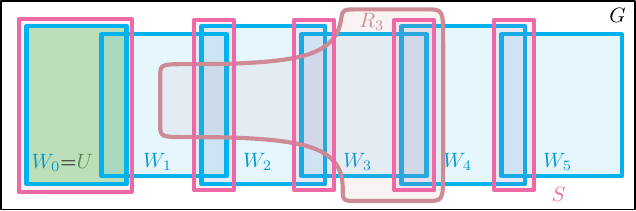} 
        \caption{
        We illustrate the statement of \Cref{lemma:excluding_a_star_patch}.
        The green area is a given set $U$, which should become $W_0$.
        The set $R_3$ has to satisfy $W_3 \cap W_4 \subset R_3 \subset W_0 \cup W_1 \cup W_2 \cup W_3 \cup W_4$ and we want $\wcol_r(G,R_3)$ to be low.
        We do not depict all the sets $R_i$ for readability.
        }
        \label{fig:excluding-star-statement}
    \end{figure}

The statement of the lemma is visualized in Figure~\ref{fig:excluding-star-statement}.

\begin{proof}
    In the proof we define many objects, they are depicted in \Cref{fig:excluding-star}.
    We proceed by induction on $|V(G)|-|U|$.
    Let $U \subseteq V(G)$ be nonempty such that $G[U]$ is connected.
    If $\mathcal{F}\vert_{G-U} = \emptyset$, then it suffices to take $W_0=W_1=U$, $\ell=1$.
    In particular, this is the case for $U = V(G)$.
    Therefore, assume $|U|<|V(G)|$ and $\mathcal{F}\vert_{G-U} \neq \emptyset$.
    Let $\mathcal{F}_0$ be the family of all the connected subgraphs $A$ of $G-U$
    such that $A$ contains a member of $\mathcal{F}$ and $V(A) \cap N_G(U) \neq \emptyset$.
    We argue that $\mathcal{F}_0 \neq \emptyset$.
    Since $\mathcal{F}\vert_{G-U} \neq \emptyset$,
    there is a component $C$ of $G-U$ containing a member of $\mathcal{F}$.
    Since $G$ is connected, $V(C) \cap N_G(U) \neq \emptyset$ and so $C \in \mathcal{F}_0$.

    Observe that any collection of $d+1$ pairwise disjoint $A_1, \dots, A_{d+1} \in \mathcal{F}_0$ yields
    an $\mathcal{F}$-rich model of $F_{2,d}$. Indeed, it suffices to take $U \cup A_{d+1}$ as the branch set 
    corresponding to the root of $F_{2,d}$ and $A_1, \dots, A_d$ as the branch sets of the remaining $d$ vertices
    of $F_{2,d}$.
    Therefore, there are no $d+1$ pairwise disjoint members of $\mathcal{F}_0$,
    and thus, by Lemma~\ref{lemma:E-P_in_Kk_minor_free_graphs}
    applied to $G$ and $\mathcal{F}_0$, there exists a set $S_0 \subseteq V(G)$ such that
    \begin{enumerate}[label={\normalfont\ref{lemma:E-P_in_Kk_minor_free_graphs}.(\makebox[\mywidth]{\alph*})}]
        \item $V(F) \cap S_0 \neq \emptyset$ for every $F \in \mathcal{F}_0$; \label{lem:EP:item:hit'_patch}
        \item $G[S_0]$ is connected; \label{lem:EP:item:connected'_patch}
        \item $\wcol_r(G,S_0) \leq \delta(k,d+1) \cdot r$ for every positive integer $r$. \label{lem:EP:item:wcol'_patch} 
    \end{enumerate}
    Since $\mathcal{F}_0 \neq \emptyset$, we have $S_0 \setminus U \neq \emptyset$.
    Let $Q$ be a $U$-$S_0$ geodesic in $G$ (possibly just a one-vertex path), and let $S_1 = S_0 \cup V(Q)$.
    Note that by~\ref{lem:EP:item:connected'_patch}, $G[S_1]$ is connected. 

    Let $\mathcal{C}_0$ be the family of all the components $C$ of $G-U-S_1$ such that $N_G(U) \cap V(C) = \emptyset$.
    Let $U' = V(G) - \bigcup_{C \in \mathcal{C}_0} V(C)$.
    Observe that $|U'| > |U|$ since $S_0 \setminus U \neq \emptyset$ and $U'$ contains $U \cup S_0$.
    Let $\mathcal{F}' = \{F \in \mathcal{F} \mid V(F) \cap U' = \emptyset\}$.
    By induction hypothesis applied to $G$, $\mathcal{F}'$ and $U'$,
    there is a path decomposition $(W'_0, \dots, W'_{\ell'})$ of $G$
    and sets $R'_2, \dots, R'_{\ell'} \subseteq V(G)$
    such that
    for $S' = U' \cup \bigcup_{i\in\{2, \dots,\ell'\}} (W'_{i-1} \cap W'_i)$,
    \begin{enumerateOurAlphPrim}
        \item $W'_0 = U'$; \label{item:excluding_a_star_patch_i'}
        \item $V(F) \cap S' \neq \emptyset$ for every $F \in \mathcal{F}'$; \label{item:excluding_a_star_patch_ii'}
        \item $G[S']$ is connected; \label{item:excluding_a_star_patch_iii'}
        \item $G[R'_i]$ is connected for every $i \in \{2, \dots, \ell'\}$; \label{item:excluding_a_star_patch_iv'}
        \item $W'_{i-1} \cap W'_i \subseteq R'_i \subseteq \bigcup_{j \in \{0, \dots, i-1\}} W'_j$ for every $i \in \{2, \dots, \ell'\}$; \label{item:excluding_a_star_patch_v'}
        \item $W'_{i}$ and $W'_{i+2}$ are disjoint for every $i \in \{0, \dots, \ell'-2\}$; and \label{item:excluding_a_star_patch_vi'}
        \item $\wcol_r(G,R'_i) \leq (\delta(k,d+1)+3) \cdot r$ for every $i \in \{2, \dots, \ell\}$ and for every positive integer~$r$. \label{item:excluding_a_star_patch_vii'}
    \end{enumerateOurAlphPrim}

    \begin{figure}[tp]
        \centering 
        \includegraphics[scale=1]{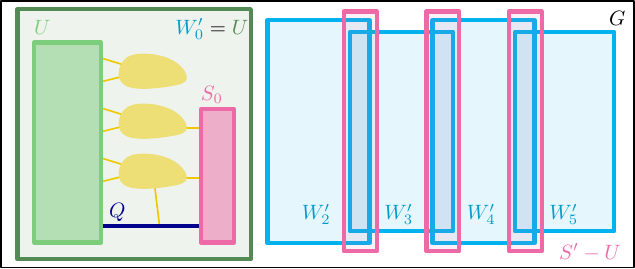} 
        \caption{
        An illustration of the considered objects in \Cref{lemma:excluding_a_star_patch}.
        Note that $S_0$ may intersect $U$.
        }
        \label{fig:excluding-star}
    \end{figure} 
    
    Let $\ell = \ell'+1$, $W_0 = U$, $W_1 = U'$, $W_2 = (W'_1 \setminus U') \cup (S_1 \setminus U)$, $W_i = W'_{i-1}$ for every $i \in \{3, \dots, \ell\}$,
    $R_2 = S_1$, and $R_i = R'_{i-1}$ for every $i \in \{3, \dots, \ell\}$.
    Note that~\ref{item:excluding_a_star_patch_i} holds by construction.
    We claim that $(W_0, \dots, W_\ell)$ is a path decomposition of $G$
    and \ref{item:excluding_a_star_patch_ii}-\ref{item:excluding_a_star_patch_vii} hold, which completes the proof of the lemma.

    Let $u \in V(G)$.
    We claim that $I = \{i \in \{0, \dots, \ell\} \mid u \in W_i\}$ is an interval.
    Since $(W'_0, \dots, W'_{\ell'})$ is a path decomposition of $G$,
    $I' = \{i \in \{0, \dots, \ell'\} \mid u \in W'_i\}$ is an interval.
    If $u \not \in U' = W'_0$, then $I =  \{i \in \{2, \dots, \ell\} \mid u \in W'_{i-1}\} = \{i+1 \mid i \in I'\}$,
    which is an interval too.
    Now suppose that $u \in U'$, and so $0 \in I'$.
    If $u \not\in S_1 \setminus U$, then $u \not \in W_2$ and $u \not\in W'_i$ for every $i \in \{2, \dots, \ell'\}$ by \ref{item:excluding_a_star_patch_i'} and \ref{item:excluding_a_star_patch_vi'}.
    Hence
    $I = \{0,1\}$ if $u \in U$, and $I = \{1\}$ otherwise,
    which is an interval in both cases.
    If $u \in S_1 \setminus U$, then $u \not\in U = W_0$,
    and so $I = \{1\} \cup \{i+1 \mid i \in I' \setminus \{0\}\}$, which is an interval.
    This proves that $I$ is an interval.

    Let $uv$ be an edge of $G$.
    We claim that there exists $i \in \{0, \dots, \ell\}$ such that $u,v \in W_i$.
    If there exists $i'\in \{2,\dots,\ell'\}$ such that $u,v \in W'_{i'}$, then $u,v \in W'_{i'} = W_{i'-1}$ and we are done.
    Now suppose that $u$ and $v$ are not both in $W'_i$ for every $i \in \{2, \dots, \ell\}$.
    Since $(W'_0, \dots, W'_{\ell'})$ is a path decomposition of $G$,
    there exists $i' \in \{0, 1\}$ such that $u,v \in W'_{i'}$.
    If $i' = 0$, then $u,v \in W'_0 = U' = W_1$.
    Now suppose that $u$ and $v$ are not both in $W'_0=U'$, and so, in particular, $i' = 1$.
    Without loss of generality assume that $v \notin U'$.
    It follows that $v \in W_1' - U' \subset W_2$.
    Let $C$ be the component of $v$ in $G-U-S_1$.
    Since $v \not\in U'$, $C$ belongs to $\mathcal{C}_0$, and so $N_G(V(C)) \cap U = \emptyset$.
    It follows that $N_G(V(C)) \cap U = \emptyset$, and so, $u \in S_1 \setminus U$.
    Therefore, $u \in W_2$, which concludes the claim.
    Furthermore, we obtained that $(W_0, \dots, W_\ell)$ is a path decomposition of $G$.

    We now prove \ref{item:excluding_a_star_patch_ii}.
    Consider $F \in \mathcal{F}$.
    If $F$ intersects $U$, then $V(F) \cap S \neq \emptyset$ since $U \subseteq S$.
    If $F$ intersects $S_1 \setminus U$, then $F$ intersects $W_1 \cap W_2 \subseteq S$.
    Now suppose that $F$ is disjoint from $U \cup S_1$.
    Let $C$ be the component containing $F$ in $G-U-S_1$.
    Since $C$ is disjoint from $S_0 \subseteq S_1$, by \ref{lem:EP:item:hit'_patch}, $C$ is not a member of $\mathcal{F}_0$.
    This implies that $N_G(U) \cap V(C) = \emptyset$, and thus, $C \in \mathcal{C}_0$.
    In particular, $F$ is disjoint from $U'$, and so, $F \in \mathcal{F}'$.
    By \ref{item:excluding_a_star_patch_ii'}, $V(F) \cap S' \neq \emptyset$,
    hence, there exists $i \in \{2, \dots, \ell'\}$ such that $V(F)$ intersects $W'_{i-1} \cap W'_i$.
    It follows that $W'_{i-1} \cap W'_i = W_{i} \cap W_{i+1}$ and so $V(F) \cap S \ne \emptyset$.
    This proves \ref{item:excluding_a_star_patch_ii}.

Let us pause to underline a simple observation that follows directly from the construction,~\ref{item:excluding_a_star_patch_i'}, and~\ref{item:excluding_a_star_patch_vi'}:
    \begin{enumerate}[label=($\star$)]
        \setcounter{enumi}{6}
        \item for every $i \in \{3,\dots,\ell\}$, we have $W_{i-1} \cap W_i = W_{i-2}' \cap W_{i-1}'$. \label{star:intersections}
    \end{enumerate}

    By \ref{item:excluding_a_star_patch_iii'}, $S' = U' \cup \bigcup_{i \in \{2, \dots, \ell'\}} (W'_{i-1} \cap W'_i)$ induces a connected subgraph of $G$.
    In particular, every component of $G[S']-U'$ has a neighbor in $N_G(V(G)\setminus U') \subseteq S_1$.
    Since $G[U \cup S_1]$ is connected, it follows that $(S' \setminus U') \cup U \cup S_1$ induces a connected subgraph of $G$.
    However, $S = (S' \setminus U') \cup U \cup S_1$ by~\ref{star:intersections}, which yields~\ref{item:excluding_a_star_patch_iii}.

    For every $i \in \{3, \dots, \ell\}$,
    $R_i = R'_{i-1}$ induces a connected subgraph of $G$ by \ref{item:excluding_a_star_patch_iv'}, and $R_2 = S_1$ induces a connected subgraph of $G$ by definition, hence~\ref{item:excluding_a_star_patch_iv} follows.

    For the proof of~\ref{item:excluding_a_star_patch_v}, first, observe that by construction,
    $\bigcup_{j \in \{0, \dots, i-1\}} W_j = \bigcup_{j \in \{0, \dots, i-2\}} W'_j$ for every $i \in \{2, \dots, \ell\}$.
    In particular, it follows that
    $R_i = R'_{i-1} \subseteq \bigcup_{j \in \{0, \dots, i-2\}} W'_j =\bigcup_{j \in \{0, \dots, i-1\}} W_j$ for every $i \in \{3,\dots,\ell\}$ by \ref{item:excluding_a_star_patch_v'}.
    Moreover, $R_2 = S_1 \subseteq W_1$.
    It remains to show that $W_{i-1} \cap W_i \subseteq R_i$ for every $i \in \{2, \dots, \ell'\}$.
    For $i=2$,
    $W_1 \cap W_2 = S_1 \setminus U \subseteq R_2$.
    For $i \in \{3, \dots, \ell\}$, $W_{i-1} \cap W_i = W'_{i-2} \cap W'_{i-1} \subseteq R'_{i-1} = R'_i$ by~\ref{star:intersections} and~\ref{item:excluding_a_star_patch_v'}.
    This gives \ref{item:excluding_a_star_patch_v}.

    For every $i \in \{3, \dots, \ell-2\}$,
    $W_i \cap W_{i+2} = W'_{i-1} \cap W'_{i+1} = \emptyset$ by \ref{item:excluding_a_star_patch_vi'}.
    Moreover, $W'_3=W_4$ is disjoint from $S_1 \setminus U \subseteq W'_0$ by \ref{item:excluding_a_star_patch_vi'}.
    Hence $W_2 \cap W_4 = W'_1 \cap W'_3 = \emptyset$ by \ref{item:excluding_a_star_patch_vi'}.
    Similarly, $W'_2=W_3$ is disjoint from $U' = W'_0$ by \ref{item:excluding_a_star_patch_vi'}.
    Hence $W_1 \cap W_3 = W'_0 \cap W'_2 = \emptyset$.
    Finally, $W_0 \cap W_2 = \emptyset$ by construction, and so, \ref{item:excluding_a_star_patch_vi} holds.

    It remains to show \ref{item:excluding_a_star_patch_vii}.
    First, for every $i \in \{3, \dots, \ell\}$,
    $R_i = R'_{i-1}$ and so $\wcol_r(G,R_i) \leq (\delta(k,d+1)+3)\cdot r$ for every positive integer $r$ by \ref{item:excluding_a_star_patch_vii'}.
    Moreover, $R_2 = S_1 = V(Q) \cup S_0$.
    Hence
    \[
    \wcol_r(G,R_2) \leq \delta(k,d+1) \cdot r + (2r+1) \leq (\delta(k,d+1)+3) \cdot r
    \]
    for every positive integer $r$, using \ref{lem:EP:item:wcol'_patch} and \Cref{obs:geodesics}.
    This shows that \ref{item:excluding_a_star_patch_vii} holds, which concludes the proof of the lemma.
\end{proof}

    \begin{figure}[tp]
        \centering 
        \includegraphics[scale=1]{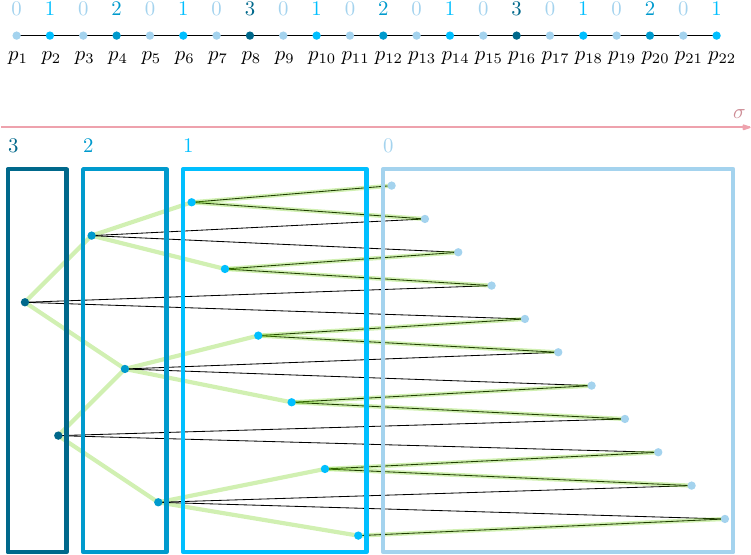} 
        \caption{
            Consider the path $P = p_1\cdots p_{22}$.
            Let $r = 7$ and $s = 3 = \left\lceil \log r \right\rceil$.
            We claim that $\wcol_r(P) \leq 2+s$.
            We mark every eighth vertex with $3$ (this corresponds to the set $I_3$ in the proof of~\Cref{lemma:excluding_a_star_wcol}).
            Then, we mark each fourth unmarked vertex with $2$ (the set $I_2 - I_3$), every second unmarked vertex with $1$ (the set $I_1 - I_2$), and finally all remaining vertices with $0$ (the set $I_0 - I_1$).
            Next, we construct an ordering on the vertices $p_1,\dots,p_{22}$.
            First, preserving the original ordering in the path, we put the vertices marked with $3$, then similarly the ones marked with $2$, with $1$, and with $0$.
            It suffices to argue that for every $u \in V(P)$, we have $|\WReach_r[P,\sigma,u]| \leq 5$.
            We build an auxiliary tree (depicted in green in the figure), where for every $i \in \{3,2,1\}$, we connect every vertex marked with $i$ with the closest vertices in the path marked with $i-1$.
            It is not hard to see that $\WReach_r[P,\sigma,u]$ consists only of the ancestors of $u$ in the auxiliary tree.
        }
        \label{fig:path}
    \end{figure} 

Now we are ready to prove Theorem~\ref{thm:main} for
graphs with no $\mathcal{F}$-rich model of a fixed star. 
This part of the argument follows ideas from the proof by \cite{JM22} that $\wcol_r(P)\leq 2+\lceil \log r \rceil$ for every path $P$ and every positive integer $r$ -- we explain the idea of this proof in \Cref{fig:path}.

\begin{lemma}\label{lemma:excluding_a_star_wcol}
    Let $k,d$ be positive integers.
    For every connected $K_k$-minor-free graph $G$, for every family $\mathcal{F}$ of connected subgraphs of $G$, if $G$ has no $\mathcal{F}$-rich model of $F_{2,d}$, then
    there is a set $S \subseteq V(G)$ such that
    \begin{enumerate}[label={\normalfont(\makebox[\mywidth]{\alph*})}]
        \item $V(F) \cap S \neq \emptyset$ for every $F \in \mathcal{F}$;\label{lemma:excluding_a_star_patch_wcol:assertion:hitting}
        \item $G[S]$ is connected;\label{lemma:excluding_a_star_patch_wcol:assertion:connected}
        \item $\wcol_r(G,S) \leq 5(\delta(k,d+1)+3) \cdot r \log r$ for every integer $r$ with $r \geq 2$. \label{lemma:excluding_a_star_patch_wcol:assertion:wcol}
    \end{enumerate}
\end{lemma}

\begin{proof}   
    Let $G$ be a connected $K_k$-minor-free graph, let $\mathcal{F}$ be a family of connected subgraphs of $G$, 
    and suppose that $G$ has no $\mathcal{F}$-rich model of $F_{2,d}$.
    Let $r$ be an integer with $r \geq 2$. 
    Let $U$ be an arbitrary singleton of a vertex in $G$.
    \Cref{lemma:excluding_a_star_patch} applied to $G$, $\mathcal{F}$, and $U$ gives
    a path decomposition $(W_0, \dots, W_\ell)$
    and sets $R_2, \dots, R_\ell \subseteq V(G)$
    such that for $S' = U \cup \bigcup_{i \in \{2, \dots, \ell\}} (W_{i-1} \cap W_i)$,
    \begin{enumerate}[label={\normalfont\ref{lemma:excluding_a_star_patch}.(\makebox[\mywidth]{\alph*})}]
        \item $W_0 = U$; \label{item:excluding_a_star_recall_i}
        \item $V(F) \cap S' \neq \emptyset$ for every $F \in \mathcal{F}$; \label{item:excluding_a_star_recall_ii}
        \item $G[S']$ is connected; \label{item:excluding_a_star_recall_iii}
        \item $G[R_i]$ is connected for every $i \in \{1, \dots, \ell\}$; \label{item:excluding_a_star_recall_iv}
        \item $W_{i-1} \cap W_i \subseteq R_{i} \subseteq \bigcup_{j \in \{0, \dots, i-1\}} W_j$ for every $i \in \{2, \dots, \ell\}$; \label{item:excluding_a_star_recall_v}
        \item $W_{i}$ and $W_{i+2}$ are disjoint for every $i \in \{0, \dots, \ell-2\}$; and \label{item:excluding_a_star_recall_vi}
        \item $\wcol_r(G,R_i) \leq (\delta(k,d+1)+3) \cdot r$ for every $i \in \{2, \dots, \ell\}$. \label{item:excluding_a_star_recall_vii}
    \end{enumerate}

    For convenience, we set $R_1 = U$.

    Let $s = \lceil \log(r+1) \rceil$.
    For every $i \in \{0, \dots, s\}$, let $I_i = \{i \in \{1, \dots, \ell\} \mid j = 0 \mod 2^i\}$.
    We construct recursively families $\set{R'_j}_{j \in \{1,\dots,\ell\}}$ and 
    $\set{S_i}_{i \in \{0,\dots,s\}}$ of subsets of $V(G)$ and a family $\{\sigma_j\}_{j \in \{1,\dots,\ell\}}$ such that 
    $\sigma_j$ is an ordering of $R_j'$ for every $j \in \{1, \dots, \ell\}$.
    For every $j \in I_s$, let 
    \[ 
        R_j' = R_j \setminus \textstyle\bigcup_{a \in  \{0,\dots , j-2^s-1\}} W_{a}
    \]
    and let $S_s = \bigcup_{ j \in I_s} R'_j$.
    Let $j \in I_s$.
    If $j < 2 \cdot 2^s$, then $j = 2^s$ and $R_j' = R_j$, and so by~\ref{item:excluding_a_star_recall_vii}, $\wcol_r(G,R_j') \leq (\delta(k,d+1)+3) \cdot r$.
    Now assume that $j \geq 2 \cdot 2^s$.
    Since $(W_0, \dots, W_\ell)$ is a path decomposition of $G$, $W_{j-2^s-1} \cap W_{j-2^s}$ separates $\bigcup_{a \in \{0, \dots, j-2^s-1\}} W_a$ and $\bigcup_{a \in \{j-2^s, \dots, \ell\}} W_a$ in $G$.
    Since $W_{j-2^s-1} \cap W_{j-2^s} \subseteq R_{j-2^s}'$ (by \ref{item:excluding_a_star_recall_v}), by \Cref{obs:wcol_components}, we obtain
    \begin{align*}
        \wcol_r\left(G-\textstyle\bigcup_{a \in \{1, \dots, j-2^s\} \cap I_s} R'_a, R'_j\right) 
        &= \wcol_r\left(G - (W_{j-2^s-1} \cap W_{j-2^s}), R'_j\right) \\
        &= \wcol_r\left(G - \textstyle\bigcup_{a \in \{0, \dots ,j-2^s-1\}} W_a, R'_j\right).
    \end{align*}
    Finally,
    \begin{align*}
        \wcol_r\left(G - \textstyle\bigcup_{a \in \{0, \dots ,j-2^s-1\}} W_a, R'_j\right)
        &\leq \wcol_r(G,R_j) && \textrm{by \Cref{obs:monotone}} \\
        &\leq (\delta(k,d+1)+3)\cdot r && \textrm{by \ref{item:excluding_a_star_recall_vii}}.
    \end{align*}
    Let $\sigma_j$ be an ordering of $R'_j$ such that
    \[
        \wcol_r\left(G - \textstyle\bigcup_{a \in \{1, \dots, j-2^s\} \cap I_s} R'_a, R'_j, \sigma_j\right) \leq (\delta(k,d+1)+3)\cdot r.
    \]

    Next, let $i \in \{0, \dots, s-1\}$ and assume that $S_{i+1}$ is defined.
    Now, for every $j \in I_i \setminus I_{i+1}$, let
    \[
    R'_j = \left(R_j \setminus \textstyle\bigcup_{a \in \{0, \dots, j-2^i-1\}} W_a \right) \setminus S_{i+1},
    \]
    and let $S_i = \bigcup_{j \in I_i} R'_j$.
    Note that $S_{i+1} \subset S_i$.
    Also note that for every $j \in I_i$, $W_{j-1} \cap W_j \subset R_j'$ by~\ref{item:excluding_a_star_recall_vi}.
    Let $j \in I_i \setminus I_{i+1}$.
    We have $j-2^i \in I_{i+1}$, and therefore,
    $W_{j-2^i-1} \cap W_{j-2^i} \subseteq R_{j-2^i}' \subset S_{i+1}$.
    Since $(W_0 ,\dots, W_\ell)$ is a path decomposition of~$G$,
    $W_{j-2^i-1} \cap W_{j-2^i}$ separates $\bigcup_{a \in \{0, \dots, j-2^i-1\}} W_a$ and $\bigcup_{a \in \{j-2^i, \dots, \ell\}} W_a$ in $G$.
    It follows by \Cref{obs:wcol_components} that 
    \begin{align*}
    \wcol_r(G-S_{i+1}, R'_j)
    &= \wcol_r\left(G-(W_{j-2^i-1} \cap W_{j-2^i}), R'_j\right) \\
    &= \wcol_r\left(G-\textstyle\bigcup_{a \in \{0, \dots, j-2^i-1\}} W_a,R'_j\right).
    \end{align*}
    Furthermore,
    \begin{align*}
        \wcol_r\left(G-\textstyle\bigcup_{a \in \{0, \dots, j-2^i-1\}} W_a,R'_j\right) 
        &\leq \wcol_r(G,R_j) && \textrm{by \Cref{obs:monotone}} \\
        &\leq (\delta(k,d+1)+3) \cdot r && \textrm{by \ref{item:excluding_a_star_recall_vii}}.
    \end{align*}
    Let $\sigma_j$ be an ordering of $R'_j$ such that
    \[
    \wcol_r(G-S_{i+1}, R'_j, \sigma_j) \leq (\delta(k,d+1)+3) \cdot r.
    \]

    We define $S = S_0$.
    Now, it suffices to show that~\ref{lemma:excluding_a_star_patch_wcol:assertion:hitting}-\ref{lemma:excluding_a_star_patch_wcol:assertion:wcol} hold.
    Since $S' \subseteq S$, \ref{lemma:excluding_a_star_patch_wcol:assertion:hitting} holds by~\ref{item:excluding_a_star_recall_ii}.
    
    Recall that $G[S']$ is connected by \ref{item:excluding_a_star_recall_iii}.
    Next, let $C$ be a component of $G[R_j']$ for some fixed $j \in I_s$.
    If $V(C) \cap (W_{j-1} \cap W_j) \neq \emptyset$, then $V(C) \cap S' \neq \emptyset$, and so, $G[S' \cup V(C)]$ is connected.
    Thus, assume that $V(C) \cap (W_{j-1} \cap W_j) = \emptyset$.
    However, $W_{j-1} \cap W_j \subset R_j'$, hence, $C$ has a neighbor in $W_{j-2^s-1}$, in particular, in $W_{j-2^s-1} \cap W_{j-2^s} \subset S'$.
    Hence, again $G[S' \cup V(C)]$ is connected.
    In particular, we have just proved that $G[S' \cup S_s]$ is connected.
    Next, suppose that $G[S' \cup S_{i+1}]$ is connected for some $i \in \{0,\dots,s-1\}$.
    Let $C$ be a component of $G[R_j']$ for some fixed $j \in I_i$.
    If $V(C) \cap (W_{j-1} \cap W_j) \neq \emptyset$, then $V(C) \cap S' \neq \emptyset$, and so, $G[S' \cup V(C)]$ is connected.
    Thus, assume that $V(C) \cap (W_{j-1} \cap W_j) = \emptyset$.
    However, $W_{j-1} \cap W_j \subset R_j'$, hence, $C$ has a neighbor in $W_{j-2^i-1} \cup S_{i+1}$, in particular, in $(W_{j-2^i-1} \cap W_{j-2^i}) \cup S_{i+1} \subset S' \cup S_{i+1}$.
    Hence, $G[S' \cup S_{i+1} \cup V(C)]$ is connected.
    Finally, $G[S'\cup S_0] = G[S]$ is connected, which yields~\ref{lemma:excluding_a_star_patch_wcol:assertion:connected}.
    
    The sets $\{R'_j\}_{j \in \{2, \dots, \ell\}}$ are pairwise disjoint, and they partition $S$.
    Let $\sigma$ be an ordering of $S$ such that
    \begin{enumerate}
        \item $\sigma$ extends $\sigma_j$, for every $j \in \{1, \dots, \ell\}$;
        \item for every $j,j' \in I_s$ with $j<j'$, for all $u \in R'_j$ and $v \in R'_{j'}$, $u<_\sigma v$; and
        \item for every $i \in \{0, \dots, s-1\}$, for all $u \in S_{i+1}$ and $v \in S_i \setminus S_{i+1}$, $u<_\sigma v$.
    \end{enumerate}
    Note the similarity of this ordering to the one described in~\Cref{fig:path}.  
    For convenience, let $R'_0 = \emptyset$ and $W_j=R'_j=\emptyset$ for every integer $j$ with $j > \ell$.

    We now show \ref{lemma:excluding_a_star_patch_wcol:assertion:wcol}.
    Let $u \in V(G)$.
    We will show that $|\WReach_r[G,S,\sigma,u]| \leq 5(\delta(k,d+1)+3)\cdot r \log r$.
    Let $j_u \in \{0, \dots, \ell\}$ be minimum such that $u \in W_{j_u}$.

    We claim that
    \[
    |\WReach_r[G,S,\sigma,u] \cap S_s| \leq 2(\delta(k,d+1)+3) \cdot r.
    \]
    Let $\alpha = \max\{0\} \cup \{a \in I_s \mid a \leq j_u\}$,
    and let $\beta = \alpha + 2^s$.
    Thus, if $\beta \leq \ell$, then $\beta \in I_s$.
    Next, we argue that
    \[
    \WReach_r[G,S,\sigma,u] \cap S_s \subseteq R'_{\alpha} \cup R'_{\beta}.
    \]
    Suppose to the contrary that there is a vertex $v \in \WReach_r[G,S,\sigma,u] \cap S_s$ with
    $v \not\in R'_{\alpha} \cup R'_{\beta}$.
    Let $\gamma \in I_s \setminus \{\alpha, \beta\}$ be such that $v \in R'_\gamma$.
    Then either $\gamma < \alpha$, or $\gamma > \beta$.
    First assume that $\gamma < \alpha$.
    Since $R'_\gamma \subseteq \bigcup_{a \in \{0, \dots, \gamma-1\}} W_a$ and because $(W_0, \dots, W_\ell)$ is a path decomposition of $G$,
    every $u$-$v$ path in $G$ intersects $W_{a-1} \cap W_{a}$ for each $a \in \{\gamma, \dots, j_u\}$.
    Since $(W_{a-1} \cap W_a)_{a \in \{1, \dots, \ell\}}$ are pairwise disjoint,
    we deduce that $\dist_G(u,v) \geq j_u - \gamma \geq \alpha - \gamma \geq 2^s > r$,
    which contradicts the fact that $v \in \WReach_r[G,S,\sigma,u]$.
    Finally, assume $\gamma > \beta$.
    Note that $\gamma \leq \ell$ since $R'_\gamma \neq \emptyset$ as $v \in R'_\gamma$.
    Since $R'_\gamma \subseteq \bigcup_{a \in \{\gamma-2^s, \dots, \gamma-1\}} W_a$,
    and because $(W_0, \dots, W_\ell)$ is a path decomposition of $G$,
    every $u$-$v$ path in $G$ intersects $W_{\beta-1} \cap W_\beta$.
    However, for every $w \in W_{\beta-1} \cap W_\beta$, we have $w<_\sigma v$, thus, $v \not\in \WReach_r[G,S,\sigma,u]$, which is a contradiction.
    We obtain that $\WReach_r[G,S,\sigma,u] \cap S_s \subseteq R'_{\alpha} \cup R'_{\beta}$.

    For every $\epsilon \in \{\alpha,\beta\}$, by definition of $\sigma$, we have
    \begin{align*}
        \WReach_r[G,S,\sigma,u] \cap R'_\epsilon
    &\subseteq \WReach_r\left[G - \textstyle\bigcup_{a \in I_s \cap \{1,\dots,\alpha-1\}}  R'_a, R'_\epsilon, \sigma_\epsilon, u\right] \\
    &\subseteq \WReach_r\left[G - \textstyle\bigcup_{a \in \{1, \dots, j_u-2^s\} \cap I_s} R'_a, R'_\epsilon, \sigma_\epsilon, u\right]
    \end{align*}
    and therefore,
    \[
    |\WReach_r[G,S,\sigma,u] \cap R'_\epsilon| \leq (\delta(k,d+1)+3) \cdot r.
    \]
    In particular,
    \begin{align*}
    |\WReach_r[G,S,\sigma,u] \cap S_s| 
    &\leq |\WReach_r[G,S,\sigma,u] \cap (R'_\alpha \cup R'_\beta)| \\
    &\leq 2(\delta(k,d+1)+3) \cdot r.
    \end{align*}
    
    Next, let $i \in \{0, \dots, s-1\}$. We claim that
    \[
    |\WReach_r[G,S,\sigma,u] \cap (S_i \setminus S_{i+1})| \leq (\delta(k,d+1)+3) \cdot r.
    \]
    Since each vertex of $S_{i+1}$ precedes each vertex of $S_i$ in $\sigma$, we have
    \[
    \WReach_r[G,S,\sigma,u] \cap (S_i \setminus S_{i+1}) \subseteq \WReach_r[G-S_{i+1}, S-S_{i+1}, \sigma \vert_{S \setminus S_{i+1}}, u] \cap (S_i \setminus S_{i+1}).
    \]
    Let $\alpha = \max\{a \in I_{i+1} \mid a \leq j_u\}$ and $\beta = \alpha + 2^{i+1}$.
    Let $C$ be the component of $u$ in $G-S_{i+1}$.
    Since $W_{\alpha-1} \cap W_\alpha, W_{\beta-1} \cap W_\beta \subseteq S_{i+1}$,
    and because $(W_0, \dots, W_\ell)$ is a path decomposition of $G$,
    $V(C) \cap S \subseteq \bigcup_{a \in \{\alpha, \dots, \beta-1\}} W_a$.
    We deduce that
    \[
    \WReach_r[G-S_{i+1},S\setminus S_{i+1}, \sigma\vert_{S \setminus S_{i+1}}, u] \cap (S_i \setminus S_{i+1})
    \subseteq \bigcup_{a \in \{\alpha, \dots, \beta-1\}} W_a.
    \]
    Since the only members of $I_i \setminus I_{i+1}$ in $\{\alpha+1, \dots, \beta-1\}$ is $\gamma = \alpha+2^i$,
    we in fact have
    \[
    \WReach_r[G-S_{i+1},S\setminus S_{i+1}, \sigma\vert_{S \setminus S_{i+1}}, u] \cap (S_i \setminus S_{i+1})
    \subseteq R'_\gamma,
    \]
    and we deduce that
    \begin{align*}
        |\WReach_r[G-S_{i+1},S\setminus S_{i+1}, \sigma\vert_{S \setminus S_{i+1}}, u] \cap (S_i \setminus S_{i+1})| &\leq \wcol_r(G - S_{i+1}, R_\gamma', \sigma_\gamma) \\
        &\leq (\delta(k,d+1)+3) \cdot r.
    \end{align*}

    For convenience let $S_{s+1} = \emptyset$.
    Since $S = S_0$, it follows that
    \begin{align*}
        |\WReach_r[G,S,\sigma,u]|
        &\leq \sum_{i \in \{0, \dots, s\}} |\WReach_r[G,S,\sigma,u] \cap (S_i \setminus S_{i+1})| \\
        &\leq (s+2) \cdot (\delta(k,d+1)+3) \cdot r \\
        &\leq 5(\delta(k,d+1)+3)\cdot r \log r. \qedhere
    \end{align*}
\end{proof}

To generalize \Cref{lemma:excluding_a_star_wcol} to graphs with no $\mathcal{F}$-rich model of $F_{h,d}$ for $h > 2$, we need the following straightforward property.

\begin{lemma}\label{lemma:rooting_a_model_of_Fhd}
    Let $h,d$ be positive integers.
    Let $G$ be a connected graph.
    If there is a model $\big(B_x \mid x \in V(F_{h,d+1})\big)$ of $F_{h,d+1}$ in $G$,
    then for every $u \in V(G)$,
    there is a model $\big(B'_x \mid x \in V(F_{h,d})\big)$ of $F_{h,d}$ in $G$ such that
    \begin{enumerateOurAlph}
        \item $u \in B'_s$, where $s$ is the root of $F_{h,d}$, and \label{lemma:rooting_a_model_of_Fhd:i}
        \item for every $x \in V(F_{h,d})$, $B_y \subseteq B'_x$ for some $y \in V(F_{h,d+1})$. \label{lemma:rooting_a_model_of_Fhd:ii}
    \end{enumerateOurAlph}
\end{lemma}

\begin{proof}
    Suppose that there is a model $\big(B_x \mid x \in V(F_{h,d+1})\big)$ of $F_{h,d+1}$ in $G$. 
    Since $G$ is connected, we can assume that $\bigcup_{x \in V(F_{h,d+1})} B_x = V(G)$.
    Let $s_0$ be the root of $F_{h,d+1}$.
    There is a subtree $T'$ of $F_{h,d+1}$ rooted in a child of $s_0$ such that $u \in \bigcup_{x \in V(T') \cup \{s_0\}} B_x$.
    Define $B_s' = \bigcup_{x \in V(T') \cup \{s_0\}} B_x$ and $B_x' = B_x$ for every $x \in V(F_{h,d+1}) \setminus (\{s_0\} \cup V(T'))$.
    The collection $\big(B'_x \mid x \in V(F_{h,d})\big)$ is a model of $F_{h,d}$ in $G$ satisfying \ref{lemma:rooting_a_model_of_Fhd:i} and \ref{lemma:rooting_a_model_of_Fhd:ii}.
\end{proof}

The following lemma is the case of \Cref{thm:main} where $t=2$.

\begin{lemma}\label{lem:new_base_case}
    Let $k,h,d$ be positive integers with $h \geq 2$.
    There is an integer $c_0(h,d,k)$ such that
    for every connected $K_k$-minor-free graph $G$,
    for every family $\mathcal{F}$ of connected subgraphs of $G$,
    if $G$ has no $\mathcal{F}$-rich model of $F_{h,d}$, then
    there is a set $S \subseteq V(G)$
    such that
    \begin{enumerateOurAlph} 
        \item $V(F) \cap S \neq \emptyset$ for every $F \in \mathcal{F}$;\label{item:new_base_case_i}
        \item $G[S]$ is connected; \label{item:new_base_case_ii}
        \item $\wcol_r(G,S) \leq c_0(h,d,k) \cdot r \log r$ for every integer $r$ with $r \geq 2$.\label{item:new_base_case_iii}
    \end{enumerateOurAlph}
\end{lemma}

\begin{proof}
    We proceed by induction on $h$.
    For $h=2$, the result is given by \Cref{lemma:excluding_a_star_wcol} setting $c_0(1,d,k) = 5(\delta(k,d+1)+3)$.
    Next, assume $h>2$ and that $c_0(h-1,d,k)$ witnesses the assertion for $h-1$.
    Let $c_0(h,d,k) = 5(\delta(k,d+1)+3) + 3 + c_0(h-1,d+1,k)$.

    Let $G$ be a connected $K_k$-minor-free graph and let $\mathcal{F}$ be a family of connected subgraphs of $G$.
    Suppose that $G$ has no $\mathcal{F}$-rich model of $F_{h,d}$.
    Let $r$ be an integer with $r \geq 2$.
    Let $\mathcal{F}'$ be the family of all the connected subgraphs $H$ of $G$ such that $H$ contains an $\mathcal{F}\vert_H$-rich model of $F_{h-1,d+1}$.
    We claim that there is no $\mathcal{F}'$-rich model of $F_{2,d}$ in $G$.
    Suppose to the contrary that $\big(B_x \mid x \in V(F_{2,d})\big)$ is such a model.
    Let $s$ be the root of $F_{2,d}$ and let $s'$ be the root of $F_{h-1,d}$.
    For every $x \in V(F_{2,d}) \setminus \{s\}$,
    by \Cref{lemma:rooting_a_model_of_Fhd},
    there is an $\mathcal{F}$-rich model $\big(C_y \mid y \in V(F_{h-1,d})\big)$ of $F_{h-1,d}$ in $G[B_x]$ such that $C_{s'}$ contains a vertex of $N_G(B_s) \cap B_x$.
    The union of these models together with $B_s$ yields an $\mathcal{F}$-rich model of $F_{h,d}$ in $G$, which is a contradiction.
    See \Cref{fig:forest-models}.
    
    \begin{figure}[tp]
        \centering 
        \includegraphics[scale=1]{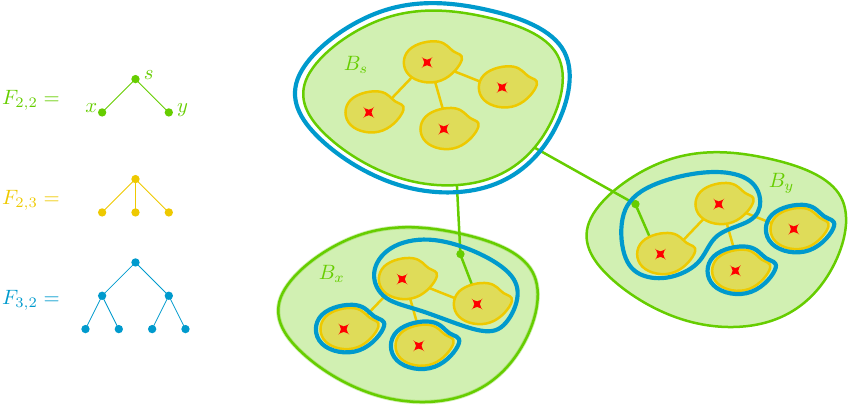} 
        \caption{
            We provide an example of the construction of an $\mathcal{F}$-rich model of $F_{h,d}$ in $G$ assuming that there is an $\mathcal{F}'$-rich model of $F_{2,d}$ in $G$ in the case where $h=3$ and $d=2$.
            In green, we depict an $\mathcal{F}'$-rich model of $F_{2,d} = F_{2,2}$ in the graph.
            Each branch set contains an $\mathcal{F}$-rich model of $F_{h-1,d+1} = F_{2,3}$.
            We depict these models in yellow and the red stars are the elements of $\mathcal{F}$.
            The obtained model of $F_{h,d} = F_{3,2}$ we depict in blue.
            Note that this model is $\mathcal{F}$-rich.
        }
        \label{fig:forest-models}
    \end{figure} 

    Since $G$ has no $\mathcal{F}'$-rich model of $F_{2,d}$, by \Cref{lemma:excluding_a_star_wcol}, there is a set $S_0 \subseteq V(G)$ such that
    \begin{enumerate}[label={\normalfont\ref{lemma:excluding_a_star_wcol}.(\makebox[\mywidth]{\alph*})}]
        \item for every $F \in \mathcal{F}'$, $V(F) \cap S_0 \neq \emptyset$; \label{lemma:star:item:hitting}
        \item $G[S_0]$ is connected; \label{lemma:star:item:connected}
        \item $\wcol_r(G,S_0) \leq 5(\delta(k,d+1)+3) \cdot r \log r$. \label{lemma:star:item:wcol}
    \end{enumerate}
    Let $C$ be a component of $G - S_0$. 
    By~\ref{lemma:star:item:hitting}, $C \notin \mathcal{F}'$, and so, $C$ has no $\mathcal{F}\vert_C$-rich model of $F_{h-1,d+1}$.
    Therefore, by induction hypothesis, there is a set $S_C \subseteq V(C)$ such that
    \begin{enumerateOurAlphPrim}
        \item $V(F) \cap S_C \neq \emptyset$ for every $F \in \mathcal{F}\vert_C$; \label{item:new_base_case_i_induction}
        \item $C[S_C]$ is connected; \label{item:new_base_case_iii_induction}
        \item $\wcol_r(C,S_C) \leq c_0(h-1,d+1,k) \cdot r \log r$.\label{item:new_base_case_ii_induction}
    \end{enumerateOurAlphPrim}
    Let $Q_C$ be an $S_C$-$N_G(S_0)$ geodesic in $G$.
    In particular, $Q_C$ is a geodesic in $C$.
    Let $\mathcal{C}$ be the family of components of $G-S_0$ and let 
    \[
    S = S_0 \cup \bigcup_{C \in \mathcal{C}} (S_C \cup V(Q_C)).
    \]
    See~\Cref{fig:excluding-forest-combine} for an illustration.
    We claim that \ref{item:new_base_case_i}-\ref{item:new_base_case_iii} hold.

    \begin{figure}[tp]
        \centering 
        \includegraphics[scale=1]{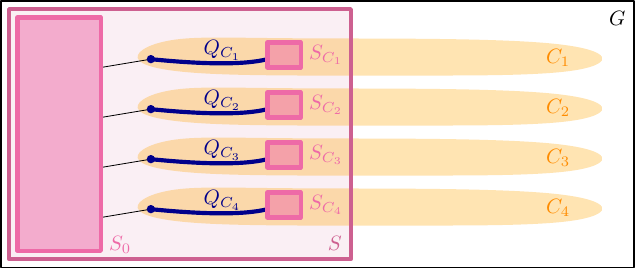} 
        \caption{
            An illustration of the construction of the set $S$ in the proof of \Cref{lem:new_base_case}.
        }
        \label{fig:excluding-forest-combine}
    \end{figure} 

    Let $F \in \mathcal{F}$.
    If $V(F) \cap S_0 = \emptyset$, then $V(F) \subset V(C)$ for some component $C$ of $G - S_0$.
    In particular, $F \in \mathcal{F}\vert_C$, and thus, by~\ref{item:new_base_case_i_induction}, $V(F) \cap S_C \neq \emptyset$, which proves~\ref{item:new_base_case_i}.
    The graph $G[S]$ is connected by construction, \ref{item:new_base_case_iii_induction} and \ref{lemma:star:item:connected}, which yields~\ref{item:new_base_case_ii}.
    The following sequence of inequalities concludes the proof of~\ref{item:new_base_case_iii} and the lemma:
    \begin{align*}
    \wcol_r(G,S) &\leq \wcol_r(G,S_0) + \wcol_r\left(G-S_0,\bigcup_{C \in \mathcal{C}} \big(S_C \cup V(Q_C)\big)\right)&&\textrm{by~\Cref{obs:wcol_union}}\\ 
    &\leq \wcol_r(G,S_0) + \max_{C \in \mathcal{C}} \wcol_r(C, S_C \cup V(Q_C))&&\textrm{by~\Cref{obs:wcol_components2}}\\ 
    &\leq \wcol_r(G,S_0) + \max_{C \in \mathcal{C}} \wcol_r(C, S_C) + (2r+1)&&\textrm{by~\Cref{obs:geodesics}}\\ 
    &\leq 5(\delta(k,d+1)+3) \cdot r\log r + c_0(h-1,d+1,k) \cdot r\log r + 3r&&\textrm{by~\ref{lemma:star:item:wcol} and \ref{item:new_base_case_ii_induction}}\\ 
    &\leq (5(\delta(k,d+1)+3) + c_0(h-1,d+1,k) + 3) \cdot r\log r\\
    &= c_0(h,d,k) \cdot r \log r. &&\qedhere
    \end{align*}    
\end{proof}

\section{Proof of the main theorem}\label{sec:induction}

In this section, we prove \Cref{thm:rooted_main}.
As already mentioned, the proof is by induction, and in \Cref{sec:base} we covered the base case.
The induction statement is encapsulated in \Cref{thm:main}.
Note that in order to obtain \Cref{thm:rooted_main} as a corollary of \Cref{thm:main} one has to apply it to each connected component of $G$ with $k = |V(X)|$ and with the family of all one-vertex subgraphs of $G$ as $\mathcal{F}$.
Note that in such a case $S$ must be equal to $V(G)$.

\begin{theorem}[\Cref{thm:main} restated]
    Let $k$ and $t$ be positive integers with $t \geq 2$.
    Let $X$ be a graph with $\rtd_2(X) \leq t$.
    There exists an integer $c(t,X,k)$ 
    such that
    for every connected $K_k$-minor-free graph $G$, 
    for every family $\mathcal{F}$ of connected subgraphs of $G$,
    if $G$ has no $\mathcal{F}$-rich model of $X$, then
    there exists $S \subseteq V(G)$ such that
    \begin{enumerateOurAlphCapital}
        \item $V(F) \cap S \neq \emptyset$ for every $F \in \mathcal{F}$;\label{thm:main:item:hit}
        \item $G[S]$ is connected; \label{thm:main:item:connected}
        \item $\wcol_r(G,S) \leq c(t,X,k) \cdot r^{t-1} \log r$ for every integer $r$ with $r \geq 2$. \label{thm:main:item:wcol}
    \end{enumerateOurAlphCapital}
\end{theorem}

\begin{proof}
    We proceed by induction on $t$.
    When $t = 2$, by \Cref{obs:char:rtd_2=2}, $X$ is a forest.
    Let $h,d$ be positive integers such that $X \subset F_{h,d}$,
    and let $c_0(h,d,k)$ be the constant given by \Cref{lem:new_base_case}.
    The assertion with $c(2,X,k) = c_0(h,d,k)$ follows by applying \Cref{lem:new_base_case}.
    Next, let $t \geq 3$, and assume that the result holds for $t-1$. 
    We refer to this property as the \emph{main induction hypothesis}. 
    
       \begin{claim1}\label{claim:adding_apices_J}
        Let $Y$ be a graph with $\rtd_2(Y) \leq t-1$.
        There is an integer $c_1(t,Y,k)$ such that
        for every connected $K_k$-minor-free graph $G$,
        for every $u \in V(G)$,
        for every family $\mathcal{F}$ of connected subgraphs of $G$,
        if $G$ has no $\mathcal{F}$-rich model of $K_1 \oplus Y$, then
        there exist $S \subseteq V(G)$, a tree $T$ rooted in $s \in V(T)$, and a tree partition $\big(T,(P_x \mid x \in V(T))\big)$
        of $G[S]$
        with $P_s = \{u\}$
        such that
        \begin{enumerateOurAlph}
            \item $V(F) \cap S \neq \emptyset$ for every $F \in \mathcal{F}$; \label{item:claim_adding_apices:hitting}
            \item $G[S]$ is connected; \label{item:claim_adding_apices:connected}
            \item for every component $C$ of $G-S$, $N_{G}(V(C)) \subset P_x \cup P_y$ for some $x,y \in V(T)$ with either $x=y$ or $xy$ is an edge in $T$;\label{item:claim_adding_apices:components}
            \item for every $x \in V(T)$, \[\wcol_r\left(G_x, P_x\right) \leq c_1(t,Y,k) \cdot r^{t-2} \log r\] for every integer $r$ with $r \geq 2$,
            where, 
            for $T_x$ being the subtree of $T$ rooted in $x$, 
            $G_x$ is the subgraph of $G$ induced by the union of $U_x = \bigcup_{y \in V(T_x)} P_y$ and all the vertex sets of the components of $G-S$ having a neighbor in $U_x$.\label{item:claim_adding_apices:wcol}
        \end{enumerateOurAlph}
    \end{claim1}


    The statement of the claim is visualized in Figure~\ref{fig:excluding-tree-statement}.

    \begin{figure}[tp]
        \centering 
        \includegraphics[scale=1]{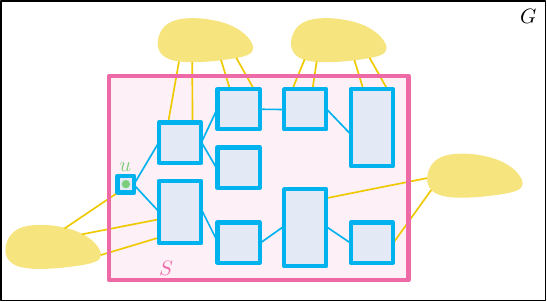} 
        \caption{
        An illustration of the statement of \Cref{claim:adding_apices_J}.
        }
        \label{fig:excluding-tree-statement}
    \end{figure} 

    \begin{proofclaim}
        Let $X' =  K_1 \sqcup Y$ and let $c_1(t,Y,k) = c(t-1, X', k)+3$.
        We proceed by induction on $|V(G)|$.
        If $|V(G)| = 1$, then since $c_1(t,Y,k) \geq 1$, the result holds.
        Next, suppose that $|V(G)| > 1$.

        First, assume that $G-\{u\}$ is not connected.
        Let $\mathcal{C}$ be the family of all the components of $G-\{u\}$.
        Consider a component $C \in \mathcal{C}$.
        By induction hypothesis applied to $G_C = G[V(C) \cup \{u\}]$, $u$, and $\mathcal{F}\vert_{G_C}$,
        there exist $S_C \subseteq V(G_C)$, a tree $T_C$ rooted in $s_C \in V(T_C)$, 
        and a tree partition $(T_C,(P_{C,x} \mid x \in V(T_C)))$ of $S_C$ in $G_C$ with $P_{C,s_C} = \{u\}$ such that
        \begin{enumerateOurAlphPrim}
            \item $V(F) \cap S_C \neq \emptyset$ for every $F \in \mathcal{F}\vert_{G[V(C) \cup \{u\}]}$; \label{item:claim_adding_apices:hitting'_0}
            \item $G[S_C]$ is connected; \label{item:claim_adding_apices:connected'_0}
            \item for every component $C'$ of $G_C-S_C$, $N_{G_C}(V(C')) \subset P_{C,x} \cup P_{C,y}$ for some $x,y \in V(T_C)$ and either $x=y$ or $xy$ is an edge in $T_C$; \label{item:claim_adding_apices:components'_0}
            \item for every $x \in V(T_C)$, \[\wcol_r\left(G_{C,x}, P_{C,x}\right) \leq c_1(t,Y,k) \cdot r^{t-2} \log r\] for every integer $r$ with $r \geq 2$,
            where,
            for $T_{C,x}$ being the subtree of $T_C$ rooted in $x$,
            $G_{C,x}$ is the subgraph of $G_C$ induced by the union of $U_{C,x} = \bigcup_{y \in V(T_{C,x})} P_{C,y}$ and all the vertex sets of the components of $G_C-S_C$ having a neighbor in $U_{C,x}$\label{item:claim_adding_apices:wcol'_0}
        \end{enumerateOurAlphPrim}
        Then let $S = \bigcup_{C \in \mathcal{C}} S_C$, let $T$ be the tree obtained from the disjoint union of all the $T_C$ for $C \in \mathcal{C}$
        by identifying all the vertices in $\{s_C \mid C \in \mathcal{C}\}$ into a single vertex $s$.
        Finally, let $P_s = \{u\}$ and $P_x = P_{C,x}$ for every $C \in \mathcal{C}$ and $x \in V(T_C) \setminus \{s_C\}$.
        Then we claim that \ref{item:claim_adding_apices:hitting}-\ref{item:claim_adding_apices:wcol} hold. 
        Indeed, for every $F\in \mathcal{F}$, $F$ is connected so either $F$ is a subgraph of $G-\set{u}$ or $F$ contains $u$. In both cases, we see that $V(F)\cap S\neq\emptyset$, so~\ref{item:claim_adding_apices:hitting} holds. 
        Item~\ref{item:claim_adding_apices:connected} holds as $u\in S_C$ and $G[S_C]$ is connected for all $C\in\mathcal{C}$. 
        Item~\ref{item:claim_adding_apices:components} holds as every component of $G-S$ is a component of $G_C-S_C$ for some $C\in\mathcal{C}$.
        Finally, for all $x\in V(T)$ with $x\neq s$ item~\ref{item:claim_adding_apices:wcol} follows directly from~\ref{item:claim_adding_apices:wcol'_0} and the construction of $T$, while when $x=s$, $\wcol_r\left(G_s, P_s\right) = 1$ for every integer $r$ with $r \geq 2$. 
        From now on we assume that $G-\{u\}$ is connected.

        Let $\mathcal{F}'$ be the family of all the connected subgraphs $H$ of $G-\{u\}$ such that $u \in N_G(V(H))$ and
        $F \subseteq H$ for some $F \in \mathcal{F}$.
        We argue that there is no $\mathcal{F}'$-rich model of $X'$ in $G-\{u\}$.
        Suppose by contradiction that it exists.
        Such a model would contain an $\mathcal{F}$-rich model of $K_1 \sqcup Y$ in $G-\{u\}$ 
        such that every branch set is adjacent to $u$.
        By adding $u$ to the branch set corresponding to $K_1$ in $K_1 \sqcup Y$,
        we obtain an $\mathcal{F}$-rich model of $K_1 \oplus Y$ in $G$, which is a contradiction.
        This proves that there is no $\mathcal{F}'$-rich model of $X'$ in $G-\{u\}$.

        Since $\rtd_2(X') \leq \max \{\rtd_2(Y),1\} \leq t-1$, by the main induction hypothesis applied to $X', G - \{u\}$, and $\mathcal{F}'$,
        there exists a set $S_0 \subseteq V(G-\{u\})$ such that
        \begin{enumerateOurAlphCapitalPrim}
            \item $V(F) \cap S_0 \neq \emptyset$ for every $F \in \mathcal{F}'$; \label{thm:main:item:hit'_J}
            \item $(G-\{u\})[S_0]$ is connected; \label{thm:main:item:connected'_J}
            \item $\wcol_r(G-\{u\},S_0) \leq c(t-1, X',k) \cdot r^{t-2} \log r$ for every integer $r$ with $r \geq 2$. \label{thm:main:item:wcol'_J}
        \end{enumerateOurAlphCapitalPrim}
        By possibly adding an arbitrary vertex of $V(G-\{u\})$ to $S_0$, we can assume $S_0 \neq \emptyset$.
        Let $Q$ be a $u$-$S_0$ geodesic in $G$ and let $S_1 = \left(S_0 \cup V(Q)\right) - \{u\}$.
        Note that $G[\{u\} \cup S_1]$ is connected by~\ref{thm:main:item:connected'_J}.

        Let $\mathcal{C}_1$ be the family of all the components $C$ of $G-(\{u\}\cup S_1)$ such that $N_G(u)\cap V(C)=\emptyset$. 
        Consider $C \in \mathcal{C}_1$.
        Since $G$ is connected and $u \not\in N_G(V(C))$, there is an edge between $V(C)$ and $S_1$ in $G$.
        Let $G_C$ be obtained from $G[V(C) \cup S_1]$ by contracting $S_1$ 
        into a single vertex $u_C$.
        Note that $|V(G_C)| < |V(G)|$ since $u \not\in V(G_C)$.
        Since $G_C$ is a minor of $G$, $G_C$ has no $\mathcal{F}\vert_C$-rich model of $K_1 \oplus Y$.
        By induction hypothesis applied to $G_C$, $u_C$, and $\mathcal{F}\vert_C$, there exist $S_C \subseteq V(G_C)$, a tree $T_C$ rooted in $s_C \in V(T_C)$, 
        and a tree partition $\big(T_C,(P_{C,x} \mid x \in V(T_C))\big)$ of $S_C$ in $G_C$ with $P_{C,s_C} = \{u_C\}$
        such that
        \begin{enumerateOurAlphPrimPrim}
            \item $V(F) \cap S_C \neq \emptyset$ for every $F \in \mathcal{F}\vert_C$; \label{item:claim_adding_apices:hitting'}
            \item $G_C[S_C]$ is connected; \label{item:claim_adding_apices:connected'}
            \item for every component $C'$ of $G_C-S_C$, $N_{G_C}(V(C')) \subset P_{C,x} \cup P_{C,y}$ for some $x,y \in V(T_C)$ and either $x=y$ or $xy$ is an edge in $T_C$; \label{item:claim_adding_apices:components'}
            \item for every $x \in V(T_C)$, \[\wcol_r\left(G_{C,x}, P_{C,x}\right) \leq c_1(t,Y,k) \cdot r^{t-2} \log r\] for every integer $r$ with $r \geq 2$,
            where,
            for $T_{C,x}$ being the subtree of $T_C$ rooted in $x$, 
            $G_x$ is the subgraph of $G$ induced by the union of $U_x = \bigcup_{y \in V(T_x)} P_y$ and all the vertex sets of the components of $G_C-S_C$ having a neighbor in $U_x$.\label{item:claim_adding_apices:wcol'}
        \end{enumerateOurAlphPrimPrim}
        Let
        \[
            S = \{u\} \cup S_1 \cup \bigcup_{C \in \mathcal{C}_1} (S_C \setminus \{u_C\}).
        \]
        Let $T$ be obtained from the disjoint union of $\{T_C \mid C\in\mathcal{C}_1\}$ by identifying the vertices $\set{s_{C}\mid C\in\mathcal{C}_1}$ into a new vertex $s'$ and by adding a new vertex $s$ adjacent to $s'$ in $T$. 
        Let $P_s=\set{u}$, $P_{s'}=S_1$, and
        for each $C\in\mathcal{C}_1$, 
        $x \in V(T_{C} \setminus \{s_{C}\})$, 
        let $P_x = P_{C,x}$.
        See \Cref{fig:excluding-tree} for the illustration of this construction.

        In order to conclude, we argue that $\big(T,(P_x \mid x \in V(T))\big)$ is a tree partition of $G[S]$ 
        and \ref{item:claim_adding_apices:hitting}-\ref{item:claim_adding_apices:wcol} hold.

    \begin{figure}[tp]
        \centering 
        \includegraphics[scale=1]{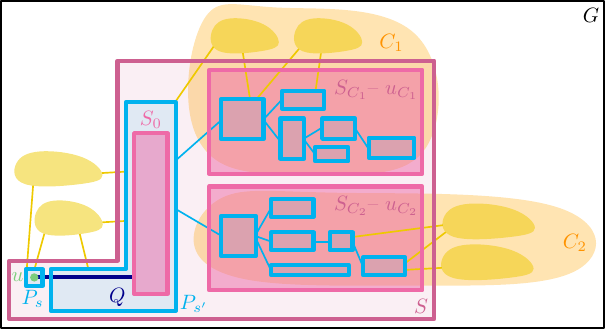} 
        \caption{
        An illustration of the construction of $S$ and its tree partition in the proof of \Cref{claim:adding_apices_J}.
        Note that in the sketched case $\mathcal{C}_1 = \{C_1,C_2\}$.
        }
        \label{fig:excluding-tree}
    \end{figure} 

        Since for every $C \in \mathcal{C}_1$, $u \notin N_G(V(C))$, every edge in $G[S]$ containing $u$ has another endpoint in $S_1 = P_{s'}$.
        Consider an edge $vw$ in $G[S]$ such that $v \in S_1$ and $w \in S_C$ for some $C \in \mathcal{C}_1$.
        Since $\big(T_C,(P_{C,x} \mid x \in V(T_C))\big)$ is a tree partition of $G_C[S_C]$ with $P_{C,s_C} = \{u_C\}$ and $u_C$ is the result of the contraction of $S_1$,
        we conclude that $w \in P_x$ for some $x \in V(T_C)$ such that $s'x$ is an edge in $T$.
        Finally, for every edge $vw$ of $G[S]$ with $v,w \not\in \{u\} \cup S_1$, $vw$ is an edge of $G[S_C \setminus \{u_C\}]$ for some component $C \in \mathcal{C}_1$,
        and so $v \in P_{C,x}$ and $w \in P_{C,y}$ for adjacent or identical vertices $v,w$ of $T_C$.
        Then $v \in P_x$ and $w \in P_y$.
        It follows that $\big(T,(P_{x} \mid x \in V(T))\big)$ is a tree partition of $G[S]$.

        Let $F \in \mathcal{F}$.
        If $V(F) \cap (\{u\} \cup S_1) \neq \emptyset$, then $V(F) \cap S \neq \emptyset$.
        Otherwise, $F \subseteq G - \{u\}$ and $V(F) \cap S_0 = \emptyset$, and therefore by \ref{thm:main:item:hit'_J}, $F \notin \mathcal{F}'$, so in particular, $u \notin N_G(V(F))$.
        In this case, there is a component $C \in \mathcal{C}_1$ such that $F \in \mathcal{F}\vert_C$, thus, $V(F) \cap S_C \neq \emptyset$ by \ref{item:claim_adding_apices:hitting'}.
        This proves~\ref{item:claim_adding_apices:hitting}.


        Item~\ref{item:claim_adding_apices:connected} holds since $G[\{u\} \cup S_1]$ is connected and for every $C \in \mathcal{C}_1$, 
        $G_C[S_C]$ is connected by~\ref{item:claim_adding_apices:connected'} and $u_C \in S_C$.

        For every component $C'$ of $G-S$, either $N_G(V(C')) \subseteq \{u\} \cup S_1 = P_s \cup P_{s'}$, or $C' \subseteq C$ for some $C \in \mathcal{C}_1$.
        In the latter case, $C'$ is a component of $G_C - S_C$, and $u \not\in N_G(V(C))$.
        By~\ref{item:claim_adding_apices:components'}, there is $x,y \in V(T_C)$ such that $N_{G_C}(V(C')) \subseteq P_{C,x} \cup P_{C,y}$, and thus, $N_G(V(C')) \subseteq P_x \cup P_y$.
        This proves \ref{item:claim_adding_apices:components}.

        Finally, we argue~\ref{item:claim_adding_apices:wcol}.
        For every $x \in V(T)$, we denote by $T_x$ the subtree of $T$ rooted in $x$, 
        and by $U_x$ the subgraph of $G$ induced by $U_x = \bigcup_{y \in V(T_x)} P_y$ with the vertex sets of all the components of $G-S$ having a neighbor in $U_x$.
        Let $r$ be an integer with $r \geq 2$ and let $x \in V(T)$.
        For $x = s$, $|P_{s}| = 1$, thus the assertion is clear.
        For $x = s'$, we have $G_{s'}$ is a union of components of $G-\{u\}$.
        By~\ref{thm:main:item:wcol'_J}, $\wcol_r(G_{s'}, S_0) = \wcol_r(G-\{u\}, S_0) \leq c(t-1,X',k) \cdot r^{t-2} \log r$.
        Since $Q$ is a geodesic in $G$, by \Cref{obs:geodesics_plus_monotone}, 
        \[
            \wcol_r(G-\{u\}, S_1) \leq c(t-1,X',k) \cdot r^{t-2} \log r + (2r+1) \leq c_1(t,Y,k) \cdot r^{t-2}\log r.
        \]
        For $x \in V(T_C - \{s_C\})$ for some $C \in \mathcal{C}_1$, we have $T_x = T_{C,x}$, thus, the asserted inequality follows from~\ref{item:claim_adding_apices:wcol'}.
        This ends the proof of the claim.
    \end{proofclaim}

    \Cref{claim:adding_apices_J} yields the following less technical statement.
    \begin{claim1}\label{claim:adding_apices_less_technical_J}
        Let $Y$ be a graph with $\rtd_2(Y) \leq t-1$.
        There is an integer $c_2(t,Y,k)$ such that for every connected $K_k$-minor-free graph $G$,
        for every family $\mathcal{F}$ of connected subgraphs of $G$,
        if $G$ has no $\mathcal{F}$-rich model of $K_1 \oplus Y$, then
        there exists $S \subseteq V(G)$ such that
        \begin{enumerateOurAlph}
            \item $V(F) \cap S \neq \emptyset$ for every $F \in \mathcal{F}$;\label{claim:adding_apices_J:hitting}
            \item $G[S]$ is connected;\label{claim:adding_apices_J:connected}
            \item $\wcol_r(G,S) \leq c_2(t,Y,k) \cdot r^{t-1} \log r$ for every integer $r$ with $r \geq 2$.\label{claim:adding_apices_J:wcol}
        \end{enumerateOurAlph}
    \end{claim1}

    \begin{proofclaim}
        Let $c_2(t,Y,k) = 2 \cdot c_1(t,Y,k)$ and let $r$ be an integer with $r \geq 2$.
        We apply \Cref{claim:adding_apices_J} with an arbitrary vertex $u \in V(G)$ to obtain $S \subseteq V(G)$, a tree $T$ rooted in $s \in V(T)$, and a tree partition $\big(T,(P_x \mid x \in V(T))\big)$
        of $S$ in $G$ with $P_s = \{u\}$
        such that
        \begin{enumerate}[label={\normalfont\ref{claim:adding_apices_J}.(\makebox[\mywidth]{\alph*})}]
            \item $V(F) \cap S \neq \emptyset$ for every $F \in \mathcal{F}$; \label{item:claim_adding_apices:hitting-called}
            \item $G[S]$ is connected; \label{item:claim_adding_apices:connected-called}
            \item for every component $C$ of $G-S$, $N_{G}(V(C)) \subset P_x \cup P_y$ for some $x,y \in V(T)$ and either $x=y$ or $xy$ is an edge in $T$;\label{item:claim_adding_apices:components-called}
            \item for every $x \in V(T)$, \[\wcol_r\left(G_x, P_x\right) \leq c_1(t,Y,k) \cdot r^{t-2} \log r\] for every integer $r$ with $r \geq 2$,
            where, for $T_x$ being the subtree of $T$ rooted in $x$, 
            the graph $G_x$ is the union of $U_x = \bigcup_{y \in V(T_x)} P_y$ with the vertex sets of all the components of $G-S$ having a neighbor in $U_x$. \label{item:claim_adding_apices:wcol-called}
        \end{enumerate}
        Items \ref{claim:adding_apices_J:hitting} and \ref{claim:adding_apices_J:connected} hold by \ref{item:claim_adding_apices:hitting-called} and \ref{item:claim_adding_apices:connected-called} respectively.
        It suffices to prove \ref{claim:adding_apices_J:wcol}.

        For each $x \in V(T)$, let $\sigma_x$ be an ordering of $P_x$ witnessing \ref{item:claim_adding_apices:wcol-called} and let $\sigma' = (x_1, \dots, x_{|V(T)|})$ be an elimination ordering of $T$.
        Finally, let $\sigma$ be the concatenation of $\sigma_{x_1},\dots,\sigma_{x_{|V(T)|}}$ in this order.

        Let $u \in V(G)$.
        To conclude the claim, we argue that
            \[|\WReach_r[G,S,\sigma,u]| \leq c_2(t,Y,k) \cdot r^{t-1} \log r.\]
        Let $x_u \in V(T)$ be such that if $u \in S$, then $u \in P_{x_u}$, and otherwise, $x_u$ is the vertex of $T$ furthest to $s$ such that $P_{x_u}$ intersects $N_G(V(C))$, where $C$ is the component of $u$ in $G-S$.
        Let $A$ be the set of all the ancestors of $x_u$ in $T$ in distance at most $r$ from $x_u$.
        In particular, $|A| \leq r+1$.
        By \Cref{obs:wcols:trees} and \ref{item:claim_adding_apices:components-called}, 
            \[\WReach_r[G,S,\sigma,u] \subset \bigcup_{y \in A} P_y. \]
        Let $y \in A$.
        Since $\sigma$ extends $\sigma_y$ and $x_u \in V(T_y)$, by~\ref{item:claim_adding_apices:wcol-called}, we have
            \[|\WReach_r[G,S,\sigma,u] \cap P_y| \leq \wcol_r\left(G_y, P_y\right) \leq c_1(t,Y,k) \cdot r^{t-2} \log r.\]
        Summarizing,
        \begin{align*}
            |\WReach_r[G,S,\sigma,u]| &= \sum_{y \in A} |\WReach_r[G,S,\sigma,u] \cap P_y| \\
            &\leq (r+1) \cdot c_1(t,Y,k) \cdot r^{t-2} \log r\\
            &\leq 2r \cdot c_1(t,Y,k) \cdot r^{t-2} \log r= c_2(t,Y,k) \cdot r^{t-1} \log r.
        \end{align*}
    This concludes the proof of~\ref{claim:adding_apices_J:wcol}, and the claim follows.
    \end{proofclaim}

    We now move on to the final step of the proof.
    \begin{claim1}\label{claim:last_claim}
        Let $Y$ be a graph with $\rtd_2(Y) \leq t-1$ and let $h,d$ be positive integers.
        There exists an integer $c_3(Y,h,d,k)$ such that
        for every connected $K_k$-minor-free graph $G$ and for every family $\mathcal{F}$ of connected subgraphs of $G$,
        if $G$ has no $\mathcal{F}$-rich model of $T'_{h,d}(Y)$\footnote{See the definition on page \pageref{def:T'hd}.}, 
        then there exists $S \subseteq V(G)$ such that 
        \begin{enumerateOurAlph}
            \item $V(F) \cap S \neq \emptyset$  for every $F \in \mathcal{F}$; \label{item:last_claim_i}
            \item $G[S]$ is connected; \label{item:last_claim_ii}
            \item $\wcol_r(G,S) \leq c_3(Y,h,d,k) \cdot r^{t-1} \log r$ for every integer $r$ with $r \geq 2$. \label{item:last_claim_iii}     
        \end{enumerateOurAlph}
    \end{claim1}

    \begin{proofclaim}
        We proceed by induction on $h$.
        When $h=1$, $T'_{h,d}(Y) = K_1 \oplus (Y \sqcup Y)$ and the result follows from the previous claim applied to $Y \sqcup Y$ (note that $\rtd_2(Y \sqcup Y) = \rtd_2(Y)\leq t-1$)
        by setting $c_3(Y,1,d,k)=c_2(t,Y \sqcup Y,k)$.

        Now assume that $h > 1$ and that the result holds for $h-1$.
        Fix a copy of $Y \sqcup Y$.
        For each $y \in V(Y \sqcup Y)$ add $d$ vertices with $y$ as a unique neighbor.
        Furthermore, add $2d$ isolated vertices.
        We call the obtained graph $Z$ -- see~\Cref{fig:combining-models-T}.
        To keep things in order, we write $V(Z) = V_Y \cup V_Z$, where $V_Y$ are the vertices of $Y \sqcup Y$ in $Z$ and $V_Z$ are all the added vertices.
        Let 
        \[
            c_3(Y,h,d,k) = c_2(t,Z,k) + c_3(Z,h-1,d,k) + 3.
        \]
        By~\ref{rtd2:item:degree1} and \ref{rtd2:item:components}, $\rtd_2(Z) \leq \max\{2,\rtd_2(Y)\} \leq t-1$ since $t \geq 3$.

        \begin{figure}[tp]
            \centering 
            \includegraphics[scale=1]{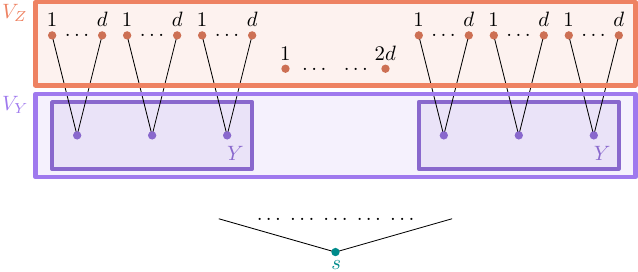} 
            \caption{
            An illustration of the graph $K_1 \oplus Z$.
            }
            \label{fig:combining-models-T}
        \end{figure} 
        
        Let $\mathcal{F}'$ be the family of all the connected subgraphs $H$ of $G$
        such that $H$ has an $\mathcal{F}\vert_H$-rich model of $T'_{h-1,d}(Y)$.
        We claim that there is no $\mathcal{F}'$-rich model of $K_1 \oplus Z$ in $G$.
        Suppose to the contrary that such a model $\big(A_y \mid y \in V(K_1 \oplus Z)\big)$ exists.
        Let $s$ be the vertex of $K_1$ in $K_1 \oplus Z$.
        In particular, $V(K_1 \oplus Z) = \{s\} \cup V_Y \cup V_Z$.
        For every vertex $z \in V_Z$, we define its parent in $K_1 \oplus Z$ in the following way.
        Note that $z$ has at most one neighbor in $Z$.
        If $z$ has a neighbor in $Z$, then the neighbor is its parent and if $z$ is isolated in $Z$, then $s$ is its parent.
        Let $z \in V_Z$ with the parent $p_z$ and let $u_z \in A_z$ be such that there is an edge between $u_z$ and a vertex in $A_{p_z}$. 
        Since the model is $\mathcal{F}'$-rich, $G[A_z]$ contains an $\mathcal{F}\vert_{G[A_z]}$-rich model of $T'_{h-1,d}(Y)$.
        Let $H_z$ be a copy of $T_{h-1,d}(Y)$ with a root $s_z$.
        By \Cref{lemma:root_a_model_of_T'hd} applied to $Y$, $G[A_{z}]$, and $u$,
        there is an $\mathcal{F}\vert_{G[A_z]}$-rich model $\big(B_{z,x} \mid x \in V(H_{z})\big)$ of $H_{z}$ in $G[A_z]$ such that $u_z \in B_{z,s_z}$.
        In particular, there is an edge between $B_{z,s_z}$ and $A_{p_z}$ in $G$.
        Finally, we construct an $\mathcal{F}$-rich model of $T'_{h,d}(Y)$ in $G$.
        Observe that the graph obtained from $(K_1 \oplus Z)[\{s\} \cup V_Y]$ (this graph is isomorphic to $K_1 \oplus (Y \sqcup Y))$ and the disjoint union of $H_z$ for each $z \in V_Z$ by identifying $s_z \in V(H_z)$ with $p_z \in \{s\} \cup V_Y$ for each $z \in V_Z$ is isomorphic to $T'_{h,d}(Y)$.
        For each $p \in \{s\} \cup V_Y$, let $P_y$ be the set of all $z \in V_Z$ such that $p$ is the parent of $z$.
        Let
        \begin{enumerate}[label={\normalfont (\roman*)}]
            \item $D_p = A_p \cup \bigcup_{z \in P_p} B_{z,s_z}$ for every $p \in \{s\} \cup V_Y$ and
            \item $D_x = B_{z,x}$ for every $z \in V_Z$ and $x \in V(H_z \setminus \{s_{z}\})$.
        \end{enumerate}
        It follows that $\big(D_x \mid x \in \{s\} \cup V_Y \cup \bigcup_{z \in V_Z} V(H_z \setminus \{s_z\}) \big)$ is an $\mathcal{F}$-rich model of $T'_{h,d}(Y)$ in $G$.
        This is a contradiction, hence, $G$ has no $\mathcal{F}'$-rich model of $K_1 \oplus Z$.

        \begin{figure}[tp]
            \centering 
            \includegraphics[scale=1]{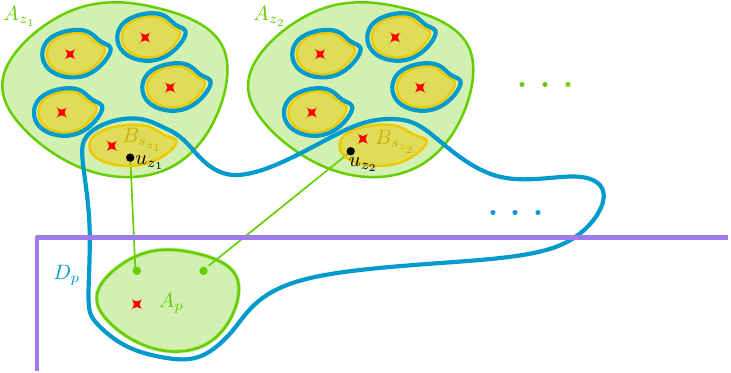} 
            \caption{
            We illustrate how the $\mathcal{F}'$-rich model of $T'_{h,d}(Y)$ in $G$ is constructed.
            We depict only one vertex $p \in \{s\} \cup V_Y$ and two vertices $z_1,z_2 \in V_Z$ such that $p$ is the parent of both $z_1$ and $z_2$.
            The yellow model is a model of $T'_{h-1,d}(Y)$ after applying \Cref{lemma:root_a_model_of_T'hd}.
            }
            \label{T-models}
        \end{figure} 

        By \Cref{claim:adding_apices_less_technical_J}, applied to $Z$, $G$, and $\mathcal{F}'$ there exists $S_0 \subseteq V(G)$ such that
        \begin{enumerate}[label={\normalfont\ref{claim:adding_apices_less_technical_J}.(\makebox[\mywidth]{\alph*})}]
            \item $V(F) \cap S_0 \neq \emptyset$ for every $F \in \mathcal{F}'$;
            \item $G[S_0]$ is connected; \label{claim:adding_apices_less_technical_J-final:connected}
            \item $\wcol_r(G,S_0) \leq c_2(t,Z,k) \cdot r^{t-1} \log r$ for every integer $r$ with $r \geq 2$. \label{claim:adding_apices_less_technical_J-final:wcol}
        \end{enumerate}
        Let $C$ be a component of $G - S_0$.
        Since $V(F) \cap S_0 \neq \emptyset$ for every $F \in \mathcal{F}'$, $C$ has no $\mathcal{F}\vert_C$-rich model of $T'_{h-1,d}(Y)$.
        Therefore, by induction hypothesis, 
        there exists $S_C \subset V(C)$ such that
        \begin{enumerateOurAlphPrim}
            \item $V(F) \cap S_C \neq \emptyset$ for every $F \in \mathcal{F}\vert_C$; \label{final_claim:item:hitting'}
            \item $C[S_C]$ is connected; \label{final_claim:item:connected'}
            \item $\wcol_r(C,S_C) \leq c_3(Y,h-1,d,k) \cdot r^{t-1} \log r$ for every integer $r$ with $r \geq 2$. \label{final_claim:item:wcol'}
        \end{enumerateOurAlphPrim}
        Let $Q_C$ be an $S_C$-$N_G(S_0)$ geodesic in $G$.
        In particular, $Q_C$ is a geodesic in $C$.
        Let $\mathcal{C}$ be the family of all the components of $G-S_0$ and let 
        \[
        S = S_0 \cup \bigcup_{C \in \mathcal{C}} (S_C \cup V(Q_C)).
        \]
        We claim that \ref{item:last_claim_i}-\ref{item:last_claim_iii} hold.
        Let $F \in \mathcal{F}$.
        If $V(F) \cap S_0 = \emptyset$, then $V(F) \subset V(C)$ for some component $C$ of $G - S_0$.
        In particular, $F \in \mathcal{F}\vert_C$, and thus, by~\ref{final_claim:item:hitting'}, $V(F) \cap S_C \neq \emptyset$, which proves~\ref{item:last_claim_i}.
        The graph $G[S]$ is connected by construction, \ref{claim:adding_apices_less_technical_J-final:connected} and \ref{final_claim:item:connected'}, which yields~\ref{item:last_claim_ii}.
        The following sequence of inequalities concludes the proof of~\ref{item:last_claim_iii} and the claim:
        \begin{align*}
            \wcol_r(G,S) &\leq \wcol_r(G,S_0) + \wcol_r\left(G-S_0,\bigcup_{C \in \mathcal{C}} \big(S_C \cup V(Q_C)\big)\right)&&\textrm{by~\Cref{obs:wcol_union}}\\ 
            &\leq \wcol_r(G,S_0) + \max_{C \in \mathcal{C}} \wcol_r(C, S_C \cup V(Q_C))&&\textrm{by~\Cref{obs:wcol_components2}}\\ 
            &\leq \wcol_r(G,S_0) + \max_{C \in \mathcal{C}} \wcol_r(C, S_C) + (2r+1)&&\textrm{by~\Cref{obs:geodesics}}\\ 
            &\leq c_2(t,Z,k) \cdot r^{t-1} \log r + c_3(Y,h-1,d,k) \cdot r^{t-1} \log r + 3r&&\textrm{by~\ref{claim:adding_apices_less_technical_J-final:wcol} and \ref{final_claim:item:wcol'}}\\ 
            &\leq \big(c_2(t,Z,k) + c_3(Y,h-1,d,k) + 3\big) \cdot r^{t-1} \log r\\
            &= c_3(Y,h,d,k) \cdot r^{t-1} \log r.  &&\qedhere
        \end{align*} 
    \end{proofclaim}
    Finally, by \Cref{lemma:Thd_universal_for_rtd}, for every for every graph $X$ with $\rtd_2(X)\leq t$, there exists a graph $Y$ with $\rtd_2(Y) \leq t-1$
    and positive integers $h,d$
    such that $X \subseteq T_{h,d}(Y) \subseteq T'_{h,d}(Y)$.
    By~\Cref{claim:last_claim}, the theorem follows with $c(t,X,k) = c_3(Y,h,d,k)$.
\end{proof}

\section{A tighter bound for graphs of bounded treewidth}\label{sec:bounded_tw}

In this section, we prove \Cref{thm:rooted_main_bounded_tw}.

\subsection{Preliminaries}\label{sec:bounded_tw:prelim}

We start by recalling the notion of tree decomposition and one of its basic properties.
Then, we introduce a refined version of tree decomposition -- natural tree decomposition.
Let $G$ be a graph.
A \emph{tree decomposition} of a graph $G$ is a pair $\mathcal{D} = \big(T,(W_x \mid x \in V(T))\big)$
where $T$ is a tree and $W_x \subseteq V(G)$ for every $x \in V(T)$ satisfying the following conditions:
\begin{enumerate}
    \item for every $u \in V(G)$, $T[\{x \in V(T) \mid u \in W_x\}]$ is a connected subtree of $T$, and
    \item for every edge $uv \in E(G)$, there exists $x \in V(T)$ such that $u,v \in W_x$.
\end{enumerate}
The sets $W_x$ are called \emph{bags} of $\mathcal{D}$.
The \emph{width} of $\mathcal{D}$ is $\max_{x \in V(T)} |W_x|-1$,
and the \emph{treewidth} of $G$, denoted by $\tw(G)$, is the minimum width of a tree decomposition of $G$.\footnote{A \emph{path decomposition} of $G$ is a tree decomposition $\big(T,(W_x \mid x \in V(T))\big)$
of $G$ where $T$ is a path.
The \emph{pathwidth} of $G$, denoted by $\pw(G)$, is the minimum width of a path decomposition of $G$.}

\begin{lemma}[{\cite[Statement (8.7)]{GM5}}]\label{lemma:helly_property_tree_decomposition}
    For every graph $G$, for every tree decomposition $\mathcal{D}$ of $G$, for every family $\mathcal{F}$ of connected subgraphs of $G$, for every positive integer $d$, either
    \begin{enumerate}[label=\normalfont(\arabic*)]
        \item there are $d$ pairwise vertex-disjoint subgraphs in $\mathcal{F}$ or
        \item there is a set $S$ that is the union of at most $d-1$ bags of $\mathcal{D}$ such that $V(F) \cap S \neq \emptyset$ for every $F \in \mathcal{F}$.
    \end{enumerate}
\end{lemma}

A tree decomposition
$\big(T,(W_x\mid x\in V(T))\big)$ of a graph $G$ is \emph{natural} if for every edge $e$ in $T$, 
for each component $T_0$ of $T-e$, the graph $G\left[\bigcup_{z\in V(T_0)} W_z\right]\ \textrm{is connected.}$ 
The following statement appeared first in~\cite{FN06}, see also~\cite{GJNW23}.
\begin{lemma}[{\cite[Theorem~1]{FN06}}]\label{lemma:natural_tree_decomposition}
    Let $G$ be a connected graph and let $\big(T,(W_x\mid x\in V(T))\big)$ be a tree decomposition of $G$. 
    There exists a natural tree decomposition $\big(T',(W'_x\mid x\in V(T'))\big)$ of $G$ such that 
    for every $x'\in V(T')$ there is $x\in V(T)$ with 
    $W'_{x'}\subseteq W_x$.
\end{lemma}

The following lemma is folklore. See e.g.~\cite[Lemma~8]{DHHJLMMRW24} for a proof.

\begin{lemma}\label{lemma:increase_X_to_have_small_interfaces}
Let $m$ be a positive integer.
Let $G$ be a graph and let $\mathcal{D}$ be a tree decomposition of $G$.
If $Y$ is the union of $m$ bags of $\mathcal{D}$, then there is a set $X$ that is the union of at most $2m-1$ bags of $\mathcal{D}$ such that $Y\subseteq X$ and for every component $C$ of $G-X$, $N_G(V(C)) \cap X$ is a subset of the union of at most two bags of $\mathcal{D}$.
Moreover, if $\mathcal{D}$ is natural, then $N_G(V(C))$ intersects at most two components of $G-V(C)$.
\end{lemma}

We need the following technical statement.
For a given set $S$, we say that a collection $S_1,\dots,S_k \subset S$ is a \emph{covering} of $S$ if $S_1 \cup \dots \cup S_k = S$.

\begin{lemma}\label{lemma:rooting_X_in_X'}
    Let $k$ be a positive integer and let $X$ be a graph.
    There exists a graph $X'$ such that $\rtd_2(X')  \leq \rtd_2(X)$ and
    for each covering $S_1, \dots,S_k$ of $V(X')$,
    there exists $i \in [k]$ such that $X'$ contains an 
    $\{H \subseteq X' \mid \text{$H$ connected and $S_i \cap V(H) \neq \emptyset$}\}$-rich model of $X$.
\end{lemma}


\Cref{lemma:rooting_X_in_X'} follows directly from \Cref{lemma:rooting_Thd(X)_in_Thd'(X')} and \Cref{lemma:Thd_universal_for_rtd}.

\begin{lemma}\label{lemma:rooting_Thd(X)_in_Thd'(X')}
    Let $h,d,k$ be positive integers and let $X$ be a graph.
    There exists a positive integer $d'$ and a graph $X'$ with $\rtd_2(X') \leq \rtd_2(X)$  
    such that for each covering $S_1, \dots,S_k$ of $V(T_{h,d'}(X'))$,
    there exists $i \in [k]$ such that $T_{h,d'}(X')$ contains an 
    $\{H \subseteq T_{h,d'}(X') \mid \text{$H$ connected and $S_i \cap V(H) \neq \emptyset$}\}$-rich model of $T_{h,d}(X)$
    whose branch set corresponding to the root of $T_{h,d}(X)$ contains the root of $T_{h,d'}(X')$.
\end{lemma}

\begin{proof}
    We proceed by induction on $(\rtd_2(X),h)$ in the lexicographic order.
    If $\rtd_2(X)=0$, then $X$ is the null graph, $T_{h,d}(X)=K_1$, and the result holds.

    Assume that $\rtd_2(X) > 0 $ and that the result holds for every graph with rooted $2$-treedepth less than $\rtd_2(X)$.
    By \Cref{lemma:Thd_universal_for_rtd}, there is a graph $Z$ such that $\rtd_2(Z) \leq \rtd_2(X) - 1$ and $X \subset T_{h_1,d_1}(Z)$ for some positive integers $h_1$ and $d_1$.

    By induction hypothesis applied to $h_1$, $d_1$, and $Z$, there exists a positive integer $d_1'$ and a graph $Y$ with $\rtd_2(Y) \leq \rtd_2(Z)$ such that for each covering $S_1, \dots,S_k$ of $V(T_{h_1,d_1'}(Y))$,
    there exists $i \in [k]$ such that $T_{h_1,d_1'}(Y)$ contains an 
    $\{H \subseteq T_{h_1,d_1'}(Y) \mid \text{$H$ connected and $S_i \cap V(H) \neq \emptyset$}\}$-rich model of $T_{h_1,d_1}(Z)$ (in particular of $X$)
    whose branch set corresponding to the root of $T_{h_1,d_1}(Z)$ contains the root of $T_{h_1,d_1'}(Y)$.

    In the case of $h=1$. 
    Let $X' = (k+1) \cdot T_{h_1,d_1'}(Y)$ and $d' = 1$.
    In particular, $T_{1,d}(X') = K_1 \oplus ((k+1) \cdot T_{h_1,d_1'}(Y))$.
    By \ref{rtd2:item:components} and \Cref{obs:rtd_and_Thd}
    \[\rtd_2(X') = \rtd_2(T_{h_1,d_1'}(Y)) = 1 + \rtd_2(Y) \leq 1 + \rtd_2(Z) \leq \rtd_2(X).\]
    Note that $T_{1,1}(X') = K_1 \oplus X' = K_1 \oplus ((k+1) \cdot T_{h_1,d_1'}(Y))$.
    Denote by $u$ the vertex of $K_1$ in $T_{1,1}(X')$ and by $H_1,\dots,H_{k+1}$ the copies of $T_{h_1,d_1'}(Y)$ in $T_{1,1}(X')$.
    Next, let $S_1,\dots,S_k$ be a covering of $T_{1,1}(X')$.
    For every $j \in [k+1]$, there exists $i_j \in [k]$ such that $H_{j}$ contains an 
    $\{H \subseteq H_j \mid \text{$H$ connected and $S_{i_j} \cap V(H) \neq \emptyset$}\}$-rich model $\big(B_{j,x} \mid x \in V(X)\big)$ of $X$.
    By the pigeonhole principle, there exist distinct $j_1,j_2 \in [k+1]$ such that $i_{j_1}=i_{j_2}$.
    Let $i=i_{j_1}=i_{j_2}$.
    Adding a branch set $\{u\} \cup \bigcup_{x \in V(X)} B_{j_2,x}$ to the model $\big(B_{j_1,x} \mid x \in V(X)\big)$ gives a model of $T_{1,d}(X) = K_1 \oplus X$ in $T_{1,1}(X')$.
    The new branch set contains the root of $T_{1,1}(X')$ and corresponds to the root of $T_{1,d}(X)$.
    Finally, the obtained model is an $\{H \subseteq T_{1,1}(X') \mid \text{$H$ connected and $S_{i} \cap V(H) \neq \emptyset$}\}$-rich model of $T_{1,d}(X)$.

    Next, suppose that $h>1$ and that the result holds for $h-1$.
    By induction hypothesis applied to $h-1$, $d$, and $X$, there exists a positive integer $d_0'$ and a graph $X'_0$ with $\rtd_2(X'_0) \leq \rtd_2(X)$ such that
    for each covering $S_1, \dots,S_k$ of $V(T_{h-1,d'}(X'_0))$,
    there exists $i \in [k]$ such that $T_{h-1,d'_0}(X'_0)$ contains an 
    $\{H \subseteq T_{h-1,d'_0}(X'_0) \mid \text{$H$ connected and $S_i \cap V(H) \neq \emptyset$}\}$-rich model of $T_{h-1,d}(X)$
    whose branch set corresponding to the root of $T_{h-1,d}(X)$ contains the root of $T_{h-1,d'}(X'_0)$.
    Let $d' = d'_0 + (dk+1)$ and $X' = X'_0 \sqcup \big((k+1) \cdot T_{h_1,d_1'}(Y)\big)$.
    We claim that $X',d'$ satisfy the conclusion of the lemma.
    By \ref{rtd2:item:components} and \Cref{obs:rtd_and_Thd}, $\rtd_2(X') \leq \rtd_2(X)$.

    Let $S_1, \dots,S_k$ be a covering of $V(T_{h,d'}(X'))$.
    Recall that $T_{h,d'}(X') = L_{d'}(K_1 \oplus X', T_{h-1,d'}(X'),s')$ where $s'$ is the root of $T_{h-1,d'}(X')$.
    For every $x \in V(K_1 \oplus X'))$, let $H_{1,x}, \dots, H_{d',x}$ be the copies of $T_{h-1,d'}(X')$ 
    such that $T_{h,d'}(X')$ is obtained from their disjoint union with $K_1 \oplus X'$
    by identifying $x$ with the copies of $s'$ in each of $H_{1,x}, \dots, H_{d',x}$.
    Note that now $x$ is the root of $H_{j,x}$ for every $j \in [d']$ and $x \in V(K_1 \oplus X'))$.
    Since $d' \geq d'_0$ and $X'_0 \subseteq X'$, 
    for every $x \in V(K_1 \oplus X')$ and for every $j \in [d']$, there exists $i_{j,x} \in [k]$ such that 
    $H_{j,x}$ contains an 
    $\{H \subseteq T_{h-1,d'}(X') \mid \text{$H$ connected and $S_{i_{j,x}} \cap V(H) \neq \emptyset$}\}$-rich model $\mathcal{M}_{j,x}$ of $T_{h-1,d}(X)$
    whose branch set of the root of $T_{h-1,d}(X)$ contains $x$.
    Since $d' \geq dk+1$, for every $x \in V(K_1 \oplus X')$, by pigeonhole principle, there exists $i_x \in [k]$ and pairwise distinct $j_{1,x}, \dots, j_{d+1,x} \in [d']$ such that $i_{j_{\ell,x},x} = i_x$ for every $\ell \in [d+1]$.

    Let $H_1,\dots,H_{k+1}$ be the copies of $T_{h_1,d_1'}(Y)$ in $X'$.
    For every $j \in [k+1]$ and for every $\ell \in [k]$, let $S_{j,\ell} = \{x \in V(H_j) \mid i_x = \ell\}$.
    For every $j \in [k+1]$, $S_{j,1},\dots,S_{j,k}$ is a covering of $V(H_j)$, therefore, there exists $\ell_j \in [k]$ such that $H_j$ contains an $\{H \subseteq H_i \mid \text{$H$ connected and $S_{j, \ell_j} \cap V(H) \neq \emptyset$}\}$-rich model $\mathcal{M}_{j}$ of $X$.
    By pigeonhole principle, there exist distinct $j_1,j_2 \in [k+1]$ and $i \in [k]$ such that $\ell_{j_1} = \ell_{j_2}$.
    Let $i = \ell_{j_1} = \ell_{j_2}$.
    Let $u$ be the vertex of $K_1$ in $k_1 \oplus X'$.
    Adding a branch set $\{u\} \cup \bigcup \mathcal{M}_{j_2}$ to the model $\mathcal{M}_{j_1}$ gives a model of
    an $\{H \subseteq K_1 \oplus X' \mid \text{$H$ connected and there is $x \in V(H)$ with $i_x = i$}\}$-rich model 
    $\big( A_y \mid y \in V(K_1 \oplus X) \big)$ of $K_1\oplus X$ in $K_1 \oplus X'$.
    For every $y \in V(K_1 \oplus X$, let $x_y \in A_y$ be such that $i_{x_y} = i$ and let $B_y = A_y \cup \bigcup \mathcal{M}_{j_{d+1},x_y}$.
    Finally, consider the model obtained from $\big( B_y \mid y \in V(K_1 \oplus X) \big)$ by adding all the models 
    $\mathcal{M}_{j_{1,x_y},x_y}, \dots, \mathcal{M}_{j_{d,x_y},x_y}$ for every $y \in V(K_1 \oplus X)$.
    We obtain an
    $\{H \subseteq T_{h,d'}(X') \mid \text{$H$ connected and $S_i \cap V(H) \neq \emptyset$}\}$-rich model of $T_{h,d}(X)$
    whose branch set of the root of $T_{h,d}(X)$ contains the root of $T_{h,d'}(X')$.
\end{proof}

\subsection{The base case}
Recall that for all positive integers $h$ and $d$, we denote by $F_{h,d}$ the (rooted) complete $d$-ary tree of vertex-height $h$.

This first lemma is a modification of a proof in~\cite[Lemma~8]{Dujmovi2023}.

\begin{lemma}\label{lemma:bounded_tw_excluding_a_star}
    Let $d$ be a positive integer.
    Let $G$ be a connected graph, let $\mathcal{D} = \big(T,(W_x \mid x \in V(T))\big)$ be a tree decomposition of $G$, and
    let $\mathcal{F}$ be a family of connected subgraphs of $G$
    such that
    $G$ has no $\mathcal{F}$-rich model of $F_{2,d}$.
    For every $u \in V(G)$, 
    there is a set $S \subseteq V(G)$ and a layering $(P_0, \dots, P_\ell)$ of $G[S]$ with $P_0 = \{u\}$ such that
    \begin{enumerateOurAlph}
        \item $V(F) \cap S \neq \emptyset$ for every $F \in \mathcal{F}$; \label{item:bounded_tw_excluding_a_star_i}
        \item for every component $C$ of $G-S$, $N_G(V(C)) \subseteq P_i \cup P_j$ for some $i,j \in \{0, \dots, \ell\}$ and either $i=j$ or $j=i+1$; \label{item:bounded_tw_excluding_a_star_v}
        \item for every $i \in [\ell]$, $P_i$ is contained in the union of at most $d$ bags of $\mathcal{D}$. \label{item:bounded_tw_excluding_a_star_iv}
    \end{enumerateOurAlph}
\end{lemma}

\begin{proof}
    We illustrate some parts of the proof in \Cref{fig:excluding-star-tw}.
    We proceed by induction on $|V(G)|$.
    If $\mathcal{F}\vert_{G-\{u\}}$ is empty, then it suffices to take $\ell=0$ and $P_0 = \{u\}$. 
     In particular, this is the case for $|V(G)| = 1$.
    Therefore, assume $|V(G)| > 1$ and $\mathcal{F}\vert_{G-\{u\}} \neq \emptyset$.
    Let $\mathcal{F}_0$ be the family of all the connected subgraphs $A$ of $G-\{u\}$ such that $A$ contains some member of $\mathcal{F}$ and
    $V(A) \cap N_G(u) \neq \emptyset$.
    Since $\mathcal{F}\vert_{G-\{u\}} \neq \emptyset$ and $G$ is connected, $\mathcal{F}_0$ is nonempty.
    
    Observe that any collection of $d+1$ pairwise disjoint $A_1, \dots, A_{d+1} \in \mathcal{F}_0$ yields an $\mathcal{F}$-rich model of $F_{2,d}$ in $G$.
    Indeed, it suffices to take $\{u\} \cup A_{d+1}$ as the branch set corresponding to the root of $F_{2,d}$ and $A_1,\ldots,A_d$ as the branch sets of the remaining $d$ vertices of $F_{2,d}$. 
    Therefore, there is at most $d$ pairwise disjoint members of $\mathcal{F}_0$, and thus, by \Cref{lemma:helly_property_tree_decomposition} applied to $G - \{u\}$ and $\mathcal{F}_0$, there exists a set  $Z \subseteq V(G-\{u\})$ included in the union of at most $d$ bags of $\mathcal{D}$ such that $Z \cap V(F) \neq \emptyset$ for every $F \in \mathcal{F}_0$.
    Suppose that $Z$ is inclusion-wise minimal for this property.
    Note that since $\mathcal{F}_0 \neq \emptyset$, $Z$ is nonempty.

        \begin{figure}[tp]
        \centering 
        \includegraphics[scale=1]{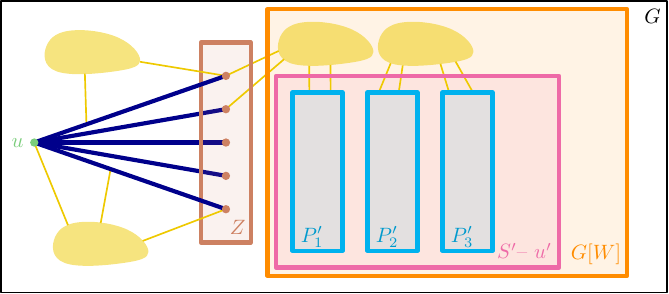} 
        \caption{
        An illustration of the proof of \Cref{lemma:bounded_tw_excluding_a_star}.
        }
        \label{fig:excluding-star-tw}
    \end{figure} 

    Let $\mathcal{C}_0$ be the family of all the components $C$ of $G-(\{u\}\cup Z)$ such that $N_G(u)\cap V(C)=\emptyset$. 
    Let $W = \bigcup_{C \in \mathcal{C}_0}V(C)$.
    Let $z \in Z$.
    By the minimality of $Z$, there exists $A_z \in \mathcal{F}_0$ with $V(A_z) \cap (Z \setminus z) = \emptyset$.
    We have $A_z - \{z\} \notin \mathcal{F}_0$, thus, $z \in V(A_z)$.
    Since $A_z \in \mathcal{F}_0$, $V(A_z) \cap N_G(u) \neq \emptyset$, and so, there is a $u$-$z$ path $Q_z$ in $G[\{u\} \cup A_z] \subset G - (Z-z)$.
    For every component $C \in \mathcal{C}_0$, we have $V(C) \cap V(Q_z) = \emptyset$, hence, $W \cap V(Q_z) = \emptyset$.

    Let $Q = \bigcup_{z \in Z} V(Q_z)$ and let $G'$ be the graph obtained from $G[W \cup Q]$ by contracting $Q$ into a single vertex $u'$.
    Note that $G'$ is a minor of $G$, and $V(G') = \{u'\} \cup W$.
    Moreover, since $Z \neq \emptyset$, we have $|V(G')| < |V(G)|$.

    For every $x \in V(T)$, let 
    \[
    W'_x =
    \begin{cases}
        W_x & \textrm{if $W_x \subseteq W$} \\
        (W_x \cap W)  \cup \{u'\} & \textrm{otherwise.}
    \end{cases}
    \]
    It follows that $\mathcal{D}' = \big(T,(W'_x \mid x \in V(T)) \big)$ is a tree decomposition of $G'$.
    By induction hypothesis applied on $G', \mathcal{D}',u',\mathcal{F}\vert_{G[W]}$, there is a set $S' \subseteq V(G')$ and a layering $(P'_0, \dots, P'_{\ell'})$ of $G'[S']$ with $P'_0=\{u'\}$ such that
    \begin{enumerateOurAlphPrim}
        \item $V(F) \cap S' \neq \emptyset$ for every $F \in \mathcal{F}\vert_{G[W]}$; \label{item:bounded_tw_excluding_a_star_i'}
        \item for every component $C$ of $G[W]-S'$, $N_G(V(C)) \subseteq P'_i \cup P'_j$ for some $i,j \in \{0, \dots, \ell'\}$ with either $i=j$ or $j=i+1$; \label{item:bounded_tw_excluding_a_star_v'}
        \item for every $i \in [\ell']$, $P'_i$ is contained in the union of at most $d$ bags of $\mathcal{D}'$. \label{item:bounded_tw_excluding_a_star_iv'}
    \end{enumerateOurAlphPrim}

    Let $S = \{u\} \cup Z \cup (S'\setminus \{u'\})$, $P_0 = \{u\}$, $P_1 = Z$, and $P_i = P'_{i-1}$ for every $i \in \{2, \dots, \ell'+1\}$ -- we set $\ell = \ell'+1$.
    We claim that $(P_0, \dots, P_\ell)$ is a layering of $G[S]$ satisfying \ref{item:bounded_tw_excluding_a_star_i}-\ref{item:bounded_tw_excluding_a_star_iv}, which will complete the proof of the lemma.

    Let $i,j \in \{0, \dots, \ell\}$ with $i<j$ and assume that there is an edge incident to a vertex in $P_i$ and a vertex in $P_j$ in $G$.
    If $i\geq 2$, then $P_i \subseteq P'_{i-1}$ and $P_j = P'_{j-1}$, which implies $|i-j|\leq 1$ since $(P_0',\dots,P_{\ell'}')$ is a layering of $G'[S']$.
    Otherwise, $i\in \{0,1\}$.
    If $i = 0$, then $j = 1$ since $u$ has no neighbors in $W$.
    If $i=1$, then $j=2$ since $N_{G'}(u') \subset P_1'$.
    It follows that $(P_0,\dots,P_\ell)$ is a layering of $G[S]$.

    Let $F \in \mathcal{F}$.
    If $V(F) \cap (\{u\} \cup Z) \neq \emptyset$, then $V(F) \cap S \neq \emptyset$.
    Otherwise, $F \subset C$ for some component of $G - (\{u\} \cup Z)$, and in particular, $C \notin \mathcal{F}_0$.
    In this case, $N_G(u) \cap V(C)$, hence, $C \in \mathcal{C}_0$, thus, $F \in \mathcal{F}\vert_{G[W]}$, and finally, $V(F) \cap S' \neq \emptyset$ by \ref{item:bounded_tw_excluding_a_star_i'}.
    This proves~\ref{item:claim_adding_apices:hitting}.

    For every component $C$ of $G-S$, either $C \cap W = \emptyset$ and so $N(V(C)) \subseteq \{u\} \cup Z = P_0 \cup P_1$,
    or $V(C) \subseteq W$, and so $C$ is a component of $G'-S'$.
    It follows that there exists $i,j \in \{0, \dots, \ell'\}$ with $|i-j| \leq 1$ such that $N_{G'}(V(C)) \subseteq P'_i \cup P'_j$
    and so $N_G(V(C)) \subseteq P_{i+1} \cup P_{j+1}$.
    This proves~\ref{item:bounded_tw_excluding_a_star_v}.

    Finally, $P_0$ is contained in one bag of $\mathcal{D}$, $P_1$ is contained in at most $d$ bags of $\mathcal{D}$ by the definition of $Z$, 
    and for every $i \in \{2, \dots, \ell\}$, $P_i$ is contained in at most $d$ bags of $\mathcal{D}$ by \ref{item:bounded_tw_excluding_a_star_iv'}.
    Therefore,~\ref{item:bounded_tw_excluding_a_star_iv} holds, which concludes the proof of the lemma.
\end{proof}

Recall that in the base case of the proof of \Cref{thm:rooted_main} (see \Cref{sec:base}), we applied the ideas required to prove that positive integer $r$ and for every path $P$, $\wcol_r(P) \leq 2+\lceil \log r \rceil$.
We illustrated these ideas in \Cref{fig:path}.
In this section, it suffices to use this result as a black box.
We state it here for reference.
\begin{lemma}[\cite{JM22}]\label{lemma:wcol_paths}
   For every positive integer $r$ and for every path $P$,
   $\wcol_r(P) \leq 2+\lceil \log r \rceil$.
\end{lemma}

\begin{lemma}\label{lemma:bounded_tw_excluding_a_star_wcol}
    Let $k,d$ be positive integers.
    For every graph $G$,
    for every tree decomposition $\mathcal{D}$ of $G$ of width at most $k-1$,
    for every family $\mathcal{F}$ of connected subgraphs of $G$,
    if $G$ has no $\mathcal{F}$-rich model of $F_{2,d}$, then
    there is a set $S \subseteq V(G)$ such that
    \begin{enumerateOurAlph}
        \item $V(F) \cap S \neq \emptyset$ for every $F \in \mathcal{F}$;\label{lemma:bounded_tw_excluding_a_star_wcol:i}
        \item for every component $C$ of $G-S$, $N_G(V(C))$ is contained in the union of at most $2d$ bags of~$\mathcal{D}$; 
        \label{lemma:bounded_tw_excluding_a_star_wcol:ii}
        \item $\wcol_r(G,S) \leq 6dk \cdot \log r$ for every integer $r$ with $r \geq 2$.\label{lemma:bounded_tw_excluding_a_star_wcol:iii}
    \end{enumerateOurAlph}
\end{lemma}

\begin{proof}
    Let $G$ be a graph, let $\mathcal{D}$ be a tree decomposition of $G$ of width at most $k-1$, let $\mathcal{F}$ be a family of connected subgraphs of $G$, and suppose that $G$ has no $\mathcal{F}$-rich model of $F_{2,d}$.
    Let $r$ be an integer with $r \geq 2$. 
    By considering the components of $G$ independently, we can assume that $G$ is connected.
    Let $u$ be an arbitrary vertex in $G$.
    \Cref{lemma:bounded_tw_excluding_a_star} applied to $G$, $\mathcal{D}$, $\mathcal{F}$, and $u$ gives
    $S \subseteq V(G)$, a layering $(P_0, \dots, P_\ell)$ of $G[S]$ with $P_0 = \{u\}$ such that
    \begin{enumerate}[label={\normalfont\ref{lemma:bounded_tw_excluding_a_star}.(\makebox[\mywidth]{\alph*})}]
        \item $V(F) \cap S \neq \emptyset$ for every $F \in \mathcal{F}$; \label{item:bounded_tw_excluding_a_star_i-call}
        \item for every component $C$ of $G-S$, $N_G(V(C)) \subseteq P_i \cup P_j$ for some $i,j \in \{0, \dots, \ell\}$ and either $i=j$ or $j=i+1$; \label{item:bounded_tw_excluding_a_star_v-call}
        \item for every $i \in [\ell]$, $P_i$ is contained in the union of at most $d$ bags of $\mathcal{D}$. \label{item:bounded_tw_excluding_a_star_iv-call}  
    \end{enumerate}

    Note that \ref{lemma:bounded_tw_excluding_a_star_wcol:i} holds by \ref{item:bounded_tw_excluding_a_star_i-call} and \ref{lemma:bounded_tw_excluding_a_star_wcol:ii} holds by~\ref{item:bounded_tw_excluding_a_star_v-call} and~\ref{item:bounded_tw_excluding_a_star_iv-call}.
    In order to conclude the proof, it suffices to show~\ref{lemma:bounded_tw_excluding_a_star_wcol:iii}.
    
    For convenience, let $P_{\ell+1} = \emptyset$.
    Consider the path $Q$ with $V(Q) = \{0,\dots,\ell+1\}$ where two numbers are connected by an edge whenever they are consecutive.
    Let $\sigma' = i_0\dots i_\ell$ be an ordering of $\{0,\dots,\ell+1\}$ given by \Cref{lemma:wcol_paths}, that is, $\wcol_r(Q,\sigma') \leq 1+\lceil \log r \rceil \leq 3 \log r$.
    For each $i \in \{0,\dots,\ell+1\}$, let $\sigma_i$ be an arbitrary ordering of $P_i$.
    Let $\sigma$ be the concatenation of $\sigma_{i_0}\dots \sigma_{i_{\ell+1}}$ in this order.

    Let $u \in V(G)$.
    To conclude the claim, we argue that
            \[|\WReach_r[G,S,\sigma,u]| \leq 6dk \cdot \log r.\]
    Let $i_u \in V(T)$ be such that if $u \in S$, then $u \in P_{i_u}$, and otherwise, $i_u \in \{0,\dots,\ell\}$ is the least value such that $P_{i_u}$ intersects $N_G(C)$, where $C$ is the component of $u$ in $G-S$.
    Let $A = \WReach_r[Q,\sigma',i_u] \cup \WReach_r[Q,\sigma',i_u+1]$.
    In particular, $|A| \leq 2\cdot \wcol_r(Q,\sigma') \leq 2 \cdot 3 \log r$
    By~\ref{item:bounded_tw_excluding_a_star_v-call}, 
        \[\WReach[G, S, \sigma,u] \subset \bigcup_{j \in A} P_j.\]
    By~\ref{item:bounded_tw_excluding_a_star_iv-call}, for every $j \in \{0,\dots,\ell+1\}$, $P_j$ is contained in the union of at most $d$ bags of $\mathcal{D}$ and since the width of $\mathcal{D}$ is at most $k-1$, we have $|P_j| \leq dk$.
    It follows that
        \[|\WReach[G, S, \sigma,u]| \leq |A| \cdot dk \leq 6dk\log r. \qedhere\]
\end{proof}

\begin{lemma}\label{lem:bounded_tw_new_base_case_0}
    Let $k,h,d$ be positive integers with $h \geq 2$.
    There is an integer $\cc_0(h,d)$ such that
    for every graph $G$, for every tree decomposition $\mathcal{D}$ of $G$ of width at most $k-1$,
    for every family $\mathcal{F}$ of connected subgraphs of $G$,
    if $G$ has no $\mathcal{F}$-rich model of $F_{h,d}$, then
    there is a set $S \subseteq V(G)$
    such that
    \begin{enumerateOurAlph}
        \item $V(F) \cap S \neq \emptyset$ for every $F \in \mathcal{F}$,\label{lem:bounded_tw_new_base_case_0:i}
        \item for every component $C$ of $G-S$, $N_G(V(C))$ is contained in the union of at most $2d(h-1) + 2\binom{h-1}{2}$ bags of $\mathcal{D}$;\label{lem:bounded_tw_new_base_case_0:item:nb_comp_interface}
        \item $\wcol_r(G,S) \leq \cc_0(h,d)k \cdot \log r$ for every integer $r$ with $r \geq 2$.\label{lem:bounded_tw_new_base_case_0:iii}
    \end{enumerateOurAlph}
\end{lemma}

\begin{proof}
    We proceed by induction on $h$.
    For $h=2$, the result holds by \Cref{lemma:bounded_tw_excluding_a_star_wcol} setting $\cc_0(h,d)=6dk$, since $2d(2-1)+2\binom{2-1}{2}=2d$.
    Next, assume $h \geq 3$ and that $\cc_0(h,d)$ witnesses the assertion for $h-1$.
    Let $\cc_0(h,d) = 6d + \cc_0(h-1,d+1)$.

    Let $G$ be a graph and let $\mathcal{D}$ be a tree decomposition of $G$ of width at most $k-1$.
    Let $\mathcal{F}$ be a family of connected subgraphs of $G$.
    Let $r$ be an integer with $r \geq 2$.
    Let $\mathcal{F}'$ be the family of all the connected subgraphs $H$ of $G$ such that $H$ contains an $\mathcal{F}\vert_H$-rich model of $F_{h-1,d+1}$.
    We claim that there is no $\mathcal{F}'$-rich model of $F_{2,d}$ in $G$.
    Suppose to the contrary that $\big(B_x \mid x \in V(F_{2,d})\big)$ is such a model.
    Let $s$ be the root of $F_{2,d}$ and let $s'$ be the root of $F_{h-1,d}$.
    For every $x \in V(F_{2,d}) \setminus \{s\}$,
    by \Cref{lemma:rooting_a_model_of_Fhd},
    there is an $\mathcal{F}$-rich model $\big(C_y \mid y \in V(F_{h-1,d})\big)$ of $F_{h-1,d}$ in $G[B_x]$ such that $C_{s'}$ contains a vertex of $N_G(B_s) \cap B_x$.
    The union of these models together with $B_s$ yields an $\mathcal{F}$-rich model of $F_{h,d}$ in $G$, which is a contradiction.
    Note that this is exactly the same argument as in \Cref{lem:new_base_case}, see \Cref{fig:forest-models}.

    Since $G$ has no $\mathcal{F}'$-rich model of $F_{2,d}$, by \Cref{lemma:bounded_tw_excluding_a_star_wcol}, there is a set $S_0 \subseteq V(G)$ such that
    \begin{enumerate}[label={\normalfont\ref{lemma:bounded_tw_excluding_a_star_wcol}.(\makebox[\mywidth]{\alph*})}]
        \item for every $F \in \mathcal{F}'$, $V(F) \cap S_0 \neq \emptyset$; \label{lemma:bounded_tw_excluding_a_star_wcol:item:hitting}
        \item for every component $C$ of $G-S_0$, $N_G(V(C))$ is contained in the union of at most $2d$ bags of~$\mathcal{D}$;\label{lemma:bounded_tw_excluding_a_star_wcol:item:connected}
        \item $\wcol_r(G,S_0) \leq 6dk \cdot \log r$.\label{lemma:bounded_tw_excluding_a_star_wcol:item:wcol}
    \end{enumerate}
    Let $C$ be a component of $G - S_0$. 
    Let $\mathcal{D}_C$ be $\mathcal{D}$ restricted to $C$.
    By~\ref{lemma:bounded_tw_excluding_a_star_wcol:item:hitting}, $C \notin \mathcal{F}'$, and so, $C$ has no $\mathcal{F}\vert_C$-rich model of $F_{h-1,d+1}$.
    Therefore, by induction hypothesis, there is a set $S_C \subseteq V(C)$ such that
    \begin{enumerateOurAlphPrim}
        \item $V(F) \cap S_C \neq \emptyset$ for every $F \in \mathcal{F}\vert_C$; \label{item:lem:bounded_tw_new_base_case_i}
        \item for every component $C'$ of $C-S_C$, $N_{C}(V(C'))$ is contained in the union of at most $2(d+1)(h-2)+2\binom{h-2}{2} = 2d(h-1)+2\binom{h-1}{2}-2d$ bags of $\mathcal{D}_C$; \label{item:lem:bounded_tw_new_base_case_ii}
        \item $\wcol_r(C,S_C) \leq \cc_0(h-1,d+1)k \cdot \log r$.\label{item:lem:bounded_tw_new_base_case_0_iii}
    \end{enumerateOurAlphPrim}
    Let $\mathcal{C}$ be the family of components of $G-S_0$ and let
    \[
    S = S_0 \cup \bigcup_{C \in \mathcal{C}} S_C.
    \]
    We claim that \ref{lem:bounded_tw_new_base_case_0:i}-\ref{lem:bounded_tw_new_base_case_0:iii} hold.
    Let $F \in \mathcal{F}$.
    If $V(F) \cap S_0 = \emptyset$, then $V(F) \subset V(C)$ for some component $C$ of $G - S_0$.
    In particular, $F \in \mathcal{F}\vert_C$, and thus, by~\ref{item:lem:bounded_tw_new_base_case_i}, $V(F) \cap S_C \neq \emptyset$, which proves~\ref{lem:bounded_tw_new_base_case_0:i}.
    For every component $C'$ of $G - S$, there exists a component $C$ of $G-S_0$ such that $C' \subset C$ and by~\ref{lemma:bounded_tw_excluding_a_star_wcol:item:connected} and \ref{item:lem:bounded_tw_new_base_case_ii}, $N_G(C')$ is contained in at most $2d + \left( 2d(h-1)+2\binom{h-1}{2}-2d \right) = 2d(h-1)+2\binom{h-1}{2}$ bags of $\mathcal{D}$, which implies \ref{lem:bounded_tw_new_base_case_0:item:nb_comp_interface}.
    
    The following sequence of inequalities concludes the proof of~\ref{lem:bounded_tw_new_base_case_0:iii} and the lemma:
    \begin{align*}
    \wcol_r(G,S) &\leq \wcol_r(G,S_0) + \wcol_r\left(G-S_0,\bigcup_{C \in \mathcal{C}} S_C \right)&&\textrm{by~\Cref{obs:wcol_union}}\\ 
    &\leq \wcol_r(G,S_0) + \max_{C \in \mathcal{C}} \wcol_r(C, S_C)&&\textrm{by~\Cref{obs:wcol_components2}}\\
    &\leq 6dk\cdot \log r + \cc_0(h-1,d+1)k \cdot \log r&&\textrm{by~\ref{lemma:bounded_tw_excluding_a_star_wcol:item:wcol} and \ref{item:lem:bounded_tw_new_base_case_0_iii}}\\ 
    &= \cc_0(h,d)k \cdot \log r. &&\qedhere
    \end{align*}  
\end{proof}

We now show that using \Cref{lemma:increase_X_to_have_small_interfaces}, 
the constant $2d(h-1)+2\binom{h-1}{2}$ in \ref{lem:bounded_tw_new_base_case_0:item:nb_comp_interface} in \Cref{lem:bounded_tw_new_base_case_0} can be reduced to two.

\begin{lemma}\label{lem:bounded_tw_new_base_case}
    Let $k,h,d$ be positive integers with $h \geq 2$.
    For every graph $G$ with $\tw(G)<k$,
    for every family $\mathcal{F}$ of connected subgraphs of $G$,
    if $G$ has no $\mathcal{F}$-rich model of $F_{h,d}$, then
    there is a set $S \subseteq V(G)$
    such that
    \begin{enumerateOurAlph}
        \item $V(F) \cap S \neq \emptyset$ for every $F \in \mathcal{F}$;\label{item:i:lem:bounded_tw_new_base_case}
        \item for every component $C$ of $G-S$, $N_G(V(C))$ intersects at most two components of $G-V(C)$;\label{item:ii:lem:bounded_tw_new_base_case}
        \item $\wcol_r(G,S) \leq \left(\cc_0(h,d)+4dh^2\right)k \cdot \log r$ for every integer $r$ with $r \geq 2$, where $\cc_0(h,d)$ is the constant from \Cref{lem:bounded_tw_new_base_case_0}. \label{item:iii:lem:bounded_tw_new_base_case} 
    \end{enumerateOurAlph}
\end{lemma}

\begin{proof}
    Let $G$ be a graph and let $\mathcal{D}$ be a tree decomposition of $G$ of width at most $k-1$.
    We may assume $G$ is connected.
    By \Cref{lemma:natural_tree_decomposition},
    we can assume that $\mathcal{D}$ is a natural tree decomposition of $G$.
    By \Cref{lemma:bounded_tw_excluding_a_star_wcol} applies to $G$, $\mathcal{D}$, and $\mathcal{F}$, there is a set $S_0 \subseteq V(G)$
    such that
    \begin{enumerate}[label={\normalfont\ref{lemma:bounded_tw_excluding_a_star_wcol}.(\makebox[\mywidth]{\alph*})}]
        \item $V(F) \cap S_0 \neq \emptyset$ for every $F \in \mathcal{F}$,\label{lem:bounded_tw_new_base_case_0:i_local}
        \item for every component $C$ of $G-S_0$, $N_G(V(C))$ is contained in the union of at most $2d(h-1)+2\binom{h-1}{2}$ bags of $\mathcal{D}$;\label{lem:bounded_tw_new_base_case_0:item:nb_comp_interface_local}
        \item $\wcol_r(G,S_0) \leq \cc_0(h,d)k \cdot \log r$ for every integer $r$ with $r \geq 2$.\label{lem:bounded_tw_new_base_case_0:iii_local}
    \end{enumerate}   
    
    Let $C$ be a component of $G - S_0$.
    By~\ref{lem:bounded_tw_new_base_case_0:item:nb_comp_interface_local} and \Cref{lemma:increase_X_to_have_small_interfaces}, there exists a family $\mathcal{B}_C$ of at most $2\left(2d(h-1)+2\binom{h-1}{2}\right)-1 \leq 4dh^2$ bags of $\mathcal{D}$ such that $N_G(V(C)) \subset \bigcup \mathcal{B}_C$ and for every component $C'$ of $G-\bigcup \mathcal{B}_C$, $N_G(V(C'))$ intersects at most two components of $G-V(C')$.
    Let $S_C = V(C) \cap \bigcup \mathcal{B}_C$.
    Then, for every component $C'$ of $G-(S_0 \cup S_C)$ intersecting $V(C)$,
    $C$ is a component of $G-\bigcup \mathcal{B}_C$, and so $N_G(V(C'))$ intersects at most two component of $G-V(C')$.

    Let $\mathcal{C}$ be the family of all the components of $G - S_0$ and let
    \[S = S_0 \cup \bigcup_{C \in \mathcal{C}} S_C.\]
    Item~\ref{item:ii:lem:bounded_tw_new_base_case} follows from the previous considerations and~\ref{item:i:lem:bounded_tw_new_base_case} follows directly from~\ref{lem:bounded_tw_new_base_case_0:i_local}.
    Now, it suffices to justify~\ref{item:iii:lem:bounded_tw_new_base_case}.
    To this end, let $r$ be an integer with $r \geq 2$.
    Then,
    \begin{align*}
    \wcol_r(G,S) &\leq \wcol_r(G,S_0) + \wcol_r\left(G-S_0,\bigcup_{C \in \mathcal{C}} S_C \right)&&\textrm{by~\Cref{obs:wcol_union}}\\ 
    &\leq \wcol_r(G,S_0) + \max_{C \in \mathcal{C}} \wcol_r(C, S_C)&&\textrm{by~\Cref{obs:wcol_components2}}\\
    &\leq \wcol_r(G,S_0) + \max_{C \in \mathcal{C}} |S_C|\\
    &\leq \cc_0(h,d)k \cdot \log r + 4dh^2 k&&\textrm{by~\ref{lem:bounded_tw_new_base_case_0:iii_local}}\\ 
    &\leq \left(\cc_0(h,d) + 4dh^2\right)k \cdot \log r. &&\qedhere
    \end{align*}  
\end{proof}

\subsection{Induction}\label{subsec:bounded_tw_induction}
We can now prove \Cref{thm:rooted_main_bounded_tw} in the following stronger version,
which is very similar to \Cref{thm:main}, with a slightly relaxed condition~\ref{thm:main:item:connected}.
Note that the following proof and the one of \Cref{thm:main} largely overlap, but since there is no non-artificial way to merge them, 
we elect to give full proofs of both theorems.

\begin{theorem}\label{thm:bounded_tw_main}
    Let $t$ be positive integers with $t \geq 2$.
    Let $X$ be a graph with $\rtd_2(X) \leq t$.
    There exist an integer $\cc(t,X)$
    such that
    for every integer $k$,
    for every graph $G$ with $\tw(G)<k$, 
    for every family $\mathcal{F}$ of connected subgraphs of $G$,
    if $G$ has no $\mathcal{F}$-rich model of $X$, then
    there exists $S \subseteq V(G)$ such that
    \begin{enumerateOurAlphCapital}
        \item $V(F) \cap S \neq \emptyset$ for every $F \in \mathcal{F}$;\label{thm:bounded_tw_main:item:hit}
        \item for every component $C$ of $G-S$, $N_G(V(C))$ intersects at most two components of $G-V(C)$; \label{thm:bounded_tw_main:item:connected}
        \item $\wcol_r(G,S) \leq \cc(t,X)k \cdot r^{t-2} \log r$ for every integer $r$ with $r \geq 2$. \label{thm:bounded_tw_main:item:wcol}
    \end{enumerateOurAlphCapital}
\end{theorem}

\begin{proof}
    We proceed by induction on $t$.
    When $t = 2$, by \Cref{obs:char:rtd_2=2}, $X$ is a forest.
    Let $h,d$ be positive integers such that $X \subset F_{h,d}$,
    and let $\cc_0(h,d)$ be the constant given by \Cref{lem:bounded_tw_new_base_case}.
    The assertion with $\cc(t,X) = \cc_0(h,d) + 4dh^2$ follows by applying \Cref{lem:bounded_tw_new_base_case}.
    Next, let $t \geq 3$, and assume that the result holds for $t-1$. 
    We refer to this property as the \emph{main induction hypothesis}.
    
    \begin{claim1}\label{claim:bounded_tw_adding_apices}
        Let $Y$ be a graph with $\rtd_2(Y) \leq t-1$.
        There is an integer $\cc_1(t,Y)$ such that
        for every positive integer $k$,
        for every connected graph $G$ with $\tw(G)<k$,
        for every nonempty set $U$ of vertices of $G$
        with $|U| \leq 2$,
        for every family $\mathcal{F}$ of connected subgraphs of $G$,
        if $G$ has no $\mathcal{F}$-rich model of $K_1 \oplus Y$, 
        then
        there exist $S \subseteq V(G)$, a tree $T$ rooted in $s \in V(T)$, and a tree partition $\big(T,(P_x \mid x \in V(T))\big)$
        of $G[S]$
        with $P_s = U$
        such that
        \begin{enumerateOurAlph}
            \item $V(F) \cap S \neq \emptyset$ for every $F \in \mathcal{F}$;\label{item:bounded_tw_claim_adding_apices_i}
            \item for every component $C$ of $G-S$, $N_G(V(C))$ intersects at most four components of $G-V(C)$;\label{item:bounded_tw_claim_adding_apices_ii}
            \item for every component $C$ of $G-S$, $N_{G}(V(C)) \subset P_x \cup P_y$ for some $x,y \in V(T)$ with either $x=y$ or $xy$ is an edge in $T$;\label{item:bounded_tw_claim_adding_apices_iii}
            \item for every $x \in V(T)$, \[\wcol_r\left(G_x, P_x\right) \leq \cc_1(t,Y) k \cdot r^{t-3} \log r\] for every integer $r$ with $r \geq 2$,
            where, for $T_x$ being the subtree of $T$ rooted in $x$,
            $G_x$ is the subgraph of $G$ induced by $U_x = \bigcup_{y \in V(T_x)} P_y$ and the vertex sets of all the components of $G-S$ having a neighbor in $U_x$. \label{item:bounded_tw_claim_adding_apices_vii}
        \end{enumerateOurAlph}
    \end{claim1}

    \begin{proofclaim}
        Let $X' =  K_1 \sqcup Y$.
        By \Cref{lemma:rooting_X_in_X'} applied to $k=2$ and $X'$, there exists a graph $X''$ with $\rtd_2(X'')=\rtd_2(X')=t$
        such that for every covering $S_1,S_2$ of $V(X'')$ there exists $i \in \{1,2\}$ such that $X''$ contains an 
        $\{H \subseteq X'' \mid \text{$H$ connected and } S_i \cap V(H) \neq \emptyset\}$-rich model of $X'$.
        Let $\cc_1(t,Y) = 2 + \cc(t-1,X'')$.
        
        We proceed by induction on $|V(G)|$.
        If $\mathcal{F}\vert_{G-U} = \emptyset$, then since $\cc_1(t,Y) \geq |U|$, the result holds.
        Now suppose that $\mathcal{F}\vert_{G-U} \neq \emptyset$ and in particular, $|V(G)-U| > 1$.

        Let $\mathcal{F}'$ be the family of all the connected subgraphs $H$ of $G-U$ such that $U \cap N_G(V(H)) \neq \emptyset$ and
        $F \subseteq H$ for some $F \in \mathcal{F}$.
        We argue that there is no $\mathcal{F}'$-rich model of $X''$ in $G-U$.
        Suppose to the contrary that there is an $\mathcal{F}'$-rich model $\big(B_x \mid x \in V(X'')\big)$ of $X''$ in $G-U$.
        For each $u \in U$, let $S_u = \{x \in V(X'') \mid u \in N_G(B_x)\}$.
        Since the model is $\mathcal{F}'$-rich, $\{S_u\}_{u \in U}$ is a covering of $V(X'')$.
        Therefore, there exists $u \in U$ such that $X''$ contains an 
        $\{H \subseteq X'' \mid \text{$H$ connected and } S_u \cap V(H) \neq \emptyset\}$-rich model $\mathcal{M}$ of $X'$.
        Moreover, by definition, each branch set of $\mathcal{M}$ contains a member of $\mathcal{F}$, and furthermore, $\mathcal{M}$ is a model of $X'$ in $G$.
        Recall that $X' = K_1 \sqcup Y$.
        If $v$ is the vertex of $K_1$ is $X'$, then the model obtained from $\mathcal{M}$ by replacing the branch set $C_v$ corresponding to $v$ by $C_v \cup \{u\}$ yields an $\mathcal{F}$-rich model of $K_1 \oplus Y$ in $G$, which is a contradiction.
        This proves that there is no $\mathcal{F}'$-rich model of $X'$ in $G-\{u\}$.

        Since $\rtd_2(X'') = \rtd_2(X') \leq \max \{\rtd_2(Y),1\} \leq t-1$, by the main induction hypothesis applied to $X'', G - U$, and $\mathcal{F}'$
        there exists a set $S_0 \subseteq V(G-U)$ such that   \begin{enumerateOurAlphCapitalPrim}
            \item $V(F) \cap S_0 \neq \emptyset$ for every $F \in \mathcal{F}'$;\label{thm:bounded_tw_main:item:hit'}
            \item for every component $C$ of $(G-U) - S_0$, $N_{G-U}(V(C))$ intersects at most two components of $(G-U)-V(C)$; \label{thm:bounded_tw_main:item:connected'}
            \item $\wcol_r(G-U,S) \leq \cc(t-1,X'')k \cdot r^{t-3} \log r$ for every integer $r$ with $r \geq 2$. \label{thm:bounded_tw_main:item:wcol'}
        \end{enumerateOurAlphCapitalPrim}
        Note that $S_0 \neq \emptyset$ since $\mathcal{F}' \neq \emptyset$.

        Let $\mathcal{C}_1$ be the family of all the components $C$ of $(G-U)- S_0$ such that $N_G(U)\cap V(C)=\emptyset$. 
        Consider $C \in \mathcal{C}_1$.
        Let $G_C$ be the graph obtained from $G[V(C) \cup N_G(V(C))]$ by contracting each component $C'$ of $(G-U) - V(C)$ into a single vertex.
        Let $U_{C}$ be the set of all the vertices resulting from these contractions.
        Note that $U_C$ is not empty since $G$ is connected, and $|U_C| \leq 2$ by \ref{thm:bounded_tw_main:item:connected'}.
        Since $N_G(U)\cap V(C)=\emptyset$,  $|V(G_C)| < |V(G)|$.
        Since $G_C$ is a minor of $G$, $G_C$ has no $\mathcal{F}\vert_C$-rich model of $K_1 \oplus Y$.
        By induction hypothesis applied to $G_C$, $u_C$, and $\mathcal{F}\vert_C$, there exist $S_C \subseteq V(G_C)$, a tree $T_C$ rooted in $s_C \in V(T_C)$, and a tree partition $(T_C,(P_{C,x} \mid x \in V(T_C)))$ of $S_C$ in $G_C$ with $P_{C,s_C} = U$
        such that
        \begin{enumerateOurAlphPrimPrim}
            \item $V(F) \cap S_C \neq \emptyset$ for every $F \in \mathcal{F}\vert_{C}$; \label{item:bounded_tw_claim_adding_apices_i-call-call}
            \item for every component $C'$ of $G_C-S_C$, $N_G(V(C'))$ intersects at most four components of $G_C-V(C')$;\label{item:bounded_tw_claim_adding_apices_ii-call-call}
            \item for every component $C'$ of $G_C-S_C$, $N_{G_C}(V(C')) \subset P_{C,x} \cup P_{C,y}$ for some $x,y \in V(T_C)$ and either $x=y$ or $xy$ is an edge in $T_C$; \label{item:bounded_tw_claim_adding_apices_iii-call-call}
            \item for every $x \in V(T_C)$, \[\wcol_r\left(G_{C,x}, P_{C,x}\right) \leq \cc_1(t,Y) k \cdot r^{t-3} \log r\] for every integer $r$ with $r \geq 2$,
            where, for $T_{C,x}$ being the subtree of $T_C$ rooted in $x$,
            $G_{C,x}$ is the subgraph of $G_C$ induced by $U_{C,x} = \bigcup_{y \in V(T_{C,x})} P_{C,y}$ and the vertex sets of all the components of $G_C-S_C$ having a neighbor in $U_{C,x}$. \label{item:bounded_tw_claim_adding_apices_vii-call-call}
        \end{enumerateOurAlphPrimPrim}

        Let
        \[
            S = U \cup S_0 \cup \bigcup_{C \in \mathcal{C}_1} (S_C \setminus U_C).
        \]
        Let $T$ be obtained from the disjoint union of $\{T_C \mid C\in\mathcal{C}_1\}$ by identifying the vertices $\set{s_{C}\mid C\in\mathcal{C}_1}$ into a new vertex $s'$ and by adding a new vertex $s$ adjacent to $s'$ in $T$. 
        Let $P_s=U$, $P_{s'}=S_0$, and
        for each $C\in\mathcal{C}_1$, 
        $x \in V(T_{C} \setminus \{s_{C}\})$, 
        let $P_x = P_{C,x}$.

        In order to conclude, we argue that $\big(T,(P_x \mid x \in V(T))\big)$ is a tree partition of $G[S]$ 
        and \ref{item:bounded_tw_claim_adding_apices_i}-\ref{item:bounded_tw_claim_adding_apices_vii} hold.

        Since for every $C \in \mathcal{C}_1$, $U \cap N_G(V(C)) = \emptyset$, every edge in $G[S]$ containing a vertex in $U$ has another endpoint in $U \cup S_0 = P_s \cup P_{s'}$.
        Consider an edge $vw$ in $G[S]$ such that $v \in S_0$ and $w \in S_C$ for some $C \in \mathcal{C}_1$.
        Since $\big(T_C,(P_{C,x} \mid x \in V(T_C))\big)$ is a tree partition of $G_C[S_C]$ with $P_{C,s_C} = U_C$ and $S_0 \subset V((G-U)-V(C))$, we conclude that $w \in P_x$ for some $x \in V(T_C)$ such that $s'x$ is an edge in $T$.
        Finally, for every edge $vw$ of $G[S]$ with $v,w \not\in U \cup S_0$, $vw$ is an edge of $G[S_C \setminus U_C]$ for some component $C \in \mathcal{C}_1$,
        and so $v \in P_{C,x}$ and $w \in P_{C,y}$ for adjacent or identical vertices $x,y$ of $T_C$.
        Then $v \in P_x$ and $w \in P_y$.
        It follows that $\big(T,(P_{x} \mid x \in V(T))\big)$ is a tree partition of $G[S]$.

        Let $F \in \mathcal{F}$.
        If $V(F) \cap (U \cup S_0) \neq \emptyset$, then $V(F) \cap S \neq \emptyset$.
        Otherwise, $F \subseteq G - U$ and $V(F) \cap S_0 = \emptyset$, and therefore by \ref{thm:bounded_tw_main:item:hit'}, $F \notin \mathcal{F}'$, so in particular, $U \cap N_G(V(F)) = \emptyset$.
        In this case, there is a component $C \in \mathcal{C}_1$ such that $F \in \mathcal{F}\vert_C$, thus, $V(F) \cap S_C \neq \emptyset$ by \ref{item:bounded_tw_claim_adding_apices_i-call-call}.
        This proves~\ref{item:bounded_tw_claim_adding_apices_i}.

        Consider a component $C'$ of $G - S$.
        If $C' \subseteq C$ for some $C \in \mathcal{C}_1$, then by~\ref{item:bounded_tw_claim_adding_apices_ii-call-call}, it follows that $N_{G_C}(V(C'))$ intersects at most four components of $G_C-V(C')$, 
        and so $N_{G}(V(C'))$ intersects at most four components of $G-V(C')$.
        Otherwise, $C'$ is a component of $(G - U) - S_0$ such that $N_G(U) \cap V(C') \neq \emptyset$.
        By~\ref{thm:bounded_tw_main:item:connected'}, $N_{G-U}(V(C'))$ intersects at most two components of $(G-U)-V(C')$, and therefore, $N_G(V(C'))$ intersects at most four components of $G-V(C')$.
        This proves~\ref{item:bounded_tw_claim_adding_apices_ii}.
        
        
        For every component $C'$ of $G-S$, either $N_G(V(C')) \subseteq U \cup S_0 = P_s \cup P_{s'}$, or $C' \subseteq C$ for some $C \in \mathcal{C}_1$.
        In the latter case, $C'$ is a component of $G_C - S_C$, and $N_G(U) \cap V(C) = \emptyset$.
        By~\ref{item:bounded_tw_claim_adding_apices_iii-call-call}, there is $x,y \in V(T_C)$ such that $N_{G_C}(V(C')) \subseteq P_{C,x} \cup P_{C,y}$, and thus, $N_G(V(C')) \subseteq P_x \cup P_y$.
        This proves \ref{item:bounded_tw_claim_adding_apices_iii}.

        Finally, we argue~\ref{item:bounded_tw_claim_adding_apices_vii}.
        Let $r$ be an integer with $r \geq 2$ and let $x \in V(T)$.
        For $x = s$, $|P_{s}| = |U| \leq 2$, thus the assertion is clear.
        For $x = s'$, we have that $G_{s'}$ is a union of components of $G-U$.
        By~\ref{thm:bounded_tw_main:item:wcol'}, 
            \[\wcol(G_{s'},S_0) = \wcol_r(G-U, S_0) \leq \cc(t-1,X'')k \cdot r^{t-3} \log r \leq \cc_1(t,Y)k \cdot r^{t-3} \log r.\]
        For $x \in V(T_C - \{s_C\})$ for some $C \in \mathcal{C}_1$, we have $G_x = G_{C,x}$, thus, the asserted inequality follows from~\ref{item:bounded_tw_claim_adding_apices_vii-call-call}.
        This ends the proof of the claim.
    \end{proofclaim}

    \Cref{claim:bounded_tw_adding_apices} yields the following less technical statement.
    \begin{claim1}\label{claim:bounded_tw_adding_apices_wcol}
        Let $Y$ be a graph with $\rtd_2(Y) \leq t-1$.
        There are integers $\cc_2(t,Y)$ such that
        for every positive integer $k$,
        for every graph $G$ with $\tw(G)<k$,
        for every family $\mathcal{F}$ of connected subgraphs of $G$,
        if $G$ has no $\mathcal{F}$-rich model of $K_1 \oplus Y$, then
        there exist $S \subseteq V(G)$ such that
        \begin{enumerateOurAlph}
            \item $V(F) \cap S \neq \emptyset$ for every $F \in \mathcal{F}$; \label{item:claim:bounded_tw_adding_apices_wcol:i}
            \item for every component $C$ of $G-S$, $N_G(V(C))$ intersects at most four components of $G-V(C)$;\label{item:claim:bounded_tw_adding_apices_wcol:ii}
            \item $\wcol_r(G,S) \leq \cc_2(t,Y)k \cdot r^{t-2} \log r$ for every integer $r$ with $r \geq 2$. \label{item:claim:bounded_tw_adding_apices_wcol:iii}
        \end{enumerateOurAlph}
    \end{claim1}

    \begin{proofclaim}
        Let $\cc_2(t,Y) = 2 \cdot \cc_1(t,Y)$ and let $r$ be an integer with $r \geq 2$.
        We apply \Cref{claim:bounded_tw_adding_apices} with an arbitrary singleton $\{u\}$ in $G$ to obtain $S \subseteq V(G)$, a tree $T$ rooted in $s \in V(T)$, and a tree partition $(T,(P_x \mid x \in V(T)))$
        of $S$ in $G$ with $P_s = \{u\}$
        such that
        \begin{enumerate}[label={\normalfont\ref{claim:bounded_tw_adding_apices}.(\makebox[\mywidth]{\alph*})}]
            \item $V(F) \cap S \neq \emptyset$ for every $F \in \mathcal{F}$;\label{item:bounded_tw_claim_adding_apices_i-called-in-next}
            \item for every component $C$ of $G-S$, $N_G(V(C))$ intersects at most four components of $G-V(C)$;\label{item:bounded_tw_claim_adding_apices_ii-called-in-next}
            \item for every component $C$ of $G-S$, $N_{G}(V(C)) \subset P_x \cup P_y$ for some $x,y \in V(T)$ with either $x=y$ or $xy$ is an edge in $T$;\label{item:bounded_tw_claim_adding_apices_iii-called-in-next}
            \item for every $x \in V(T)$, \[\wcol_r\left(G_x, P_x\right) \leq \cc_1(t,Y) k \cdot r^{t-3} \log r\] for every integer $r$ with $r \geq 2$,
            for $T_x$ being the subtree of $T$ rooted in $x$,
            $G_x$ is the subgraph of $G$ induced by $\bigcup_{y \in V(T_x)} P_y$ and the vertex sets of all the components of $G-S$ having a neighbor in $U_x$. \label{item:bounded_tw_claim_adding_apices_vii-called-in-next}
        \end{enumerate}
        Items \ref{item:claim:bounded_tw_adding_apices_wcol:i} and \ref{item:claim:bounded_tw_adding_apices_wcol:ii} hold by \ref{item:bounded_tw_claim_adding_apices_i-called-in-next} and \ref{item:bounded_tw_claim_adding_apices_ii-called-in-next} respectively.
        It suffices to prove \ref{item:claim:bounded_tw_adding_apices_wcol:iii}.

        For each $x \in V(T)$, let $\sigma_x$ be an ordering of $P_x$ witnessing \ref{item:bounded_tw_claim_adding_apices_vii-called-in-next} and let $\sigma' = (x_1, \dots, x_{|V(T)|})$ be an elimination ordering of $T$.
        Finally, let $\sigma$ be the concatenation of $\sigma_{x_1},\dots,\sigma_{x_{|V(T)|}}$ in this order.

        Let $u \in V(G)$.
        To conclude the claim, we argue that
        \[
            |\WReach_r[G,S,\sigma,u]| \leq \cc_2(t,Y)k \cdot r^{t-1} \log r.
        \]
        Let $x_u \in V(T)$ be such that if $u \in S$, then $u \in P_{x_u}$, and otherwise, $x_u$ is the vertex of $T$ furthest to $s$ such that $P_{x_u}$ intersects $N_G(V(C))$, where $C$ is the component of $u$ in $G-S$.
        Let $A$ be the set of all the ancestors of $x_u$ in $T$ in distance at most $r$ from $x_u$.
        In particular, $|A| \leq r+1$.
        By \Cref{obs:wcols:trees} and \ref{item:bounded_tw_claim_adding_apices_iii-called-in-next},
        \[
            \WReach_r[G,S,\sigma,u] \subset \bigcup_{y \in A} P_y.
        \]
        Let $y \in A$.
        Since $\sigma$ extends $\sigma_y$ and $x \in T_y$, by~\ref{item:bounded_tw_claim_adding_apices_vii-called-in-next}, we have
        \[
            |\WReach_r[G,S,\sigma,u] \cap P_y| \leq \wcol_r\left(G_y, P_y\right) \leq \cc_1(t,X)k \cdot r^{t-2} \log r.
        \]
        Summarizing,
        \begin{align*}
            |\WReach_r[G,S,\sigma,u]| &= \sum_{y \in A} |\WReach_r[G,S,\sigma,u] \cap P_y| \\
            &\leq (r+1) \cdot \cc_1(t,X)k \cdot r^{t-2} \log r\\
            &\leq 2r \cdot \cc_1(t,X)k \cdot r^{t-2} \log r= \cc_2(t,Y)k \cdot r^{t-1} \log r.
        \end{align*}
    This concludes the proof of~\ref{item:claim:bounded_tw_adding_apices_wcol:iii}, and the claim follows.       
    \end{proofclaim}

    \begin{claim1}\label{claim:bounded_tw_last_claim}
        Let $Y$ be a graph with $\rtd_2(Y) \leq t-1$ and let $h,d$ be positive integers.
        There exists an integer $\cc_3(Y,h,d)$ such that 
        for every positive integer $k$,
        for every graph $G$ with $\tw(G)<k$,
        for every family $\mathcal{F}$ of connected subgraphs of $G$,
        if $G$ has no $\mathcal{F}$-rich model of $T'_{h,d}(Y)$, then there exists
        $S \subseteq V(G)$ such that
        \begin{enumerateOurAlph}
            \item $V(F) \cap S \neq \emptyset$ for every $F \in \mathcal{F}$; \label{item:bounded_tw_last_claim_i}
            \item for every component $C$ of $G-S$, $N_G(V(C))$ intersects at most $4h$ components of $G-V(C)$;\label{item:bounded_tw_last_claim_ii}
            \item $\wcol_r(G,S) \leq \cc_3(Y,h,d)k \cdot r^{t-2}\log r$ for every integer $r$ with $r \geq 2$. \label{item:bounded_tw_last_claim_iii}
        \end{enumerateOurAlph}
    \end{claim1}

    \begin{proofclaim}
        We proceed by induction on $h$.
        When $h=1$, $T'_{h,d}(Y) = K_1 \oplus (Y \sqcup Y)$ and the result follows from the previous claim applied to $Y \sqcup Y$ (note that $\rtd_2(Y \sqcup Y) = \rtd_2(Y)\leq t-1$)
        by setting $\cc_3(Y,1,d)=\cc_2(t,Y \sqcup Y)$.

        Now assume that $h > 1$ and that the result holds for $h-1$.
        Fix a copy of $Y \sqcup Y$.
        For each $y \in V(Y \sqcup Y)$ add $d$ vertices with $y$ as a unique neighbor.
        Furthermore, add $2d$ isolated vertices.
        We call the obtained graph $Z$ -- see~\Cref{fig:combining-models-T}.
        To keep things in order, we write $V(Z) = V_Y \cup V_Z$, where $V_Y$ are the vertices of $Y \sqcup Y$ in $Z$ and $V_Z$ are all the added vertices.
        Let
        \[
            \cc_3(Y,h,d) = \cc_2(t,Z) + \cc_3(Y,h-1,d).
        \]
        By~\ref{rtd2:item:degree1} and \ref{rtd2:item:components}, $\rtd_2(Z) \leq \max\{2,\rtd_2(Y)\} \leq t-1$ since $t \geq 3$.

        Let $\mathcal{F}'$ be the family of all the connected subgraphs $H$ of $G$
        such that $H$ has an $\mathcal{F}\vert_H$-rich model of $T'_{h-1,d}(Y)$.
        We claim that there is no $\mathcal{F}'$-rich model of $K_1 \oplus Z$ in $G$.
        Suppose to the contrary that such a model $\big(A_y \mid y \in V(K_1 \oplus Z)\big)$ exists.
        Let $s$ be the vertex of $K_1$ in $K_1 \oplus Z$.
        In particular, $V(K_1 \oplus Z) = \{s\} \cup V_Y \cup V_Z$.
        For every vertex $z \in V_Z$, we define its parent in $K_1 \oplus Z$ in the following way.
        Note that $z$ has at most one neighbor in $Z$.
        If $z$ has a neighbor in $Z$, then the neighbor is its parent and if $z$ is isolated in $Z$, then $s$ is its parent.
        Let $z \in V_Z$ with the parent $p_z$ and let $u_z \in A_z$ be such that there is an edge between $u_z$ and a vertex in $A_{p_z}$. 
        Since the model is $\mathcal{F}'$-rich, $G[A_z]$ contains an $\mathcal{F}\vert_{G[A_z]}$-rich model of $T'_{h-1,d}(Y)$.
        Let $H_z$ be a copy of $T_{h-1,d}(Y)$ with a root $s_z$.
        By \Cref{lemma:root_a_model_of_T'hd} applied to $Y$, $G[A_{z}]$, and $u$,
        there is an $\mathcal{F}\vert_{G[A_z]}$-rich model $\big(B_{z,x} \mid x \in V(H_{z})\big)$ of $H_{z}$ in $G[A_z]$ such that $u_z \in B_{z,s_z}$.
        In particular, there is an edge between $B_{z,s_z}$ and $A_{p_z}$ in $G$.
        Finally, we construct an $\mathcal{F}$-rich model of $T'_{h,d}(Y)$ in $G$.
        Observe that the graph obtained from $(K_1 \oplus Z)[\{s\} \cup V_Y]$ (this graph is isomorphic to $K_1 \oplus (Y \sqcup Y))$ and the disjoint union of $H_z$ for each $z \in V_Z$ by identifying $s_z \in V(H_z)$ with $p_z \in \{s\} \cup V_Y$ for each $z \in V_Z$ is isomorphic to $T'_{h,d}(Y)$.
        For each $p \in \{s\} \cup V_Y$, let $P_y$ be the set of all $z \in V_Z$ such that $p$ is the parent of $z$.
        Let
        \begin{enumerate}[label={\normalfont (\roman*)}]
            \item $D_p = A_p \cup \bigcup_{z \in P_p} B_{z,s_z}$ for every $p \in \{s\} \cup V_Y$ and
            \item $D_x = B_{z,x}$ for every $z \in V_Z$ and $x \in V(H_z \setminus \{s_{z}\})$.
        \end{enumerate}
        It follows that $\big(D_x \mid x \in \{s\} \cup V_Y \cup \bigcup_{z \in V_Z} V(H_z \setminus \{s_z\}) \big)$ is an $\mathcal{F}$-rich model of $T'_{h,d}(Y)$ in $G$.
        This is a contradiction, hence, $G$ has no $\mathcal{F}'$-rich model of $K_1 \oplus Z$.

        By \Cref{claim:bounded_tw_adding_apices_wcol}, applied to $Z$, $G$, and $\mathcal{F}'$ there exists $S_0 \subseteq V(G)$ such that
        \begin{enumerate}[label={\normalfont\ref{claim:bounded_tw_adding_apices_wcol}.(\makebox[\mywidth]{\alph*})}]
            \item $V(F) \cap S_0 \neq \emptyset$ for every $F \in \mathcal{F}'$;
            \item for every component $C$ of $G-S_0$, $N_G(V(C))$ intersects at four components of $G-V(C)$;\label{claim:bounded_tw_adding_apices_less_technical_J-final:connected}
            \item  $\wcol_r(G,S_0) \leq \cc_2(t,1,Z)k \cdot r^{t-2} \log r$ for every integer $r$ with $r \geq 2$. \label{claim:bounded_tw_adding_apices_less_technical_J-final:wcol}
        \end{enumerate}
        Let $C$ be a component of $G - S_0$.
        Since $V(F) \cap S_0 \neq \emptyset$ for every $F \in \mathcal{F}'$, $C$ has no $\mathcal{F}\vert_C$-rich model of $T'_{h-1,d}(Y)$.
        Therefore, by induction hypothesis, 
        there exists $S_C \subset V(C)$ such that
        \begin{enumerateOurAlphPrim}
            \item $V(F) \cap S_C \neq \emptyset$ for every $F \in \mathcal{F}\vert_C$; \label{bounded_tw_final_claim:item:hitting'}
            \item for every component $C'$ of $C-S$, $N_C(V(C'))$ intersects at most $4(h-1)$ components of $C-V(C')$;\label{bounded_tw_final_claim:item:connected'}
            \item $\wcol_r(C,S_C) \leq \cc_3(Y,h-1,d)k \cdot r^{t-2} \log r$ for every integer $r$ with $r \geq 2$. \label{bounded_tw_final_claim:item:wcol'}
        \end{enumerateOurAlphPrim}

        Let $\mathcal{C}$ be the family of all the components of $G-S_0$ and let 
        \[
        S = S_0 \cup \bigcup_{C \in \mathcal{C}} S_C.
        \]
        We claim that \ref{item:bounded_tw_last_claim_i}-\ref{item:bounded_tw_last_claim_iii} hold.
        Let $F \in \mathcal{F}$.
        If $V(F) \cap S_0 = \emptyset$, then $V(F) \subset V(C)$ for some component $C$ of $G - S_0$.
        In particular, $F \in \mathcal{F}\vert_C$, and thus, by~\ref{bounded_tw_final_claim:item:hitting'}, $V(F) \cap S_C \neq \emptyset$, which proves~\ref{item:bounded_tw_last_claim_i}.
        Let $C$ be a component of $G-S$, and let $C'$ be the component of $G-S_0$ intersecting $V(C)$.
        By \ref{bounded_tw_final_claim:item:connected'}, $N_C(V(C'))$ intersects at most $4(h-1)$ components of $C-V(C')$, and by \ref{claim:adding_apices_less_technical_J-final:connected},
        $N_G(V(C'))$ intersects at most fours components of $G-V(C')$.
        Hence, $N_G(V(C))$ intersects at most $4h$ components of $G-V(C)$, which yields~\ref{item:bounded_tw_last_claim_ii}.
        The following sequence of inequalities concludes the proof of~\ref{item:bounded_tw_last_claim_iii} and the claim:
        \begin{align*}
        \wcol_r(G,S) &\leq \wcol_r(G,S_0) + \wcol_r\left(G-S_0,\bigcup_{C \in \mathcal{C}} S_C\right)&&\textrm{by~\Cref{obs:wcol_union}}\\ 
        &\leq \wcol_r(G,S_0) + \max_{C \in \mathcal{C}} \wcol_r(C, S_C)&&\textrm{by~\Cref{obs:wcol_components2}}\\ 
        &\leq \cc_2(t,Z)k \cdot r^{t-2} \log r + \cc_3(Y,h-1,d)k \cdot r^{t-2} \log r &&\textrm{by~\ref{claim:bounded_tw_adding_apices_less_technical_J-final:wcol} and \ref{bounded_tw_final_claim:item:wcol'}}\\ 
        &\leq \big(\cc_2(t,Z) + \cc_3(Y,h-1,d)\big)k \cdot r^{t-1} \log r\\
        &= \cc_3(Y,h,d)k \cdot r^{t-1} \log r. 
        \end{align*}
    \end{proofclaim}

    Finally, by \Cref{lemma:Thd_universal_for_rtd}, for every for every graph $X$ with $\rtd_2(X)\leq t$, there exists a graph $Y$ with $\rtd_2(Y) \leq t-1$
    and positive integers $h,d$
    such that $X \subseteq T_{h,d}(Y) \subseteq T'_{h,d}(Y)$.
    Let $G$ be a graph with $\tw(G) < k$ and let $\mathcal{F}$ be a family of connected subgraphs of $G$.
    Suppose that $G$ has no $\mathcal{F}$-rich model of $X$.
    By \Cref{claim:bounded_tw_last_claim},
    there exists $S_0 \subseteq V(G)$ such that
    \begin{enumerate}[label={\normalfont\ref{claim:bounded_tw_last_claim}.(\makebox[\mywidth]{\alph*})}]
            \item $V(F) \cap S_0 \neq \emptyset$ for every $F \in \mathcal{F}$; \label{final-final-i}
            \item for every component $C$ of $G-S_0$, $N_G(V(C))$ intersects at most $4h$ components of $G-V(C)$;\label{final-final-ii}
            \item $\wcol_r(G,S) \leq \cc_3(Y,h,d)k \cdot r^{t-2}\log r$ for every integer $r$ with $r \geq 2$. \label{final-final-iii}
    \end{enumerate}
    Let $\mathcal{C}$ be the family of components of $G - S_0$.
    Consider a component $C \in \mathcal{C}$.
    By~\ref{final-final-ii}, $N_G(V(C))$ intersects at most $4h$ components of $G - V(C)$.
    Let $C'_1, \dots, C'_a$ be the components of $G-V(C)$. For every $i \in [a-1]$, let $Q_i$ be a shortest $V(C'_i)$-$V(C'_{i+1})$ path in $G$.
    Now, let $\mathcal{Q}_C$ be the family $\{Q_1, \dots, Q_{a-1}\}$.
    Note that for every component $C'$ of $C-\bigcup_{Q \in \mathcal{Q}_C}V(Q)$, $N_G(V(C'))$ intersects at most one component of $G-V(C')$.

    Let
    \[
        S = S_0 \cup \bigcup_{C \in \mathcal{C}} \bigcup_{Q \in \mathcal{Q}_C} V(Q).
    \]
    Item~\ref{thm:bounded_tw_main:item:hit} follows from~\ref{final-final-i}.
    Let $C'$ be a component of $G-S$, and let $C$ be the component of $G-S_0$ intersecting $V(C')$.
    Then $N_G(V(C')) \subseteq N_G(V(C)) \cup \bigcup_{Q \in \mathcal{Q}_C} V(Q)$, and since $\bigcup_{Q \in \mathcal{Q}_C} V(Q)$ 
    induces a connected subgraph of $C$ and contains a neighbor
    of every component of $G-V(C)$ having a neighbor in $V(C)$, we deduce that
    $N_G(V(C'))$ intersects at most one component of $G-V(C')$.
    This proves that \ref{thm:bounded_tw_main:item:connected} holds.
    And finally,~\ref{thm:bounded_tw_main:item:wcol} is true once we set $\cc(t,X) = \cc_3(Y,h,d) + (4h-1) \cdot 3$ since
    \begin{align*}
        \wcol_r(G,S) &\leq \wcol_r(G,S_0) + \wcol_r\left(G-S_0,\bigcup_{C \in \mathcal{C}} \bigcup_{Q \in \mathcal{Q}_C} V(Q) \right)&&\textrm{by~\Cref{obs:wcol_union}}\\ 
        &\leq \wcol_r(G,S_0) + \max_{C \in \mathcal{C}} \wcol_r\left(C, \bigcup_{Q \in \mathcal{Q}_C} V(Q)\right)&&\textrm{by~\Cref{obs:wcol_components2}}\\
        &\leq \wcol_r(G,S_0) + \max_{C \in \mathcal{C}} |\mathcal{Q}_C|(2r+1)&&\textrm{by~\Cref{obs:geodesics}}\\
        &\leq \cc_3(Y,h,d)k \cdot r^{t-2}\log r + (4h-1)(2r+1)&&\textrm{by~\ref{final-final-iii}}\\ 
        &\leq \left(\cc_3(Y,h,d) + (4h-1) \cdot 3\right)k \cdot r^{t-2}\log r. &&\qedhere
    \end{align*}      
\end{proof}

\bibliographystyle{plain}
\bibliography{biblio}

\providecommand{\noopsort}[1]{}
\begin{thebibliography}{10}

\bibitem{CS93}
Guantao~T. Chen and Richard~H. Schelp.
\newblock Graphs with linearly bounded {R}amsey numbers.
\newblock {\em Journal of Combinatorial Theory, Series B}, 57(1):138--149, 1993.

\bibitem{DHHJLMMRW24}
Vida Dujmovic, Robert Hickingbotham, Jędrzej Hodor, Gwena{\"{e}}l Joret, Hoang La, Piotr Micek, Pat Morin, Cl{\'{e}}ment Rambaud, and David~R. Wood.
\newblock The grid-minor theorem revisited.
\newblock In David~P. Woodruff, editor, {\em Proceedings of the 2024 {ACM-SIAM} Symposium on Discrete Algorithms, {SODA} 2024, Alexandria, VA, USA, January 7-10, 2024}, pages 1241--1245. {SIAM}, 2024.
\newblock \href{https://arxiv.org/abs/2307.02816}{arXiv:2307.02816}.

\bibitem{Dujmovi2023}
Vida Dujmović, Robert Hickingbotham, Gwenaël Joret, Piotr Micek, Pat Morin, and David~R. Wood.
\newblock The excluded tree minor theorem revisited.
\newblock {\em Combinatorics, Probability and Computing}, 33(1):85–90, 2023.
\newblock \href{https://arxiv.org/abs/2303.14970}{arXiv:2303.14970}.

\bibitem{Dvorak13}
Zdeněk Dvořák.
\newblock Constant-factor approximation of the domination number in sparse graphs.
\newblock {\em European Journal of Combinatorics}, 34(5):833--840, 2013.
\newblock \href{https://arxiv.org/abs/1110.5190}{arXiv:1110.5190}.

\bibitem{EGKKPRS17}
Kord Eickmeyer, Archontia~C. Giannopoulou, Stephan Kreutzer, O-joung Kwon, Michal Pilipczuk, Roman Rabinovich, and Sebastian Siebertz.
\newblock {Neighborhood Complexity and Kernelization for Nowhere Dense Classes of Graphs}.
\newblock In Ioannis Chatzigiannakis, Piotr Indyk, Fabian Kuhn, and Anca Muscholl, editors, {\em 44th International Colloquium on Automata, Languages, and Programming (ICALP 2017)}, volume~80 of {\em Leibniz International Proceedings in Informatics (LIPIcs)}, pages 63:1--63:14, Dagstuhl, Germany, 2017. Schloss Dagstuhl -- Leibniz-Zentrum f{\"u}r Informatik.
\newblock \href{https://arxiv.org/abs/1612.08197}{arXiv:1612.08197}.

\bibitem{FN06}
Pierre Fraigniaud and Nicolas Nisse.
\newblock Connected treewidth and connected graph searching.
\newblock In Jos{\'{e}}~R. Correa, Alejandro Hevia, and Marcos~A. Kiwi, editors, {\em Proc. 7th Latin American Symposium on Theoretical Informatics, {LATIN} 2006}, volume 3887 of {\em Lecture Notes in Comput. Sci.}, pages 479--490. Springer, 2006.

\bibitem{GJNW23}
Carla Groenland, Gwena\"{e}l Joret, Wojciech Nadara, and Bartosz Walczak.
\newblock Approximating pathwidth for graphs of small treewidth.
\newblock {\em ACM Transactions on Algorithms}, 19(2), 2023.
\newblock \href{https://arxiv.org/abs/2008.00779}{arXiv:2008.00779}.

\bibitem{Grohe15}
Martin Grohe, Stephan Kreutzer, Roman Rabinovich, Sebastian Siebertz, and Konstantinos Stavropoulos.
\newblock Coloring and covering nowhere dense graphs.
\newblock {\em SIAM Journal on Discrete Mathematics}, 32(4):2467--2481, 2018.
\newblock \href{https://arxiv.org/abs/1602.05926}{arXiv:1602.05926}.

\bibitem{GroheKreutzerSiebertz17}
Martin Grohe, Stephan Kreutzer, and Sebastian Siebertz.
\newblock Deciding first-order properties of nowhere dense graphs.
\newblock {\em Journal of the ACM}, 64(3), 2017.
\newblock \href{https://arxiv.org/abs/1311.3899}{arXiv:1311.3899}.

\bibitem{vandenHeuvel2018}
Jan van~den Heuvel and David~R. Wood.
\newblock Improper colourings inspired by {H}adwiger’s conjecture.
\newblock {\em Journal of the London Mathematical Society}, 98(1):129–148, 2018.
\newblock \href{https://arxiv.org/abs/1704.06536}{arXiv:1704.06536}.

\bibitem{vdHetal17}
Jan {\noopsort{Heuvel}}{van den Heuvel}, Patrice {Ossona de Mendez}, Daniel Quiroz, Roman Rabinovich, and Sebastian Siebertz.
\newblock On the generalised colouring numbers of graphs that exclude a fixed minor.
\newblock {\em European Journal of Combinatorics}, 66:129--144, 2017.
\newblock \href{https://arxiv.org/abs/1602.09052}{arXiv:1602.09052}.

\bibitem{HJMSW22}
Tony Huynh, Gwena{\"{e}}l Joret, Piotr Micek, Michal~T. Seweryn, and Paul Wollan.
\newblock Excluding a ladder.
\newblock {\em Combinatorica}, 42(3):405--432, 2022.
\newblock \href{https://arxiv.org/abs/2002.00496}{arXiv:2002.00496}.

\bibitem{JM22}
Gwena{\"{e}}l Joret and Piotr Micek.
\newblock Improved bounds for weak coloring numbers.
\newblock {\em The Electronic Journal of Combinatorics}, 29(1), 2022.
\newblock \href{https://arxiv.org/abs/2102.10061}{arXiv:2102.10061}.

\bibitem{KY03}
Hal~A. Kierstead and Daqing Yang.
\newblock Orderings on graphs and game coloring number.
\newblock {\em Order}, 20(3):255--264, 2003.

\bibitem{sparsity}
Jaroslav Ne\v{s}et\v{r}il and Patrice {Ossona de Mendez}.
\newblock {\em Sparsity --- {G}raphs, {S}tructures, and {A}lgorithms}, volume~28 of {\em Algorithms and {C}ombinatorics}.
\newblock Springer, 2012.

\bibitem{notes}
Marcin Pilipczuk, Micha\l{} Pilipczuk, and Sebastian Siebertz.
\newblock Lecture notes for the course ``{S}parsity'' given at {F}aculty of {M}athematics, {I}nformatics, and {M}echanics of the {U}niversity of {W}arsaw, Winter semesters 2017/18 and 2019/20.
\newblock \url{https://www.mimuw.edu.pl/~mp248287/sparsity2}.

\bibitem{PS19}
Michal Pilipczuk and Sebastian Siebertz.
\newblock Polynomial bounds for centered colorings on proper minor-closed graph classes.
\newblock {\em {Journal of Combinatorial Theory, Series B}}, 151:111--147, 2021.
\newblock \href{https://arxiv.org/abs/1807.03683}{arXiv:1807.03683}.

\bibitem{ReidlSullivan23}
Felix Reidl and Blair~D. Sullivan.
\newblock A color-avoiding approach to subgraph counting in bounded expansion classes.
\newblock {\em Algorithmica}, 85(8):2318–2347, 2023.
\newblock \href{https://arxiv.org/abs/2001.05236}{arXiv:2001.05236}.

\bibitem{GM1}
Neil Robertson and Paul Seymour.
\newblock {Graph minors. I. Excluding a forest}.
\newblock {\em Journal of Combinatorial Theory, Series B}, 35(1):39–61, 1983.

\bibitem{GM5}
Neil Robertson and Paul Seymour.
\newblock Graph minors. {V}. {Excluding a planar graph}.
\newblock {\em Journal of Combinatorial Theory, Series B}, 41(1):92--114, 1986.

\bibitem{Zhu09}
Xuding {Zhu}.
\newblock {Colouring graphs with bounded generalized colouring number.}
\newblock {\em {Discrete Mathematics}}, 309(18):5562--5568, 2009.

\end{thebibliography}
\end{document}